\documentclass[fontsize=11pt]{article}

\usepackage[dvipsnames]{xcolor}

\usepackage[utf8]{inputenc}
\usepackage[T1]{fontenc}
\usepackage{amsmath}
\usepackage{amsfonts}
\usepackage{amssymb}
\usepackage[total={6.5in, 9in}, heightrounded]{geometry}

\usepackage{graphicx}
\graphicspath{{graphics/}}

\usepackage{scalerel}
\usepackage{tikz}
\usepackage{exscale}
\usepackage{wallpaper}
\usepackage{tcolorbox}
\usepackage{relsize}
\usepackage{float}
\usepackage{stmaryrd}
\usepackage{MnSymbol}
\usepackage{mdframed}
\usepackage{titlesec}
\usepackage{tabu}
\usepackage{caption}
\usepackage{enumitem}
\usepackage{stackengine}
\usepackage{blkarray}
\usepackage{xifthen}
\usepackage{ytableau}
\usepackage{pdfpages}
\usepackage{multicol}
\usepackage{array}
\usepackage{mathtools}
\usepackage{soul}
\usepackage{hyperref}

\hypersetup
{
	colorlinks = true,
	allcolors = OliveGreen
}

\titleformat{\section}
{\normalfont \Large \bfseries \centering}{\arabic{section}. }{0pt}{}

\setenumerate[0]{label=\alph*)}

\usepackage{amsthm}

\newcommand{\bra}[1]{\left< #1 \right| \,}
\newcommand{\ket}[1]{\,\left| #1 \right>}

\newcommand{\leftbra}[1]{\left< #1 \left| \,}
\newcommand{\rightket}[1]{\,\right| #1 \right>}

\newtheoremstyle{mytheoremstyle} 
  {8pt}                    
  {0}                      
  {}                       
  {}                       
  {\bfseries}              
  {.}                      
  {.5em}                   
  {}                       

\theoremstyle{mytheoremstyle}
\newtheorem{Thm}{Theorem}[section]
\newtheorem{Def}[Thm]{Definition}
\newtheorem{Prop}[Thm]{Proposition}
\newtheorem{Lem}[Thm]{Lemma}

\newtheorem{Ex}[Thm]{Example}

\makeatletter
\newenvironment{Proof}[1][\proofname]{\par
  \pushQED{\qed}%
  \normalfont \partopsep=\z@skip \topsep=\z@skip
  \trivlist
  \item[\hskip\labelsep
        \itshape
    #1\@addpunct{.}]\ignorespaces
}{%
  \popQED\endtrivlist\@endpefalse
}
\makeatother

\newcommand{\tableaubox}[4][]{
  \draw[#1] ({#3/2},{-#2/2}) -- ({#3/2},{-(#2/2+0.5)}) -- ({#3/2+0.5},{-(#2/2+0.5)}) -- ({#3/2+0.5},{-#2/2}) -- cycle;
  \node at ({#3/2+0.25},{-(#2/2+0.25)}) {#4};
}

\newcommand{\tableaulabel}[4][]{
  \node[#1] at ({#3/2+0.25},{-(#2/2+0.25)}) {#4};
}

\newcommand{\tableauline}[5][]{
  \draw[#1] ({(#3+0.5)/2},{-(#2+0.5)/2}) -- ({(#5+0.5)/2},{-(#4+0.5)/2});
}

\newcommand{\tableaudraw}[2][]{
  \draw[#1][cm={0, -0.5, 0.5, 0, (0.25, -0.25)}] #2;
}

\newcommand{\tableauedgelabelleft}[4][]{
  \node[#1,align=right,text width=.45cm,minimum width=.45cm] at ({(#3-1)/2+0.25},{-(#2/2+0.25)}) {\tiny #4};
}

\newcommand{\tableauedgelabeltop}[4][]{
  \node[#1] at ({#3/2+0.25},{-((#2-1+.3)/2+0.25)}) {\tiny #4};
}

\counterwithin{figure}{section}

\newcommand*{\figref}[2][]{%
  \hyperref[{#2}]{%
    Figure~\ref*{#2}%
    \ifx\\#1\\%
    \else
      \,#1%
    \fi
  }%
}

\renewcommand*{\eqref}[2][]{%
  \hyperref[{#2}]{%
    Equation~(\ref*{#2}%
    \ifx\\#1\\%
    \else
      \,#1%
    \fi
  )}%
}

\newcommand*{\defref}[2][]{%
  \hyperref[{#2}]{%
    Definition~\ref*{#2}%
    \ifx\\#1\\%
    \else
      \,#1%
    \fi
  }%
}

\newcommand*{\lemref}[2][]{%
  \hyperref[{#2}]{%
    Lemma~\ref*{#2}%
    \ifx\\#1\\%
    \else
      \,#1%
    \fi
  }%
}

\newcommand*{\propref}[2][]{%
  \hyperref[{#2}]{%
    Proposition~\ref*{#2}%
    \ifx\\#1\\%
    \else
      \,#1%
    \fi
  }%
}

\newcommand*{\thmref}[2][]{%
  \hyperref[{#2}]{%
    Theorem~\ref*{#2}%
    \ifx\\#1\\%
    \else
      \,#1%
    \fi
  }%
}

\newcommand*{\corref}[2][]{%
  \hyperref[{#2}]{%
    Corollary~\ref*{#2}%
    \ifx\\#1\\%
    \else
      \,#1%
    \fi
  }%
}

\newcommand*{\conjref}[2][]{%
  \hyperref[{#2}]{%
    Conjecture~\ref*{#2}%
    \ifx\\#1\\%
    \else
      \,#1%
    \fi
  }%
}

\newcommand*{\secref}[2][]{%
  \hyperref[{#2}]{%
    Section~\ref*{#2}%
    \ifx\\#1\\%
    \else
      \,#1%
    \fi
  }%
}

\newcommand*{\tblref}[2][]{%
  \hyperref[{#2}]{%
    Table~\ref*{#2}%
    \ifx\\#1\\%
    \else
      \,#1%
    \fi
  }%
}

\DeclareGraphicsExtensions{.pdf}

\usepackage[style=alphabetic,backend=bibtex,sorting=ynt]{biblatex}

\title{Bijectivizing the PT--DT Correspondence}
\author{Cruz Godar and Benjamin Young}

\begin{document}
\maketitle

\begin{abstract}
  \noindent Pandharipande--Thomas theory and Donaldson--Thomas theory (PT and DT) are two branches of enumerative geometry in which particular generating functions arise that count plane-partition-like objects. That these generating functions differ only by a factor of MacMahon's function was proven recursively in \cite{jenneWebbYoung} using the double dimer model. We bijectivize two special cases of the result by formulating these generating functions using vertex operators and applying a particular type of local involution known as a \textit{toggle}, first introduced in the form we use in \cite{pak}.
\end{abstract}

\section{Introduction}

Plane partitions, along with their many generalizations and their cousins the standard and semistandard Young tableaux, are frequently-studied objects in combinatorics (see, for instance, \cite{bressoud}). As is often the case with combinatorial objects of interest, they arise in other areas of mathematics; relevant to our work is a particular use in computing the equivariant Calabi--Yau topological vertex in Pandharipande--Thomas theory and Donaldson--Thomas theory (PT and DT). This generating function is a local contribution to the generating function of DT (resp.\ PT) invariants for a Calabi--Yau 3-fold with a torus action. One may use Atiyah--Bott localization to reduce the computation to one on the fixed loci of the torus action; the vertex represents the contribution of one fixed point. We refer the interested reader to \cite{MNOP} for further details. In \cite{PT}, the authors conjecture that two generating functions that count DT and PT objects, which are generalizations of standard and reverse plane partitions, respectively, are equal up to a factor of MacMahon's function $M(q)$ (i.e.\ the generating function for standard plane partitions). The fully general conjecture was proven in \cite{jenneWebbYoung} using the double-dimer model; however, this last proof is strikingly involved, and no combinatorial proof has been given for the general case or any special case. We give a combinatorial proof for two special cases, with the goal of eventually extending our methods to the fully general case.

The generating function for reverse plane partitions was first derived by Stanley \cite{stanleyThesis}, and later bijectivized by Hillman and Grassl \cite{hillmanGrassl}. Their map places reverse plane partitions in bijection with tableaux of the same shape containing nonnegative integers, and the lack of the plane partition inequalities greatly simplifies the process of working with them. A second bijection was introduced by Pak \cite{pak} (and later independently by Sulzgruber \cite{sulzgruber}), using local operations on the diagonals of reverse plane partitions; these local moves were later independently introduced in \cite{hopkinsNotes}, and they are the main tool we use.

In \secref{sec:planePartitions}, we discuss plane partitions and a pair of vertex operators that allow us to build MacMahon's generating function $M(q)$ one diagonal at a time. In \secref{sec:toggles}, we examine the toggle operation on diagonals and use it to introduce bijective proofs of the commutation relations of the vertex operators. In \secref{sec:oneLegObjects}, we define a special case of the objects counted by the PT and DT generating functions (so-called \textit{one-leg} objects), and we give a bijective proof that the PT--DT correspondence holds for them. In \secref{sec:twoLegObjects}, we do the same for a more general case (\textit{two-leg} objects). Finally, in \secref{sec:futureDirections}, we define the fully general (i.e.\ \textit{three-leg}) case of the objects, give an explanation as to why the correspondence is fundamentally more difficult in this case, and discuss how we hope to generalize our methods in the future.

\textbf{Acknowledgments:} The authors thank Tatyana Benko and Ava Bamforth for helpful conversations regarding the labeling conditions of 3-leg RPPs.

\section{Plane Partitions} \label{sec:planePartitions}

We begin with straightforward integer partitions, which give rise to Young diagrams and thereby plane partitions. We review several definitions; for further details, we refer the interested reader to \cite{stanley}, \cite{bressoud}, or \cite{andrews}.

\begin{Def}
  A \textbf{partition} $\lambda$ is a sequence of nonnegative integers $\lambda = \left( \lambda_1, \lambda_2, \ldots \right)$ such that $\lambda_i \geq \lambda_{i + 1}$ for all $i \in \mathbb{N}$\footnote{We use $\mathbb{N}$ to denote the set $\left\{ 1, 2, ... \right\}$ of natural numbers, and $\mathbb{N}_{\geq 0}$ for $\left\{ 0, 1, 2, ... \right\}$.} and only finitely many $\lambda_i$ are nonzero. Since every partition ends in zeros, we typically write only the nonzero entries. The \textbf{weight} of $\lambda$ is $ \left| \lambda \right| = \sum_{i = 1}^\infty \lambda_i$, and the \textbf{Young diagram} corresponding to $\lambda$ is the subset $\left\{ (i, j) \in \mathbb{N}^2 \mid 1 \leq j \leq \lambda_i \right\}$, which we represent as a set of top-left-justified squares whose row lengths weakly decrease. Finally, the \textbf{conjugate partition} to $\lambda$ is the partition $\lambda'$ given by the column lengths of the Young diagram corresponding to $\lambda$.
\end{Def}

The following definition is slightly unusual, but it is a special case of a more natural generalization of plane partitions that we introduce later.

\begin{Def}
  Let $\lambda \subset \mathbb{N}^2$ be a Young diagram. A \textbf{Young diagram asymptotic to $\lambda$} is the collection of boxes $\mathbb{N}^2 \setminus \lambda$.
\end{Def}

While this is a slight abuse of notation, since Young diagrams asymptotic to a partition $\lambda$ are not in fact Young diagrams, it is a convenient definition for objects we will define in \secref{sec:oneLegObjects}. We also recall the usual notion of a hook of a plane partition, which is a right-angled collection of cells extending toward the boundary of a Young diagram.

\begin{Def}
  Let $\lambda \subset \mathbb{N}^2$ be a Young diagram and let $(i, j) \in \lambda$. The \textbf{arm} and \textbf{leg} of $(i, j)$ are
  \begin{align*}
    \operatorname{arm}(i, j) &= \left\{ (i, j') \mid j < j' \leq \lambda_i \right\}\\
    \operatorname{leg}(i, j) &= \left\{ (i', j) \mid i < i' \leq \lambda'_j \right\}.
  \end{align*}
  If $(i, j) \in \mathbb{N}^2 \setminus \lambda$, then
  \begin{align*}
    \operatorname{arm}(i, j) &= \left\{ (i, j') \mid \lambda_i < j' < j \right\}\\
    \operatorname{leg}(i, j) &= \left\{ (i', j) \mid \lambda'_j < i' < i \right\}.
  \end{align*}
  The \textbf{hook} of $(i, j)$ is
  $$
    \operatorname{hook}(i, j) = \operatorname{arm}(i, j) \cup \operatorname{leg}(i, j) \cup \left\{ (i, j) \right\},
  $$
  and the \textbf{hook length} is $h(i, j) = \left| \operatorname{hook}(i, j) \right|$. We often refer to \textbf{$n$-hooks} of a diagram, which are hooks with length $n$, and the \textbf{pivot} of a hook, which is the corner $(i, j)$ from which it is defined.
\end{Def}

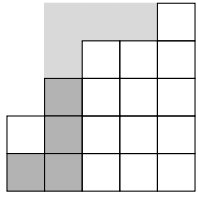
\begin{figure}
  \centering
  \begin{tikzpicture}
    
      \begin{scope}[shift={(0,0)}]
        \tableaulabel{0}{0}{}
\tableaubox[fill=gray!30,gray!30]{0}{1}{}
\tableaubox[fill=gray!30,gray!30]{0}{2}{}
\tableaubox[fill=gray!30,gray!30]{0}{3}{}
\tableaubox{0}{4}{}
\tableaulabel{1}{0}{}
\tableaubox[fill=gray!30,gray!30]{1}{1}{}
\tableaubox{1}{2}{}
\tableaubox{1}{3}{}
\tableaubox{1}{4}{}
\tableaulabel{2}{0}{}
\tableaubox[fill=gray!60]{2}{1}{}
\tableaubox{2}{2}{}
\tableaubox{2}{3}{}
\tableaubox{2}{4}{}
\tableaubox{3}{0}{}
\tableaubox[fill=gray!60]{3}{1}{}
\tableaubox{3}{2}{}
\tableaubox{3}{3}{}
\tableaubox{3}{4}{}
\tableaubox[fill=gray!60]{4}{0}{}
\tableaubox[fill=gray!60]{4}{1}{}
\tableaubox{4}{2}{}
\tableaubox{4}{3}{}
\tableaubox{4}{4}{}

      \end{scope}
    
  \end{tikzpicture} \caption{A $4$-hook with pivot $(5, 2)$ in an asymptotic Young diagram $\mathbb{N}^2 \setminus \lambda$ (dark gray), and a $4$-hook with pivot $(1, 2)$ in the Young diagram $\lambda$ (light gray).}
  \label{fig:hookExample}
\end{figure}

With Young diagrams, we can express partitions as both one-dimensional lists of weakly decreasing nonnegative integers and as two-dimensional sets of boxes; our primary object of study generalizes partitions by increasing both of these dimensions by one.

\begin{figure}
  \begin{tabu}{X[m,c]X[m,c]}
    \begin{tikzpicture}
      
      \begin{scope}[shift={(0,0)}]
        \tableaubox{0}{0}{5}
\tableaubox{0}{1}{4}
\tableaubox{0}{2}{3}
\tableaubox{0}{3}{3}
\tableaubox{1}{0}{4}
\tableaubox{1}{1}{4}
\tableaubox{1}{2}{2}
\tableaubox{2}{0}{2}
\tableaubox{2}{1}{1}
\tableaubox{3}{0}{2}
\tableaubox{3}{1}{1}

      \end{scope}
    
    \end{tikzpicture} & \includegraphics[height=1.5in]{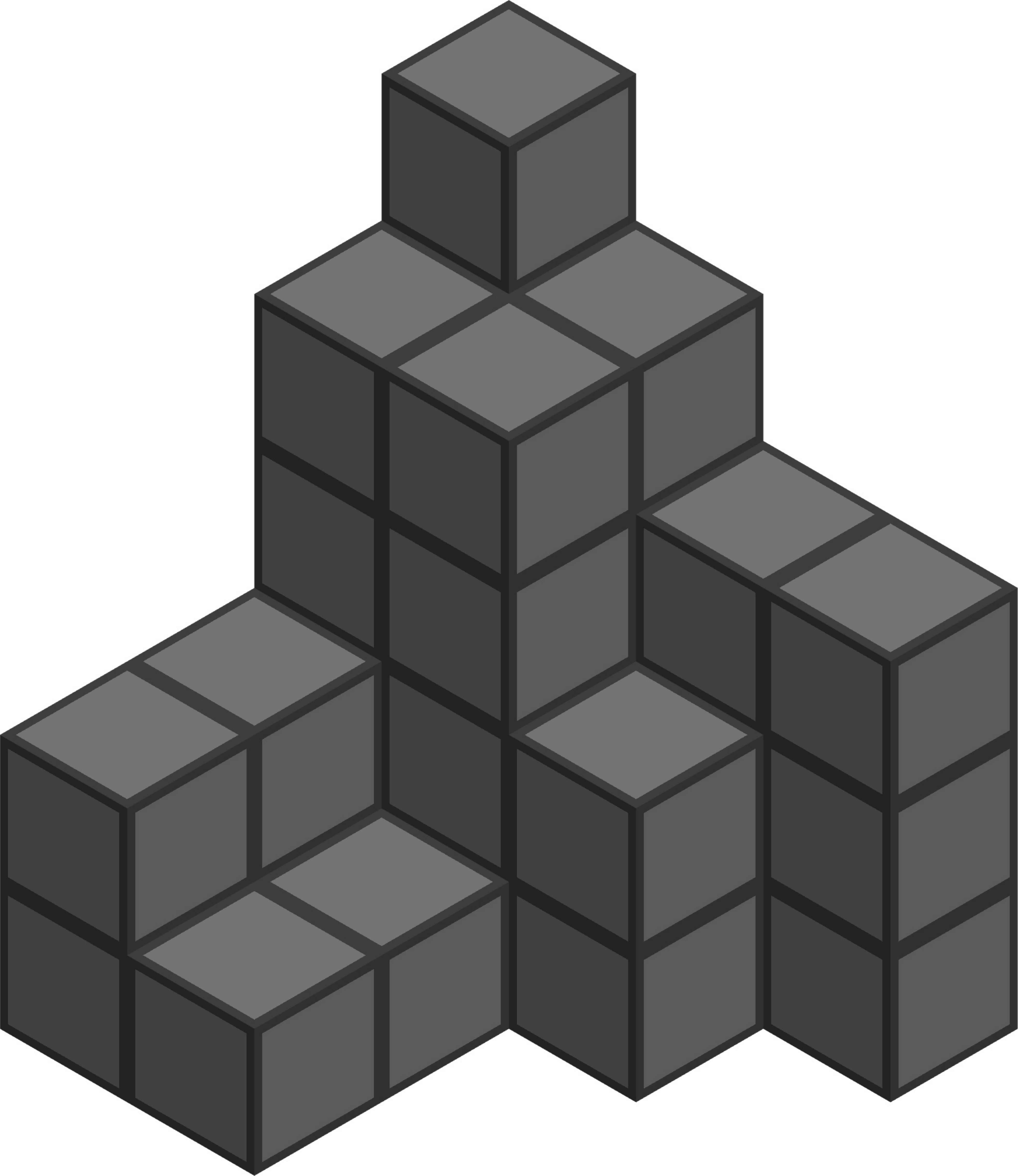}
  \end{tabu} \caption{A plane partition of weight 31, visualized both as a grid of numbers and a stack of 31 blocks.}
  \label{fig:planePartition}
\end{figure}

\begin{Def}
  A \textbf{plane partition} is a function $\pi : \mathbb{N}^2 \to \mathbb{N}_{\geq 0}$ such that $\pi(i, j) \geq \pi(i + 1, j)$ and $\pi(i, j) \geq \pi(i, j + 1)$ for all $i, j \in \mathbb{N}$ and only finitely many $\pi(i, j)$ are nonzero. The \textbf{weight} of $\pi$ is $\left| \pi \right| = \sum_{i, j \in \mathbb{N}} \pi(i, j)$. Analogous to a Young diagram, we can associate a plane partition $\pi$ with the three-dimensional stack of blocks $\left\{ (i, j, k) \in \mathbb{N}^3 \mid 1 \leq k \leq \pi(i, j) \right\}$, as in \figref{fig:planePartition}.
\end{Def}

The rows and columns of a plane partition $\pi$ are themselves interrelated partitions, but the \textit{diagonals} of $\pi$ turn out to be a more fruitful source of information. We use established techniques and terms in the remainder of this section: the definitions and proofs from here to \secref{sec:toggles} are due to \cite{okounkovReshetikhin, okounkovReshetikhin2, randomPartitions1, randomPartitions2, okounkovReshetikhinVafa}. We first define a partial order on partitions that governs when two partitions can sit next to one another in a plane partition.

\begin{Def}
  Let $\lambda$ and $\mu$ be partitions. We say $\lambda$ \textbf{interlaces} $\mu$, written $\lambda \succ \mu$, if
  $$
    \lambda_1 \geq \mu_1 \geq \lambda_2 \geq \mu_2 \geq \lambda_3 \cdots.
  $$
  In other words, $\lambda \succ \mu$ if and only if $\lambda$ and $\mu$ satisfy the plane partition inequalities when written side-by-side as diagonal slices with $\lambda$ coming first, as in \figref{fig:interlacing}.
\end{Def}

\begin{figure}
  \centering
  \begin{tikzpicture}
    
      \begin{scope}[shift={(0,0)}]
        \tableaubox{0}{0}{5}
\tableaubox{0}{1}{3}
\tableaubox{0}{2}{}
\tableaubox{0}{3}{}
\tableaubox{0}{4}{}
\tableaubox{1}{0}{}
\tableaubox{1}{1}{3}
\tableaubox{1}{2}{2}
\tableaubox{1}{3}{}
\tableaubox{1}{4}{}
\tableaubox{2}{0}{}
\tableaubox{2}{1}{}
\tableaubox{2}{2}{1}
\tableaubox{2}{3}{1}
\tableaubox{2}{4}{}
\tableaubox{3}{0}{}
\tableaubox{3}{1}{}
\tableaubox{3}{2}{}
\tableaubox{3}{3}{1}
\tableaubox{3}{4}{0}
\tableaubox{4}{0}{}
\tableaubox{4}{1}{}
\tableaubox{4}{2}{}
\tableaubox{4}{3}{}
\tableaubox{4}{4}{0}
\tableaulabel{2}{5}{$\cdots$}
\tableaulabel{5}{2}{$\vdots$}

      \end{scope}
    
  \end{tikzpicture} \caption{By placing them as diagonals in a plane partition, we see that $(5, 3, 1, 1) \succ (3, 2, 1)$.}
  \label{fig:interlacing}
\end{figure}

Given a partition $\mu$, the set of partitions $\lambda$ with $\lambda \succ \mu$ is a poset product of intervals in $\mathbb{N}_{\geq 0}$: each $\lambda_i$ is bounded below by $\mu_i$ and above by $\mu_{i - 1}$ (or unbounded for $i = 1$), and critically is independent from every other $\lambda_i$. It is largely for this reason that decomposing a plane partition into diagonals is more useful than into rows or columns, and it enables us to express the generating function for plane partitions in terms of operators on the diagonals.

\begin{Def}
  We define a formal $\mathbb{Q}$-vector space $\Lambda$ whose basis consists of vectors $\ket{\lambda}$ for partitions $\lambda$. We also denote corresponding elements of the dual basis by $\bra{\lambda}$.

  On this vector space, we define a weighing operator $Q(q) : \Lambda \to \Lambda$, given by $Q(q)\ket{\lambda} = q^{|\lambda|}\ket{\lambda}$. We also define vertex operators $\Gamma_+$ and $\Gamma_-$, each of which accepts a formal variable as a parameter: $\Gamma_\pm(q) : \Lambda \to \Lambda$ is given by
  \begin{align*}
    \Gamma_+(q)\ket{\lambda} &= \sum_{\mu \succ \lambda}q^{|\mu|-|\lambda|}\ket{\mu}\\
    \Gamma_-(q)\ket{\lambda} &= \sum_{\mu \prec \lambda}q^{|\lambda|-|\mu|}\ket{\mu}.
  \end{align*}
  The convention of which operator is denoted $\Gamma_+$ and which is denoted $\Gamma_-$ differs between \cite{okounkovReshetikhin} and the other cited papers in which these operators appear; our convention matches \cite{okounkovReshetikhin}. We think of the subscript $+$ as indicating that $\Gamma_+$ produces larger partitions, while $\Gamma_-$ produces smaller ones.

  Note that the definition of $\Gamma_+$ contains an infinite sum, which may fail to converge even in the ring of formal power series; for example, $\bra{\emptyset} \Gamma_-(1) \Gamma_+(1) \ket{\emptyset}$ is not defined. In practice, we avoid these issues by using the weighing operator $Q$ to produce sums with a formal variable $q$ that do converge. We refer the interested reader to \cite[p. 11--13]{formalPowerSeries} for further background.

  A product of $\Gamma$ operators whose arguments are all of the form $q^i$ for some $i$ therefore gives a generating function for plane-partition-like objects, since one partition interlacing another is equivalent to the two being able to sit next to one another in a plane partition. The shape of the objects counted is determined by the order of $\Gamma_+$ operators and $\Gamma_-$ operators, and its weight is determined by the arguments of the operators; this weight may or may not be equal to the sum of the entries in the object.
  \label{def:gammaOperators}
\end{Def}

These operators are defined and derived in appendix B of \cite{okounkovReshetikhin} (and ultimately from \cite{kac}) where they are used primarily to compute Schur functions, and they are also used to great effect in \cite{vertexOperators} to compute the generating functions of various objects. We first demonstrate their use in \cite{vertexOperators} to express MacMahon's function $M(q)$, i.e.\ the generating function for plane partitions.

Let $N \in \mathbb{N}$; then the generating function for plane partitions whose nonzero entries are contained in the $N \times N$ square Young diagram is
\begin{align}
  \leftbra{\emptyset} \left( \prod_{i = 1}^N Q \Gamma_-(1) \right) Q \left( \prod_{i = 1}^N \Gamma_+(1) Q \right) \rightket{\emptyset}.
  \label{eq:macMahonsFunction}
\end{align}
The presence of the weighing operators combined with the lack of variables in the interlacing operators seems to suggest that there is room for improvement in this formula, and indeed the two types of operator have a simple commutation relation that will lead to just that.

\begin{Lem}
  Let $\lambda$ and $\mu$ be partitions. Then
  \begin{align*}
    Q(q) \Gamma_+(a) &= \Gamma_+(qa) Q(q)\\
    \Gamma_-(a) Q(q) &= Q(q) \Gamma_-(qa).
  \end{align*}
\end{Lem}

\begin{Proof}
  Fix partitions $\lambda$ and $\mu$ with $\mu \succ \lambda$; then
  \begin{align*}
    \bra{\mu} Q(q) \Gamma_+(a) \ket{\lambda} &= q^{|\mu|}a^{|\mu| - |\lambda|}\\
    &= (qa)^{|\mu| - |\lambda|}q^{|\lambda|}\\
    &= \bra{\mu} \Gamma_+(qa) Q(q) \ket{\lambda},
  \end{align*}
  and similarly for the $\Gamma_-$ relation.
\end{Proof}

We are now prepared to construct a simpler expression for $M(q)$. We commute each $Q$ outward from the center in \eqref{eq:macMahonsFunction}, splitting the middle $Q$ into $Q^{1/2}Q^{1/2}$ to preserve symmetry. The $Q$ operators are annihilated at both the $\ket{\emptyset}$ and the $\bra{\emptyset}$, and so \eqref{eq:macMahonsFunction} reduces to
$$
  \leftbra{\emptyset} \Gamma_-\left( q^{\frac{2N - 1}{2}} \right) \cdots \Gamma_-\left( q^{\frac{3}{2}} \right) \Gamma_-\left( q^{\frac{1}{2}} \right) \Gamma_+\left( q^{\frac{1}{2}} \right) \Gamma_+\left( q^{\frac{3}{2}} \right) \cdots \Gamma_+\left( q^{\frac{2N - 1}{2}} \right) \rightket{\emptyset}.
$$
Since any plane partition can have only finitely many nonzero entries, taking the limit as $N \to \infty$ produces the following expression for $M(q)$:
\begin{align}
  M(q) &= \leftbra{\emptyset} \cdots \Gamma_-\left( q^{\frac{3}{2}} \right) \Gamma_-\left( q^{\frac{1}{2}} \right) \Gamma_+\left( q^{\frac{1}{2}} \right) \Gamma_+\left( q^{\frac{3}{2}} \right) \cdots \rightket{\emptyset}.
  \label{eq:macMahonWithArguments}
\end{align}
This vertex operator form generalizes to nearly every other object we will discuss, and by understanding local operations on the $\Gamma$ operators --- the content of the next section --- we can derive substantially simpler and more useful expressions.

\section{Toggles} \label{sec:toggles}

Having explored the commutation relations between the $Q$ and $\Gamma$ operators, we now turn to how the $\Gamma$ operators commute with one another. To explain the relationship bijectively, we will require a particular local move called a \textbf{toggle}, described independently in \cite{pak} and \cite{hopkinsNotes}.

\begin{Def}
  \label{def:toggle1}
  Given three partitions $\lambda \succ \nu \succ \mu$, the toggle of $\nu$ relative to $\lambda$ and $\mu$ is a partition $T(\nu)$ defined in the following manner:
  $$
    T(\nu)_i = \min\left\{ \lambda_i, \mu_{i - 1} \right\} + \max\left\{ \lambda_{i + 1}, \mu_i \right\} - \nu_i,
  $$
  where we take $\mu_0 = \infty$. Intuitively, when we write $\lambda$, $\nu$, and $\mu$ side-by-side, each $\nu_i$ is bounded above by the minimum of the entries immediately above and to the left and below by the maximum of the entries immediately below and to the right; toggling is then just the unique map that is an involution on each entry's poset of possible values.
\end{Def}

We can also define toggles when the middle partition either interlaces both of the others or is interlaced by them.

\begin{Def}
  If $\lambda \prec \nu \succ \mu$, then the toggle of $\nu$ relative to $\lambda$ and $\mu$ is a map that produces a pair $(T(\nu), n)$, where $\lambda \succ T(\nu) \prec \mu$ and $n \in \mathbb{N}_{\geq 0}$. Similarly to \defref{def:toggle1}, $T(\nu)$ is given by
  $$
    T(\nu)_i = \min\left\{ \lambda_i, \mu_i \right\} + \max\left\{ \lambda_{i + 1}, \mu_{i + 1} \right\} - \nu_{i - 1}
  $$
  for $i \geq 2$. We handle $\nu_1$ separately: we say that the toggle \textbf{pops off} the the value $n = \nu_1 - \max\left\{ \lambda_1, \mu_1 \right\}$. Similarly, if $\nu$ is a partition with with $\lambda \succ \nu \prec \mu$, and $n \in \mathbb{N}_{\geq 0}$, then the toggle of $\nu$ relative to $\lambda$ and $\mu$ sends $(\nu, n)$ to a partition $T(\nu, n)$ with $\lambda \prec T(\nu, n) \succ \mu$, defined by
  $$
    T(\nu, n)_i = \min\left\{ \lambda_i, \mu_i \right\} + \max\left\{ \lambda_{i + 1}, \mu_{i + 1} \right\} - \nu_i
  $$
  for $i \geq 2$ and $T(\nu, n)_1 = n + \min\left\{ \lambda_1, \mu_1 \right\}$.
\end{Def}

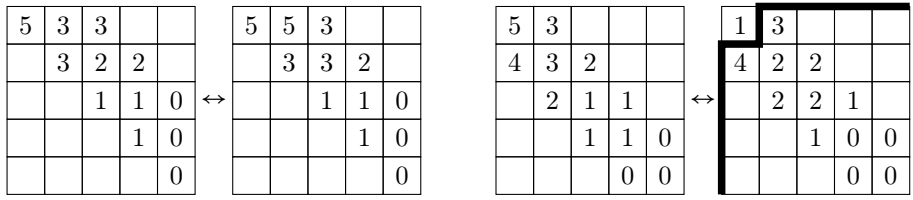
\begin{figure}
  \centering
  \begin{tikzpicture}
    
      \begin{scope}[shift={(0,0)}]
        \tableaubox{0}{0}{5}
\tableaubox{0}{1}{3}
\tableaubox{0}{2}{3}
\tableaubox{0}{3}{}
\tableaubox{0}{4}{}
\tableaulabel{0}{5}{}
\tableaubox{0}{6}{5}
\tableaubox{0}{7}{5}
\tableaubox{0}{8}{3}
\tableaubox{0}{9}{}
\tableaubox{0}{10}{}
\tableaubox{1}{0}{}
\tableaubox{1}{1}{3}
\tableaubox{1}{2}{2}
\tableaubox{1}{3}{2}
\tableaubox{1}{4}{}
\tableaulabel{1}{5}{}
\tableaubox{1}{6}{}
\tableaubox{1}{7}{3}
\tableaubox{1}{8}{3}
\tableaubox{1}{9}{2}
\tableaubox{1}{10}{}
\tableaubox{2}{0}{}
\tableaubox{2}{1}{}
\tableaubox{2}{2}{1}
\tableaubox{2}{3}{1}
\tableaubox{2}{4}{0}
\tableaulabel{2}{5}{$\leftrightarrow$}
\tableaubox{2}{6}{}
\tableaubox{2}{7}{}
\tableaubox{2}{8}{1}
\tableaubox{2}{9}{1}
\tableaubox{2}{10}{0}
\tableaubox{3}{0}{}
\tableaubox{3}{1}{}
\tableaubox{3}{2}{}
\tableaubox{3}{3}{1}
\tableaubox{3}{4}{0}
\tableaulabel{3}{5}{}
\tableaubox{3}{6}{}
\tableaubox{3}{7}{}
\tableaubox{3}{8}{}
\tableaubox{3}{9}{1}
\tableaubox{3}{10}{0}
\tableaubox{4}{0}{}
\tableaubox{4}{1}{}
\tableaubox{4}{2}{}
\tableaubox{4}{3}{}
\tableaubox{4}{4}{0}
\tableaulabel{4}{5}{}
\tableaubox{4}{6}{}
\tableaubox{4}{7}{}
\tableaubox{4}{8}{}
\tableaubox{4}{9}{}
\tableaubox{4}{10}{0}

      \end{scope}
    
      \begin{scope}[shift={(6.5,0)}]
        \tableaubox{0}{0}{5}
\tableaubox{0}{1}{3}
\tableaubox{0}{2}{}
\tableaubox{0}{3}{}
\tableaubox{0}{4}{}
\tableaulabel{0}{5}{}
\tableaubox{0}{6}{1}
\tableaubox{0}{7}{3}
\tableaubox{0}{8}{}
\tableaubox{0}{9}{}
\tableaubox{0}{10}{}
\tableaubox{1}{0}{4}
\tableaubox{1}{1}{3}
\tableaubox{1}{2}{2}
\tableaubox{1}{3}{}
\tableaubox{1}{4}{}
\tableaulabel{1}{5}{}
\tableaubox{1}{6}{4}
\tableaubox{1}{7}{2}
\tableaubox{1}{8}{2}
\tableaubox{1}{9}{}
\tableaubox{1}{10}{}
\tableaubox{2}{0}{}
\tableaubox{2}{1}{2}
\tableaubox{2}{2}{1}
\tableaubox{2}{3}{1}
\tableaubox{2}{4}{}
\tableaulabel{2}{5}{$\leftrightarrow$}
\tableaubox{2}{6}{}
\tableaubox{2}{7}{2}
\tableaubox{2}{8}{2}
\tableaubox{2}{9}{1}
\tableaubox{2}{10}{}
\tableaubox{3}{0}{}
\tableaubox{3}{1}{}
\tableaubox{3}{2}{1}
\tableaubox{3}{3}{1}
\tableaubox{3}{4}{0}
\tableaulabel{3}{5}{}
\tableaubox{3}{6}{}
\tableaubox{3}{7}{}
\tableaubox{3}{8}{1}
\tableaubox{3}{9}{0}
\tableaubox{3}{10}{0}
\tableaubox{4}{0}{}
\tableaubox{4}{1}{}
\tableaubox{4}{2}{}
\tableaubox{4}{3}{0}
\tableaubox{4}{4}{0}
\tableaulabel{4}{5}{}
\tableaubox{4}{6}{}
\tableaubox{4}{7}{}
\tableaubox{4}{8}{}
\tableaubox{4}{9}{0}
\tableaubox{4}{10}{0}

\tableaudraw[line width=1mm]{(4.5, 5.5) -- (0.5, 5.5) -- (0.5, 6.5) -- (-0.5, 6.5) -- (-0.5, 10.5)}
      \end{scope}
    
  \end{tikzpicture} \caption{Toggling $(3, 2, 1)$ relative to $(5, 3, 1, 1)$ and $(3, 2)$ (middle- and far-left), and toggling $(5, 3, 1, 1)$ relative to $(3, 2, 1)$ and $(4, 2, 1)$ (middle- and far-right).}
  \label{fig:toggles}
\end{figure}

\begin{Ex}
  Let $\lambda = (5, 3, 1, 1)$, $\mu = (3, 2)$, and $\nu = (3, 2, 1)$, so that $\lambda \succ \nu \succ \mu$. Then the toggle of $\nu$ relative to $\lambda$ and $\mu$ is $T(\nu) = (5, 3, 1)$, as shown in the left half of \figref{fig:toggles}.

  On the other hand, with $\eta = (4, 2, 1)$, $\eta \prec \lambda \succ \nu$, and we can toggle $\lambda$ relative to $\eta$ and $\nu$ to produce the partition $(2, 2)$ and the popped-off value $1$. This is shown in the right half of \figref{fig:toggles} --- we write the value of $1$ where it was popped off and separate it from the remaining diagram by a bold border. As we will soon see, this slightly strange notation will enable us to perform multiple subsequent toggles.
\end{Ex}

Toggling is a useful local move in many areas of combinatorics surrounding plane partitions and similar objects (for example, they are used to give an alternate definition for the RSK algorithm in \cite{hopkinsNotes}), but we need them to fully explain how the $\Gamma$ operators commute with one another. The commutation relation $\Gamma_-(b)\Gamma_+(a) = \frac{1}{1 - ab}\Gamma_+(a)\Gamma_-(b)$ is given in \cite[Appendix B.2]{okounkovReshetikhin}, and it turns out to be an algebraic equivalent to toggling, with the particular type of toggle dependent on the signs of the $\Gamma$ operators.

\begin{Prop}
  \label{prop:oppositeSignGammaCommutation}
  Let $\lambda$ and $\mu$ be partitions. Then there is a bijection between partitions $\nu$ with $\mu \prec \nu \succ \lambda$ and pairs $(\nu', n)$ of partitions $\nu'$ with $\mu \succ \nu' \prec \lambda$ and nonnegative integers $n$, given by toggling $\nu$ with respect to $\lambda$ and $\mu$. Moreover, the bijection preserves weight in the following manner:
  \begin{align}
    \begin{split}
      |\nu| - |\lambda| &= |\lambda| - |T(\nu)| + n \text{ and}\\
      |\nu| - |\mu| &= |\mu| - |T(\nu)| + n,
      \label{eq:weightPreservation}
    \end{split}
  \end{align}
  where $n$ is the entry popped off in the toggle.
\end{Prop}

\begin{Proof}
  Given such a $\nu$, we first verify that the toggled partition $T(\nu)$ is interlaced by both $\lambda$ and $\mu$. Let $i \geq 2$ and consider $\nu_i$. Since $\mu \prec \nu \succ \lambda$, we have the following inequalities:
  $$
    \arraycolsep=1pt\def\arraystretch{.75}
    \begin{array}{ccccc}
      \nu_{i - 1} & \geq & \lambda_{i - 1} & & \\
      \rotatebox{90}{$\leq$} & & \rotatebox{90}{$\leq$} & &\\
      \mu_{i - 1} & \geq & \nu_i & \geq & \lambda_i\\
      & & \rotatebox{90}{$\leq$} & & \rotatebox{90}{$\leq$}\\
      & & \mu_i & \geq & \nu_{i + 1}
    \end{array}
  $$
  Equivalently, $\min(\lambda_{i - 1}, \mu_{i - 1}) \geq \nu_i \geq \max(\lambda_i, \mu_i)$. By negating the inequality and performing some simple algebraic moves, we find that
  \begin{align*}
    -\min(\lambda_{i - 1}, \mu_{i - 1}) \leq &-\nu_i \leq -\max(\lambda_i, \mu_i)\\
    \Rightarrow 0 \leq \min(\lambda_{i - 1}, \mu_{i - 1})&-\nu_i \leq \min(\lambda_{i - 1}, \mu_{i - 1}) - \max(\lambda_i, \mu_i)\\
    \Rightarrow \max(\lambda_i, \mu_i) \leq \min(\lambda_{i - 1}, \mu_{i - 1}) + \max(\lambda_i, \mu_i) &-\nu_i \leq \min(\lambda_{i - 1}, \mu_{i - 1})\\
    \Rightarrow \max(\lambda_i, \mu_i) \leq & T(\nu)_{i - 1} \leq \min(\lambda_{i - 1}, \mu_{i - 1}).
  \end{align*}
  Putting this back into grid form, we have that $\mu \succ T(\nu) \prec \lambda$, as required:
  $$
    \arraycolsep=1pt\def\arraystretch{.75}
    \begin{array}{ccccc}
      T(\nu)_{i - 2} & \geq & \lambda_{i - 1} & & \\
      \rotatebox{90}{$\leq$} & & \rotatebox{90}{$\leq$} & &\\
      \mu_{i - 1} & \geq & T(\nu)_{i - 1} & \geq & \lambda_i\\
      & & \rotatebox{90}{$\leq$} & & \rotatebox{90}{$\leq$}\\
      & & \mu_i & \geq & T(\nu)_{i}
    \end{array}
  $$

  We now show that toggling is weight-preserving in the manner specified in \eqref{eq:weightPreservation}. The crucial observation is that $\min(\lambda_i, \mu_i) + \max(\lambda_i, \mu_i) = \lambda_i + \mu_i$. Therefore, the weight of $T(\nu)$ is
  \begin{align*}
    |T(\nu)| &= \sum_{i \geq 2} \min(\lambda_{i - 1}, \mu_{i - 1}) + \max(\lambda_i, \mu_i) - \nu_i\\
    &= -\min(\lambda_1, \mu_1) + \sum_{i \geq 1} (\lambda_i + \mu_i) - \sum_{i \geq 2} \nu_i\\
    &= -\min(\lambda_1, \mu_1) + |\lambda| + |\mu| - (|\nu| - \nu_1)\\
    &= |\lambda| + |\mu| - |\nu| + n,
  \end{align*}
  where $n$ is the entry popped off in the process of toggling. This proves both of the claims regarding weight.
\end{Proof}

A similar commutation relation holds for $\Gamma$ operators of the same sign and is also given in \cite[Appendix B.2]{okounkovReshetikhin}: $\Gamma_+(a)\Gamma_+(b) = \Gamma_+(b)\Gamma_+(a)$, and identically for $\Gamma_-$. Again, we give a bijective proof with toggles.

\begin{Prop}
  \label{prop:sameSignGammaCommutation}
  Let $\lambda$ and $\mu$ be partitions. Then there is a bijection between partitions $\nu$ with $\mu \succ \nu \succ \lambda$ and partitions $\nu'$ with $\mu \succ \nu' \succ \lambda$, given by toggling $\nu$ with respect to $\lambda$ and $\mu$, and it is weight-preserving in the following manner:
  \begin{align*}
    |\nu| - |\lambda| &= |\mu| - |T(\nu)|\\
    |\mu| - |\nu| &= |T(\nu)| - |\lambda|
  \end{align*}
\end{Prop}

\begin{Proof}
  Given such a $\nu$, the interlacing fact is immediate: the toggling operation necessarily preserves both directions of interlacing. For the weight preservation, the toggled partition $T(\nu)$ is defined by
  $$
    T(\nu)_i = \min(\lambda_i, \mu_{i - 1}) + \max(\lambda_{i + 1}, \mu_i) - \nu_i,
  $$
  where we take $\mu_0 = \infty$, so
  $$
    |T(\nu)| = \sum_i (\lambda_i + \mu_i - \nu_i) = |\lambda| + |\mu| - |\nu|,
  $$
  and so
  \begin{align*}
    b^{|\mu| - |T(\nu)|} a^{|T(\nu)| - |\lambda|} &= b^{|\mu| - |\lambda| - |\mu| + |\nu|} a^{|\lambda| + |\mu| - |\nu| - |\lambda|}\\
    &= b^{|\nu| - |\lambda|} a^{|\mu| - |\nu|},
  \end{align*}
  as required.
\end{Proof}

To apply these propositions to produce explicit bijections on objects counted by generating functions expressed as products of vertex operators, we require one more technical result.

\begin{Lem}
  \label{lem:bijectivizingToggles}
  Let $z_1(q)$ be a generating function given by
  $$
    z_1(q) = \leftbra{\emptyset}\,\prod_{n \in \mathbb{Z}} \Gamma_{s_n}\left(q^{p_n}\right)\,\rightket{\emptyset},
  $$
  where each $s_n = \pm 1$, and we write $\Gamma_{\pm 1}$ to mean $\Gamma_{\pm}$, respectively. Let $m \in \mathbb{Z}$, and let $z_2(q)$ be equal to $z_1(q)$, except with the order of $\Gamma_{s_m}\left(q^{p_m}\right)$ and $\Gamma_{s_{m + 1}}\left(q^{p_{m + 1}}\right)$ swapped in the product. Let $S_1$ and $S_2$ be the sets of objects counted by $z_1$ and $z_2$, with weight marked by $q$, as in \defref{def:gammaOperators}. Since $\Gamma_{s_m}\left(q^{p_m}\right)$ and $\Gamma_{s_{m + 1}}\left(q^{p_{m + 1}}\right)$ are adjacent in $z_1$, the objects in $S_1$ have some diagonal that both of these $\Gamma$ operators count; call it $d_m$. Then there is a weight-preserving bijection $f : S_1 \to S_2$, and it is given by toggling this diagonal $d_m$ with respect to the two diagonals adjacent to it.
\end{Lem}

\begin{Proof}
  We show the case where $s_m = -1$ and $s_{m + 1} = 1$; the other cases are exactly analogous. 
  Write $z_1(q)$ as
  $$
    z_1(q) = \sum_{\mu \prec \nu \succ \lambda} \leftbra{\emptyset}\,\prod_{n \leq m} \Gamma_{s_n}\left(q^{p_n}\right)\,\rightket{\nu} \leftbra{\nu}\,\prod_{n \geq m + 1} \Gamma_{s_n}\left(q^{p_n}\right)\,\rightket{\emptyset}.
  $$
  In this presentation, $\nu$ is the diagonal counted by both $\Gamma_{s_m}\left(q^{p_m}\right)$ and $\Gamma_{s_{m + 1}}\left(q^{p_{m + 1}}\right)$. Applying \propref{prop:oppositeSignGammaCommutation}, toggling that diagonal gives a bijection between the partitions $\nu$ counted in the sum to pairs $(\nu', k)$ of partitions $\nu'$ satisfying $\mu \succ \nu' \prec \lambda$ and integers $k \geq 0$ that satisfy \eqref{eq:weightPreservation}. The generating function that counts these new objects is then
  $$
    \frac{1}{1 - q^{p_m + p_{m + 1}}} \sum_{\mu \succ \nu' \prec \lambda} \leftbra{\emptyset}\,\left( \prod_{n < m} \Gamma_{s_n}\left(q^{p_n}\right)\right) \Gamma_{s_{m + 1}}\left(q^{p_{m + 1}} \right)\,\rightket{\nu'} \leftbra{\nu'}\,\Gamma_{s_m}\left(q^{p_m}\right) \left( \prod_{n > m + 1} \Gamma_{s_n}\left(q^{p_n}\right) \right)\,\rightket{\emptyset},
  $$
  which is exactly $z_2(q)$.
\end{Proof}

This lemma bijectivizes the commutation of the $\Gamma$ operators: given a vertex operator expression for a generating function $f(q)$, the commutation of two same-sign operators to form a new generating function $g(q)$ is a weight-preserving bijection of the objects counted by $f$ to those counted by $g$, given by toggling the diagonal corresponding to the commuted operators. Similarly, commuting a $\Gamma_+$ to the left past a $\Gamma_-$ is a weight-preserving bijection given by toggling, and it also produces a nonnegative integer based on the arguments of the $\Gamma$ operators. These bijections form the bedrock of many of our subsequent results: if the generating function for a plane-partition-like object has a vertex operator expansion and an identity can be proven algebraically by performing successive commutations, then it can be bijectivized by interpreting each successive commutation as a toggle of a particular diagonal. Our first application of this method is to compute a bijectivization of the product expansion of MacMahon's function.

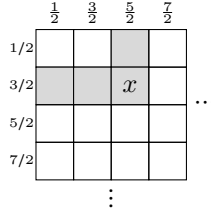
\begin{figure}
  \centering
  \begin{tikzpicture}
    
      \begin{scope}[shift={(0,0)}]
        \tableaubox{0}{0}{}
\tableaubox{0}{1}{}
\tableaubox[fill=gray!30]{0}{2}{}
\tableaubox{0}{3}{}
\tableaubox[fill=gray!30]{1}{0}{}
\tableaubox[fill=gray!30]{1}{1}{}
\tableaubox[fill=gray!30]{1}{2}{$x$}
\tableaubox{1}{3}{}
\tableaubox{2}{0}{}
\tableaubox{2}{1}{}
\tableaubox{2}{2}{}
\tableaubox{2}{3}{}
\tableaubox{3}{0}{}
\tableaubox{3}{1}{}
\tableaubox{3}{2}{}
\tableaubox{3}{3}{}
\tableaulabel{1.5}{4}{$\cdots$}
\tableaulabel{4}{1.5}{$\vdots$}

\tableauedgelabelleft{0}{0}{$1/2$}
\tableauedgelabelleft{1}{0}{$3/2$}
\tableauedgelabelleft{2}{0}{$5/2$}
\tableauedgelabelleft{3}{0}{$7/2$}
\tableauedgelabeltop{-0.25}{0}{$\frac{1}{2}$}
\tableauedgelabeltop{-0.25}{1}{$\frac{3}{2}$}
\tableauedgelabeltop{-0.25}{2}{$\frac{5}{2}$}
\tableauedgelabeltop{-0.25}{3}{$\frac{7}{2}$}
      \end{scope}
    
  \end{tikzpicture} \caption{A cell $x \in \mathbb{N}^2$ and its hook. The edges of the diagram are labeled with the exponents of the $\Gamma$ operators they correspond to. The hook of $x$ intersects the boundary at edges corresponding to operators $\Gamma_+\left( q^{5/2} \right)$ and $\Gamma_-\left( q^{3/2} \right)$, and the hook length of $x$ is $h(x) = 4 = \frac{5}{2} + \frac{3}{2}$.}
  \label{fig:planePartitionEdgeLabels}
\end{figure}

\begin{Thm}
  \label{thm:macMachonsFunctionByToggling}
  There is a weight-preserving bijection $\tau$ between plane partitions and tableaux of shape $\mathbb{N}^2$ (i.e.\ functions $\mathbb{N}^2 \to \mathbb{N}_{\geq 0}$ with no restrictions on outputs) that are weighted by hook length.
\end{Thm}

\begin{Proof}
  Our bijection is functionally identical to the Pak--Sulzgruber algorithm \cite{pak, sulzgruber}, except described on standard plane partitions instead of reverse ones. What is new is the method by which we derive it, which is generalizable to the objects we discuss in future sections. We follow the methods in \cite[page 11]{vertexOperators} used to prove that the vertex operator expansion
  \begin{align}
    M(q) = \leftbra{\emptyset} \cdots \Gamma_-\left(q^{5/2}\right) \Gamma_-\left(q^{3/2}\right) \Gamma_-\left(q^{1/2}\right) \Gamma_+\left(q^{1/2}\right) \Gamma_+\left(q^{3/2}\right) \Gamma_+\left(q^{5/2}\right) \cdots \rightket{\emptyset},
    \label{eq:macmahonFromVertexOperators}
  \end{align}
  can be iteratively commuted into the expression
  $$
    M(q) = \prod_{\square \in \mathbb{N}^2} \frac{1}{1 - q^{h(\square)}}.
  $$
  
  Consider the first commutation in \eqref{eq:macmahonFromVertexOperators}. Define
  $$
    f_0 = \leftbra{\emptyset} \cdots \Gamma_-\left(q^{5/2}\right) \Gamma_-\left(q^{3/2}\right) \Gamma_-\left(q^{1/2}\right) \Gamma_+\left(q^{1/2}\right) \Gamma_+\left(q^{3/2}\right) \Gamma_+\left(q^{5/2}\right) \cdots \rightket{\emptyset} = M(q)
  $$
  and
  $$
    f_1 = \leftbra{\emptyset} \cdots \Gamma_-\left(q^{5/2}\right) \Gamma_-\left(q^{3/2}\right) \Gamma_+\left(q^{1/2}\right) \Gamma_-\left(q^{1/2}\right) \Gamma_+\left(q^{3/2}\right) \Gamma_+\left(q^{5/2}\right) \cdots \rightket{\emptyset},
  $$
  where the middle two operators have been commuted. Denote the sets of objects counted by these two generating functions by $P_0$ and $P_1$, respectively. While $P_0$ is simply the set of plane partitions, weighted as usual, $P_1$ is less well-behaved. Its elements are of the form of the far-right object in \figref{fig:toggles}: each is a placement of nonnegative integers into $\mathbb{N}^2 \setminus \left\{ (1, 1) \right\}$ that weakly decrease along rows and columns, and the weight of such an object is given by $f_1$. By \lemref{lem:bijectivizingToggles}, the map $\tau_0 : P_0 \to P_1 \times \mathbb{N}_{\geq 0}$ given by toggling the main diagonal is a bijection.

  Define $f_n$, $P_n$, and $\tau_n : P_n \to P_{n + 1} \times \mathbb{N}_{\geq 0}$ for $n \geq 1$ similarly; we may enumerate them in any order compatible with the order in which corners appear when toggling. We will define our bijection $\tau$ from the set of plane partitions to the set of tableaux of shape $\mathbb{N}^2$ effectively by composing every $\tau_k$, but the notation makes this slightly cumbersome. Let $(s_n) \subset \mathbb{N}^2$ be the sequence of cells such that $s_n$ is popped off by $\tau_n$. The sequence depends on our ordering of the $\tau_n$; one possible ordering is by off-diagonals, i.e.
  $$
    (s_n)_{n \in \mathbb{N}} = ((1, 1), (1, 2), (2, 1), (1, 3), (2, 2), (3, 1), (1, 4), ...),
  $$
  since each subsequent off-diagonal consists entirely of corners after every cell in the previous off-diagonal has been popped off. More generally, we may choose any ordering so that for any $k \in \mathbb{N}$, the tableau
  $$
    t : \left\{ s_n \mid n \in \left\{ 1, 2, ..., k \right\} \right\} \to \mathbb{N}
  $$
  given by $t(s_n) = n$ is a standard Young tableau.

  Let $B$ be the set of tableaux of shape $\mathbb{N}^2$ (i.e.\ functions $\mathbb{N}^2 \to \mathbb{N}_{\geq 0}$) and define functions $\tau'_n : P_n \times B \to P_{n + 1} \times B$, where $\tau'_n$ applies $\tau_n$ to its first argument and leaves its second argument unchanged, except for setting the value in cell $s_n$ to the value popped off by $\tau_n$. For a plane partition $\pi$ and the zero tableau $0$, the composition $\left( \cdots \circ \tau'_2 \circ \tau'_1 \circ \tau'_0 \right)(\pi, 0)$ then produces a tuple $(\emptyset, \beta)$ for a tableau $\beta$ of shape $\mathbb{N}^2$; we define the map $\tau : P_0 \to B$ by setting $\tau(\pi) = \beta$. This map $\tau$ is well-defined since every plane partition has only finitely many nonzero entries --- given a plane partition, we may stop toggling diagonals once there are no longer any nonzero entries left. Moreover, since each $\tau_i$ is a bijection, $\tau$ is also a bijection.

  It remains to show that $\tau$ is weight-preserving when its output tableaux are weighted by hook length. Specifically, if $\tau'_n(\pi, \beta) = (\pi', \beta')$ and the number popped off in the toggle is $a$, then we wish to show that $|\pi| = |\pi'| + a \cdot h(s_n)$, where the weights of $\pi$ and $\pi'$ are as measured by $f_n$ and $f_{n + 1}$, respectively. Now $f_n$ and $f_{n + 1}$ differ only by a single commutation that moves some $\Gamma_-(q^i)$ to the right past some $\Gamma_+(q^j)$, producing a factor of $\frac{1}{1 - q^{i + j}}$ that corresponds to the number popped off in the toggle. Our claim of weight-preservation then reduces to the claim that $h(s_n) = i + j$.

  Set $s_n = (k, l)$. In the original square diagram, $\operatorname{hook}(s_n)$ intersects the boundary of $\mathbb{N}^2$ at a horizontal edge corresponding to $\Gamma_+\left( q^{l - \frac{1}{2}} \right)$ and a vertical one corresponding to $\Gamma_-\left( q^{k - \frac{1}{2}} \right)$, as in \figref{fig:planePartitionEdgeLabels}. When we toggle a diagonal and commute a $\Gamma_-$ to the right past a $\Gamma_+$, the corresponding edges swap places in the diagram --- in effect, the horizontal edge corresponding to $\Gamma_+$ moves down, and the vertical edge corresponding to $\Gamma_-$ moves right. Therefore, no matter the order in which we toggle diagonals before reaching $s_n$, the horizontal edge corresponding to $\Gamma_+\left( q^{l - \frac{1}{2}} \right)$ is always above $s_n$ and the vertical one corresponding to $\Gamma_-\left( q^{k - \frac{1}{2}} \right)$ is always to its left. When we finally toggle the diagonal containing $s_n$ after $n - 1$ previous toggles, those two edges are bordering $s_n$, and so the values of $i$ and $j$ from the previous paragraph are in fact $k - \frac{1}{2}$ and $l - \frac{1}{2}$. Since $h(s_n) = k + l - 1 = i + j$, our claim is proven.
\end{Proof}

In short, this map $\tau$ sends a plane partition to a tableau weighted by hook length, given by toggling diagonals until the diagram is empty (i.e.\ the $\Gamma$ operators have been completely commuted) and recording the popped numbers in the locations from which they were removed. The map is functionally identical to that described in \cite{pak, sulzgruber}, but the vertex operator description provides an alternate lens through which to view the bijection and a clear explanation of why it is independent of the order in which we toggle diagonals (specifically, the commutators of $\Gamma$ operators that are produced can be freely factored out of the entire expression).

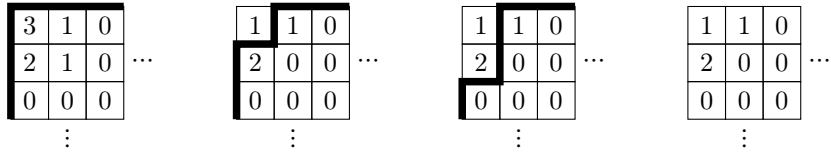
\begin{figure}
  \centering
  \begin{tikzpicture}
    
      \begin{scope}[shift={(0,0)}]
        \tableaubox{0}{0}{3}
\tableaubox{0}{1}{1}
\tableaubox{0}{2}{0}
\tableaubox{1}{0}{2}
\tableaubox{1}{1}{1}
\tableaubox{1}{2}{0}
\tableaubox{2}{0}{0}
\tableaubox{2}{1}{0}
\tableaubox{2}{2}{0}
\tableaulabel{1}{3}{$\cdots$}
\tableaulabel{3}{1}{$\vdots$}
\tableaudraw[line width=1mm]{(-0.5, 2.5) -- (-0.5, -0.5) -- (0.5, -0.5) -- (0.5, -0.5) -- (1.5, -0.5) -- (1.5, -0.5) -- (2.5, -0.5)}
      \end{scope}
    
      \begin{scope}[shift={(3,0)}]
        \tableaubox{0}{0}{1}
\tableaubox{0}{1}{1}
\tableaubox{0}{2}{0}
\tableaubox{1}{0}{2}
\tableaubox{1}{1}{0}
\tableaubox{1}{2}{0}
\tableaubox{2}{0}{0}
\tableaubox{2}{1}{0}
\tableaubox{2}{2}{0}
\tableaulabel{1}{3}{$\cdots$}
\tableaulabel{3}{1}{$\vdots$}
\tableaudraw[line width=1mm]{(-0.5, 2.5) -- (-0.5, 0.5) -- (0.5, 0.5) -- (0.5, -0.5) -- (1.5, -0.5) -- (1.5, -0.5) -- (2.5, -0.5)}
      \end{scope}
    
      \begin{scope}[shift={(6,0)}]
        \tableaubox{0}{0}{1}
\tableaubox{0}{1}{1}
\tableaubox{0}{2}{0}
\tableaubox{1}{0}{2}
\tableaubox{1}{1}{0}
\tableaubox{1}{2}{0}
\tableaubox{2}{0}{0}
\tableaubox{2}{1}{0}
\tableaubox{2}{2}{0}
\tableaulabel{1}{3}{$\cdots$}
\tableaulabel{3}{1}{$\vdots$}
\tableaudraw[line width=1mm]{(-0.5, 2.5) -- (-0.5, 0.5) -- (0.5, 0.5) -- (0.5, 0.5) -- (1.5, 0.5) -- (1.5, -0.5) -- (2.5, -0.5)}
      \end{scope}
    
      \begin{scope}[shift={(9,0)}]
        \tableaubox{0}{0}{1}
\tableaubox{0}{1}{1}
\tableaubox{0}{2}{0}
\tableaubox{1}{0}{2}
\tableaubox{1}{1}{0}
\tableaubox{1}{2}{0}
\tableaubox{2}{0}{0}
\tableaubox{2}{1}{0}
\tableaubox{2}{2}{0}
\tableaulabel{1}{3}{$\cdots$}
\tableaulabel{3}{1}{$\vdots$}

      \end{scope}
    
  \end{tikzpicture} \caption{A weight-7 plane partition $\pi$ (far left) being mapped bijectively to a weight-7 tableau $\tau(\pi)$ that is weighted by hook length (far right). The center-left and center-right figures are the first two steps in the bijection.}
  \label{fig:planePartitionToggles}
\end{figure}

\begin{Ex}
  Let $\pi$ be the plane partition given in \figref{fig:planePartitionToggles}. In the expression for $M(q)$ given in \eqref{eq:macmahonFromVertexOperators}, $\pi$ is represented by a term of $q^7$. Commuting the $\Gamma_-\left(q^{1/2}\right)$ with the $\Gamma_+\left(q^{1/2}\right)$ results in
  $$
    M(q) = \frac{1}{1 - q} \leftbra{\emptyset} \cdots \Gamma_-\left(q^{5/2}\right) \Gamma_-\left(q^{3/2}\right) \Gamma_+\left(q^{1/2}\right) \Gamma_-\left(q^{1/2}\right) \Gamma_+\left(q^{3/2}\right) \Gamma_+\left(q^{5/2}\right) \cdots \rightket{\emptyset},
  $$
  and the term of $q^7$ now represents the object second from left in \figref{fig:planePartitionToggles}: a $1 \times 1$ tableau with the entry $1$, and an object in $P_1$. We draw the two superimposed with a bold dividing border to emphasize how the process preserves the overall shape of $\pi$. The $1 \times 1$ tableau is weighted by hook length (i.e.\ it has weight $1$), and the object on the right has weight $6$ as counted by $f_1$ (in the language of \defref{def:gammaOperators}), so the combined weight of $7$ is preserved. Note that this weight of $6$ is \textit{not} the sum of the entries; in general, only the starting plane partition has weight equal to the sum of its entries.

  We now have two choices of corner to commute --- if we choose to commute $\Gamma_-\left(q^{3/2}\right)$ with $\Gamma_+\left(q^{1/2}\right)$, the resulting pair of objects is second from right in \figref{fig:planePartitionToggles}. The tableau now has weight $5$, since the box containing $2$ has hook length $2$, while the object in $P_2$ has weight $2$ as counted by $f_2$. After additionally toggling the diagonal containing the $1$, there are no more nonzero entries in the plane-partition-like object, and so all future toggles place a zero into the tableau. The result is the final tableau $\tau(\pi)$ on the far right, whose weight (accounting for hook length) is correctly equal to $7$.
\end{Ex}

\section{One-Leg Objects} \label{sec:oneLegObjects}

The techniques of the previous section apply to far more than just plane partitions: any object whose generating function is expressible as a product of these $\Gamma$ operators is a candidate for a decomposition given by iterative toggling. We begin with two definitions of plane-partition-like objects --- the standard notion of reverse plane partitions, along with a less standard notion of a skew plane partition (see \cite[Section 3.4]{vertexOperators}) --- before connecting them with a result that we bijectivize.

\begin{Def}
  Let $\lambda \subset \mathbb{N}^2$ be a Young diagram.
  \begin{enumerate}[label=\arabic*]
    \item A \textbf{reverse plane partition}, or RPP, of shape $(\emptyset, \emptyset, \lambda)$ is a function $\rho : \lambda \to \mathbb{N}_{\geq 0}$ such that $\rho(i, j) \leq \rho(i + 1, j)$ and $\rho(i, j) \leq \rho(i, j + 1)$ for all choices $(i, j)$ where those quantities are defined. We write the generating function for RPPs of shape $(\emptyset, \emptyset, \lambda)$ as $W_{(\emptyset, \emptyset, \lambda)}(q)$, where $q$ marks the weight (i.e.\ the sum of all entries), following the notation of \cite{PT}.
    \item A \textbf{skew plane partition}, or SPP, of shape $(\emptyset, \emptyset, \lambda)$ is a function $\sigma : \mathbb{N}^2 \setminus \lambda \to \mathbb{N}_{\geq 0}$ such that $\sigma(i, j) \geq \sigma(i + 1, j)$ and $\sigma(i, j) \geq \sigma(i, j + 1)$ for all $(i, j) \in \mathbb{N}^2 \setminus \lambda$ and only finitely many $\sigma(i, j)$ are nonzero. We write the generating function for SPPs of shape $(\emptyset, \emptyset, \lambda)$ as $V_{(\emptyset, \emptyset, \lambda)}(q)$, where $q$ once again marks the weight.
  \end{enumerate}
  \label{def:rppAndSpp}
\end{Def}

We call these \textbf{one-leg} objects, again following \cite{PT}, since only one of the three entries is nonempty in the shape term $(\emptyset, \emptyset, \lambda)$. The notation suggests further generalizations are possible, and this is indeed the case; we discuss the case where at most two of the entries of the shape term are nonempty in \secref{sec:twoLegObjects}, and we summarize the fully general case in \secref{sec:futureDirections}.

For these one-leg objects, the PT--DT correspondence states that
\begin{align}
  V_{(\emptyset, \emptyset, \lambda)}(q) = M(q)W_{(\emptyset, \emptyset, \lambda)}(q),
  \label{eq:PTDT1}
\end{align}
or that combinatorially, every SPP of weight $n$ and shape $(\emptyset, \emptyset, \lambda)$ corresponds to an RPP of shape $(\emptyset, \emptyset, \lambda)$ and a plane partition whose weights sum to $n$. This equation was stated in \cite{PT}; throughout this section, we supply details that were omitted and bijectivize the equation.

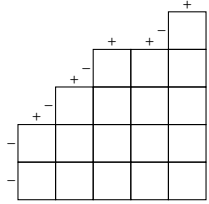
\begin{figure}
  \centering
  \begin{tikzpicture}
    
      \begin{scope}[shift={(0,0)}]
        \tableaulabel{0}{0}{}
\tableaulabel{0}{1}{}
\tableaulabel{0}{2}{}
\tableaulabel{0}{3}{}
\tableauedgelabelleft{0}{4}{$-$}
\tableauedgelabeltop{0}{4}{$+$}
\tableaubox{0}{4}{}
\tableaulabel{1}{0}{}
\tableaulabel{1}{1}{}
\tableauedgelabelleft{1}{2}{$-$}
\tableauedgelabeltop{1}{2}{$+$}
\tableaubox{1}{2}{}
\tableauedgelabeltop{1}{3}{$+$}
\tableaubox{1}{3}{}
\tableaubox{1}{4}{}
\tableaulabel{2}{0}{}
\tableauedgelabelleft{2}{1}{$-$}
\tableauedgelabeltop{2}{1}{$+$}
\tableaubox{2}{1}{}
\tableaubox{2}{2}{}
\tableaubox{2}{3}{}
\tableaubox{2}{4}{}
\tableauedgelabelleft{3}{0}{$-$}
\tableauedgelabeltop{3}{0}{$+$}
\tableaubox{3}{0}{}
\tableaubox{3}{1}{}
\tableaubox{3}{2}{}
\tableaubox{3}{3}{}
\tableaubox{3}{4}{}
\tableauedgelabelleft{4}{0}{$-$}
\tableaubox{4}{0}{}
\tableaubox{4}{1}{}
\tableaubox{4}{2}{}
\tableaubox{4}{3}{}
\tableaubox{4}{4}{}

      \end{scope}
    
  \end{tikzpicture} \caption{The edge sign sequence for $\lambda = (4, 2, 1)$ is $e_\lambda = (\ldots, -, -, +, -, +, -, +, +, -, +, \ldots)$.}
  \label{fig:edgeSignSequence}
\end{figure}

We begin by focusing on skew plane partitions. To express $V_{(\emptyset, \emptyset, \lambda)}(q)$ in terms of the $\Gamma$ operators, the middle section of the sequence of $\Gamma$s contains a mix of $\Gamma_+$ and $\Gamma_-$ according to the sequence of horizontal and vertical edges in the diagram. With $\lambda = (4, 2, 1)$ for instance, $V_{(\emptyset, \emptyset, \lambda)}$ has the following generating function, where the signs of the $\Gamma$ operators are determined from the order of edges in \figref{fig:edgeSignSequence}. To make this sequence more legible, the $Q$ operator corresponding to the middle diagonal of the diagram is bolded.
$$
  \bra{\emptyset}\,\cdots\,Q\,\Gamma_-(1)\,Q\,\Gamma_-(1)\,Q\,\Gamma_+(1)\,Q\,\Gamma_-(1)\,Q\,\Gamma_+(1)\,\mathbf{Q}\,\Gamma_-(1)\,Q\,\Gamma_+(1)\,Q\,\Gamma_+(1)\,Q\,\Gamma_-(1)\,Q\,\Gamma_+(1)\,Q\,\cdots\,\ket{\emptyset},
$$
After commuting every $Q$ outward, this results in
$$
  \leftbra{\emptyset}\,\cdots\Gamma_-\left( q^{\frac{9}{2}} \right)\Gamma_-\left( q^{\frac{7}{2}} \right)\Gamma_+\left( q^{-\frac{5}{2}} \right)\Gamma_-\left( q^{\frac{3}{2}} \right)\Gamma_+\left( q^{-\frac{1}{2}} \right)\Gamma_-\left( q^{-\frac{1}{2}} \right)\Gamma_+\left( q^{\frac{3}{2}} \right)\Gamma_+\left( q^{\frac{5}{2}} \right)\Gamma_-\left( q^{-\frac{7}{2}} \right)\Gamma_+\left( q^{\frac{9}{2}} \right)\cdots\,\rightket{\emptyset}.
$$
The pattern of increasing odd-numerator fractions of $2$ persists, but now the ``out of place'' $\Gamma$ operators have negative-exponent arguments. This construction will prove quite useful, and it will be helpful to formalize it.

\begin{Def}
  Let $\lambda \subset \mathbb{N}^2$ be a Young diagram and label the edges of $\mathbb{N}^2 \setminus \lambda$ from bottom-left to top-right with consecutive integers so that the two edges adjacent to the main diagonal are labeled $-1$ and $0$. The \textbf{edge sign sequence} of $\lambda$, denoted $e_\lambda$, is the doubly infinite sequence whose entries are $\pm 1$, where $e_\lambda(n) = 1$ if the edge labeled $n$ is horizontal and $-1$ if it is vertical. We also define the \textbf{edge power sequence} of $\lambda$ as $p_\lambda(n)$, where
  $$
    p_\lambda(n) = \begin{cases} \left| n + \frac{1}{2} \right|, & e_\lambda(n) = \operatorname{sign}\left(n + \frac{1}{2}\right) \\ - \left| n + \frac{1}{2} \right|, & e_\lambda(n) \neq \operatorname{sign}\left(n + \frac{1}{2}\right) \end{cases}.
  $$
\end{Def}

The edge sign sequence is identical to the edge sequence defined in \cite[Definition 7.3.1]{orbifoldTopologicalVertex}, except with opposite signs (the authors also use the opposite convention of $\Gamma_+$ and $\Gamma_-$). It also bears a resemblance to the word associated with steep domino tilings in \cite[Proposition 1]{steepTilings}, as well as the sign sequence in \cite[Definition 2.1]{railYardGraphs}; both also encode binary data related to square-shaped objects into signs.

Given a Young diagram $\lambda$, its edge sign sequence $e(n) = e_\lambda(n)$, and its edge power sequence $p(n) = p_\lambda(n)$, the generating function $V_{(\emptyset, \emptyset, \lambda)}(q)$ is
\begin{align}
  V_{(\emptyset, \emptyset, \lambda)}(q) = \leftbra{\emptyset}\,\prod_{n \in \mathbb{Z}} \Gamma_{e(n)}\left(q^{p(n)}\right)\,\rightket{\emptyset},
  \label{eq:oneLegSppVertexOperators}
\end{align}
where $\Gamma_{1}$ and $\Gamma_{-1}$ are written to mean $\Gamma_+$ and $\Gamma_-$, respectively. The exponents serve a greater purpose than merely defining the shape --- they have a more direct interpretation in terms of hook lengths in the tableau. We first prove a technical lemma regarding the exponent sign sequence, and then a more substantial result.

\begin{Lem}
  \label{lem:effectOfRemovingCornerOnExponentSequence}
  Let $\lambda \subset \mathbb{N}^2$ be a Young diagram with edge power sequence $p_\lambda(n)$, let $b$ be a corner of $\lambda$ (i.e.\ a box $(i, j)$ with $(i + 1, j) \notin \lambda$ and $(i, j + 1) \notin \lambda$), and suppose that the bottom and right edges of $b$ have labels $p_\lambda(k - 1)$ and $p_\lambda(k)$ for some $k \in \mathbb{Z}$. Let $\mu$ be the Young diagram given by removing $b$ from $\lambda$, so that the left and top edges of $b$ have labels $p_\mu(k - 1)$ and $p_\mu(k)$. Then $p_\mu(k - 1) = p_\lambda(k) + 1$.
\end{Lem}

\begin{Proof}
  All edge power sequences are the same in absolute value, so $\left| p_\mu(k - 1) \right| = \left| p_\lambda(k - 1) \right|$ and $\left| p_\mu(k) \right| = \left| p_\lambda(k) \right|$. Since the two edge sign sequences $e_\lambda$ and $e_\mu$ differ only in that $e_\mu(k - 1) = -e_\lambda(k - 1)$ and $e_\mu(k) = -e_\lambda(k)$, we must have that $p_\mu(k - 1) = -p_\lambda(k - 1)$ and $p_\mu(k) = -p_\lambda(k)$.

  To finish the proof, we relate $p_\lambda(k - 1)$ to $p_\lambda(k)$. We know $e_\lambda(k - 1) = 1$ and $e_\lambda(k) = -1$ since $b$ is a corner, and what remains is a straightforward computation with three cases. If $k = 0$, then $p_\lambda(k - 1) = p_\lambda(k) = -\frac{1}{2}$. If $k > 0$, then $p_\lambda(k - 1) = k - 1 + \frac{1}{2}$ and $p_\lambda(k) = -\left( k + \frac{1}{2} \right)$. Finally, if $k < 0$, then $p_\lambda(k - 1) = -\left( k - 1 + \frac{1}{2} \right)$ and $p_\lambda(k) = k + \frac{1}{2}$. In every case, $p_\lambda(k - 1) + p_\lambda(k) = -1$. Therefore,
  \begin{align*}
    p_\mu(k - 1) &= -p_\lambda(k - 1)\\
    &= -\left( -1 - p_\lambda(k) \right)\\
    &= p_\lambda(k) + 1,
  \end{align*}
  as required.
\end{Proof}

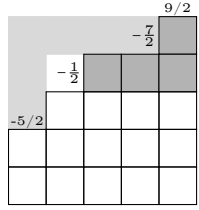
\begin{figure}
  \centering
  \begin{tikzpicture}
    
      \begin{scope}[shift={(0,0)}]
        \tableaubox[fill=gray!30,gray!30]{0}{0}{}
\tableaubox[fill=gray!30,gray!30]{0}{1}{}
\tableaubox[fill=gray!30,gray!30]{0}{2}{}
\tableaubox[fill=gray!30,gray!30]{0}{3}{}
\tableauedgelabelleft{0}{4}{$-\frac{7}{2}$}
\tableauedgelabeltop{0}{4}{9/2}
\tableaubox[fill=gray!60]{0}{4}{}
\tableaubox[fill=gray!30,gray!30]{1}{0}{}
\tableaulabel{1}{1}{}
\tableauedgelabelleft{1}{2}{$-\frac{1}{2}$}
\tableaubox[fill=gray!60]{1}{2}{}
\tableaubox[fill=gray!60]{1}{3}{}
\tableaubox[fill=gray!60]{1}{4}{}
\tableaubox[fill=gray!30,gray!30]{2}{0}{}
\tableaubox{2}{1}{}
\tableaubox{2}{2}{}
\tableaubox{2}{3}{}
\tableaubox{2}{4}{}
\tableauedgelabeltop{3}{0}{-5/2}
\tableaubox{3}{0}{}
\tableaubox{3}{1}{}
\tableaubox{3}{2}{}
\tableaubox{3}{3}{}
\tableaubox{3}{4}{}
\tableaubox{4}{0}{}
\tableaubox{4}{1}{}
\tableaubox{4}{2}{}
\tableaubox{4}{3}{}
\tableaubox{4}{4}{}

      \end{scope}
    
  \end{tikzpicture} \caption{Hooks and their boundary edges in $\lambda = (4, 2, 1)$ and the Young diagram asymptotic to it.}
  \label{fig:exponentSequence}
\end{figure}

\begin{Thm}
  \label{thm:hookLengthsInAsymptoticDiagrams}
  Let $\lambda \subset \mathbb{N}^2$ be a Young diagram with edge power sequence $p_\lambda(n)$, let $(i, j) \in \mathbb{N}^2 \setminus \lambda$, and let $k$ and $l$ be the labels of the vertical edge in row $i$ and the horizontal edge in column $j$ of the boundary of $\lambda$, respectively. Then the hook length $h_\lambda(i, j)$ satisfies the formula
  $$
    p_\lambda(k) + p_\lambda(l) = h_\lambda(i, j)
  $$
  If $(i, j) \in \lambda$ and $k$ and $l$ are as before, then the same result holds, but with an additional sign:
  $$
    p_\lambda(k) + p_\lambda(l) = -h_\lambda(i, j).
  $$
  In \figref{fig:exponentSequence}, for example, the hook with pivot $(2, 5)$ has boundary edges with labels $-\frac{1}{2}$ and $\frac{9}{2}$ and hook length $-\frac{1}{2} + \frac{9}{2} = 4$, and the hook with pivot $(1, 1)$ has boundary edges with labels $-\frac{5}{2}$ and $-\frac{7}{2}$ and hook length $-\left( -\frac{5}{2} - \frac{7}{2} \right) = 6$.
\end{Thm}

\begin{Proof}
  We prove the result by induction on the size of $\lambda$. When $\lambda = \emptyset$, there are no boxes inside of $\lambda$, so we need only show the result for boxes outside. Given $(i, j) \in \mathbb{N}^2$, its hook length is $h_\lambda(i, j) = i + j - 1$, and its hook meets the boundary of $\lambda$ --- i.e.\ the boundary of $\mathbb{N}^2$ --- at a vertical edge with label $k = -i$ and a horizontal edge with label $l = j - 1$ (recall that the edges bordering the main diagonal have labels $-1$ and $0$). Both $p_\lambda(k)$ and $p_\lambda(l)$ are then positive, and in particular,
  \begin{align*}
    p_\lambda(k) + p_\lambda(l) &= \left| -i + \frac{1}{2} \right| + \left| j - 1 + \frac{1}{2} \right|\\
    &= -\left( -i + \frac{1}{2} \right) + \left( j - \frac{1}{2} \right)\\
    &= i + j - 1\\
    &= h_\lambda(i, j).
  \end{align*}
  Now suppose the proposition holds for Young diagrams with at most $n - 1$ boxes, let $\lambda$ be one with $n$ boxes, and let $(i, j) \in \mathbb{N}^2 \setminus \lambda$. If $\lambda_i = \lambda'_j = 0$ (i.e.\ the hook meets the boundary of $\mathbb{N}^2$ itself), then
  $$
    p_\lambda(k) + p_\lambda(l) = i + j - 1 = h_\lambda(i, j)
  $$
  by identical logic to the base case. Otherwise, suppose without loss of generality that $\lambda_i \neq 0$. Then the left leg of the hook of $(i, j)$ meets the boundary of $\lambda$ at the box $b = (i, \lambda_i)$. If $\lambda_{i + 1} < \lambda_i$, then $b$ is a corner, and we may remove it from the Young diagram to produce a diagram $\mu$ with $n - 1$ boxes. By the induction hypothesis,
  $$
    p_\mu(k - 1) + p_\mu(l) = h_\mu(i, j).
  $$
  Now $h_\mu(i, j) = h_\lambda(i, j) + 1$, and by \lemref{lem:effectOfRemovingCornerOnExponentSequence}, $p_\mu(k - 1) = p_\lambda(k) + 1$. Moreover, the edge labeled $l$ in $\lambda$ is unchanged in $\mu$ and still labeled $l$, so $p_\mu(l) = p_\lambda(l)$. In total,
  $$
    p_\lambda(k) + 1 + p_\lambda(l) = h_\lambda(i, j) + 1,
  $$
  proving the result.

  If $b = (i, \lambda_i)$ is not a corner, the argument is much simpler: let $b'$ be the corner in the same column in $b$, which is then necessarily strictly below $b$, and let $\mu$ be the Young diagram given by removing it. Then the induction hypothesis guarantees that
  $$
    p_\mu(k) + p_\mu(l) = h_\mu(i, j),
  $$
  but $p_\mu(k) = p_\lambda(k)$, $p_\mu(l) = p_\lambda(l)$, and $h_\mu(i, j) = h_\lambda(i, j)$.

  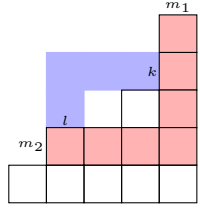
\begin{figure}
    \centering
    \begin{tikzpicture}
      
      \begin{scope}[shift={(0,0)}]
        \tableaulabel{0}{0}{}
\tableaulabel{0}{1}{}
\tableaulabel{0}{2}{}
\tableaulabel{0}{3}{}
\tableauedgelabeltop{0}{4}{$m_1$}
\tableaubox[fill=red!30]{0}{4}{}
\tableaulabel{1}{0}{}
\tableaubox[fill=blue!30,blue!30]{1}{1}{}
\tableaubox[fill=blue!30,blue!30]{1}{2}{}
\tableaubox[fill=blue!30,blue!30]{1}{3}{}
\tableauedgelabelleft{1}{4}{$k$}
\tableaubox[fill=red!30]{1}{4}{}
\tableaulabel{2}{0}{}
\tableaubox[fill=blue!30,blue!30]{2}{1}{}
\tableaulabel{2}{2}{}
\tableaubox{2}{3}{}
\tableaubox[fill=red!30]{2}{4}{}
\tableaulabel{3}{0}{}
\tableauedgelabelleft{3}{1}{$m_2$}
\tableauedgelabeltop{3}{1}{$l$}
\tableaubox[fill=red!30]{3}{1}{}
\tableaubox[fill=red!30]{3}{2}{}
\tableaubox[fill=red!30]{3}{3}{}
\tableaubox[fill=red!30]{3}{4}{}
\tableaubox{4}{0}{}
\tableaubox{4}{1}{}
\tableaubox{4}{2}{}
\tableaubox{4}{3}{}
\tableaubox{4}{4}{}

      \end{scope}
    
    \end{tikzpicture} \caption{A hook inside $\lambda = (4, 4, 3, 1)$ (blue) and a corresponding hook in $\mathbb{N}^2 \setminus \lambda$ intersecting its boundary edges.}
    \label{fig:hooksInExponentSequenceProp}
  \end{figure}
  
  It remains to show that the formula holds for a box $(i, j) \in \lambda$. The labels $k$ and $l$ are now on the vertical edge of $(i, \lambda_i + 1)$ and the horizontal edge of $(\lambda'_j + 1, j)$, respectively. Both of those boxes, as well as $(\lambda'_j + 1, \lambda_i + 1)$, lie outside $\lambda$, so the first part of the proposition applies to them. Suppose the horizontal and vertical boundary edges of the hook corresponding to the box $(\lambda'_j + 1, \lambda_i + 1)$ are labeled $m_1$ and $m_2$, respectively, as in \figref{fig:hooksInExponentSequenceProp}. Then
  \begin{align*}
    h_\lambda(i, \lambda_i + 1) &= p_\lambda(k) + p_\lambda(m_1)\\
    h_\lambda(\lambda'_j + 1, j) &= p_\lambda(m_2) + p_\lambda(l)\\
    h_\lambda(\lambda'_j + 1, \lambda_i + 1) &= p_\lambda(m_2) + p_\lambda(m_1).
  \end{align*}
  On the other hand,
  \begin{align*}
    h_\lambda(\lambda'_j + 1, \lambda_i + 1) - h_\lambda(i, \lambda_i + 1) - h_\lambda(\lambda'_j + 1, j) &= (\lambda'_j + 1 - i - 1) + (\lambda_i + 1 - j - 1) + 1\\
    &= (\lambda'_j - i) + (\lambda_i - j) + 1\\
    &= h_\lambda(i, j).
  \end{align*}
  Replacing every hook length expression with a sum of entries in the edge power sequence results in
  \begin{align*}
    p_\lambda(m_1) + p_\lambda(m_2) - \left( p_\lambda(k) + p_\lambda(m_1) + p_\lambda(l) + p_\lambda(m_2) \right) &= h_\lambda(i, j)\\
    \Rightarrow h_\lambda(i, j) = -p_\lambda(k) - p_\lambda(l),
  \end{align*}
  as required.
  
  This shows the result for every box in $\lambda$, proving the proposition.
\end{Proof}

In the Young diagram asymptotic to the empty partition (i.e.\ $\mathbb{N}^2$), there are exactly $n$ distinct $n$-hooks: their pivots are the boxes along the off-diagonal from $(n, 1)$ to $(1, n)$. In the more general case of a Young diagram asymptotic to a partition $\lambda$, the number of $n$-hooks depends also on the number of (down-right) $n$-hooks in $\lambda$. The following proposition uses the standard notion of $n$-quotients of a partition; for a thorough reference, see e.g.\ \cite{quotientsWriteUp}. Here, we give a brief treatment of $n$-quotients that integrates well with our notion of edge sign sequence.

\begin{figure}
  \centering
  \begin{tikzpicture}
    
      \begin{scope}[shift={(0,0)}]
        \tableaulabel{0}{0}{}
\tableaubox[fill=gray!30,gray!30]{0}{1}{}
\tableaubox[fill=gray!30,gray!30]{0}{2}{}
\tableaubox[fill=gray!30,gray!30]{0}{3}{}
\tableaubox{0}{4}{}
\tableaulabel{1}{0}{}
\tableaubox[fill=gray!30,gray!30]{1}{1}{}
\tableaubox{1}{2}{}
\tableaubox{1}{3}{}
\tableaubox{1}{4}{}
\tableaulabel{2}{0}{}
\tableaubox[fill=gray!60]{2}{1}{}
\tableaubox{2}{2}{}
\tableaubox{2}{3}{}
\tableaubox{2}{4}{}
\tableaubox{3}{0}{}
\tableaubox[fill=gray!60]{3}{1}{}
\tableaubox{3}{2}{}
\tableaubox{3}{3}{}
\tableaubox{3}{4}{}
\tableaubox[fill=gray!60]{4}{0}{}
\tableaubox[fill=gray!60]{4}{1}{}
\tableaubox{4}{2}{}
\tableaubox{4}{3}{}
\tableaubox{4}{4}{}
\tableaulabel{5.5}{2}{$\lambda = (4, 2, 1)$}
\tableauline[line width=1mm,blue]{-0.5}{3.5}{0.5}{3.5}
\tableauline[line width=1mm,green]{1.5}{0.5}{1.5}{1.5}
\tableauline[line width=1mm,red]{3.5}{-0.5}{4.5}{-0.5}
      \end{scope}
    
      \begin{scope}[shift={(3.5,0)}]
        \tableaubox[fill=gray!30,gray!30]{0}{0}{}
\tableaubox{0}{1}{}
\tableaubox{0}{2}{}
\tableaubox{0}{3}{}
\tableaubox{0}{4}{}
\tableaubox[fill=gray!60]{1}{0}{}
\tableaubox{1}{1}{}
\tableaubox{1}{2}{}
\tableaubox{1}{3}{}
\tableaubox{1}{4}{}
\tableaubox{2}{0}{}
\tableaubox{2}{1}{}
\tableaubox{2}{2}{}
\tableaubox{2}{3}{}
\tableaubox{2}{4}{}
\tableaubox{3}{0}{}
\tableaubox{3}{1}{}
\tableaubox{3}{2}{}
\tableaubox{3}{3}{}
\tableaubox{3}{4}{}
\tableaubox{4}{0}{}
\tableaubox{4}{1}{}
\tableaubox{4}{2}{}
\tableaubox{4}{3}{}
\tableaubox{4}{4}{}
\tableaulabel{5.5}{2}{$\lambda_{4, 3} = (1)$}
\tableauline[line width=1mm,blue]{-0.5}{0.5}{0.5}{0.5}
\tableauline[line width=1mm,green]{0.5}{-0.5}{0.5}{0.5}
\tableauline[line width=1mm,red]{0.5}{-0.5}{1.5}{-0.5}
      \end{scope}
    
  \end{tikzpicture} \caption{A $4$-hook in a skew plane partition of shape $(\emptyset, \emptyset, \lambda)$ and the corresponding $4$-hook in $\lambda$ (left). The corresponding inner and outer corners in the $4$-quotient $\lambda_{4, 3}$ (right).}
  \label{fig:nQuotientExample}
\end{figure}
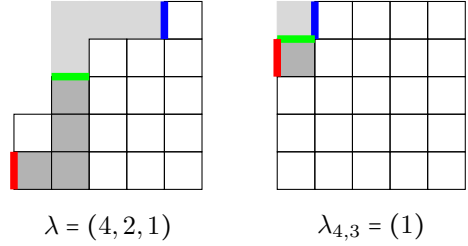

\begin{Def}
  Let $\lambda$ be a partition. Given the edge sign sequence $e_\lambda(m)$ for $m \in \mathbb{Z}$, the subsequence $e_\lambda(nm + i)$ consisting of every $n$th edge sign from $e_\lambda(m)$ corresponds to a distinct partition $\lambda_{n, i}$ for each $i \in \{0, \ldots, n - 1 \}$, called the \textbf{$n$-quotients} of $\lambda$. For example, the object on the right in \figref{fig:nQuotientExample} is the $4$-quotient that selects the red, green, and blue edges. 
\end{Def}

\begin{Prop}
  Let $\lambda$ be a partition and let $n \in \mathbb{N}$. If there are $k$ distinct $n$-hooks in $\lambda$, then there are $n + k$ distinct $n$-hooks in the Young diagram asymptotic to $\lambda$.
\end{Prop}

\begin{Proof}
  Every $n$-hook in $\lambda$ corresponds to a $+$ in $e_\lambda$ followed by a $-$ exactly $n$ terms later, and every $n$-hook in $\mathbb{N}^2 \setminus \lambda$ corresponds to a $+$ edge followed $n$ terms later by a $-$. For example, let $\lambda = (4, 2, 1)$. Then the box $(1, 2) \in \lambda$ has hook length $4$ and corresponds to the green and blue edges in \figref{fig:nQuotientExample}, and the box $(1, 5)$ also has hook length $4$ and corresponds to the green and red edges. There is then a one-to-one correspondence between corners of $\lambda_{n, i}$ for all $i$ and $n$-hooks in $\lambda$, and similarly between corners of $\mathbb{N}^2 \setminus \lambda_{n, i}$ and $n$-hooks in $\mathbb{N}^2 \setminus \lambda$.

  On the other hand, every Young diagram asymptotic to any partition $\mu$ contains exactly one more corner than $\mu$, and so fixing $i$, there is a one-to-one correspondence between corners of $\lambda_{n, i}$ and all the corners of $\mathbb{N}^2 \setminus \lambda_{n, i}$ except for one. Since the two types of corners alternate from bottom-left to top-right, we can make this correspondence explicit by associating each corner of $\lambda_{n, i}$ with the corner of $\mathbb{N}^2 \setminus \lambda_{n, i}$ immediately preceding it in the edge sign sequence. For example, the light and dark gray corners in $\lambda_{4, 3}$ in \figref{fig:nQuotientExample} are in correspondence.

  In total, each of the $k$ distinct $n$-hooks in $\lambda$ corresponds to a corner in $\lambda_{n, i}$ for a unique choice of $i$, which corresponds to a corner in $\mathbb{N}^2 \setminus \lambda_{n, i}$ and therefore an $n$-hook in $\mathbb{N}^2 \setminus \lambda$. However, there is exactly one unmatched corner in $\mathbb{N}^2 \setminus \lambda_{n, i}$ for each $i$, each of which contributes another $n$-hook in $\mathbb{N}^2 \setminus \lambda$, resulting in a total of $n + k$ distinct $n$ hooks in $\mathbb{N}^2 \setminus \lambda$.
\end{Proof}

We aim to bijectivize the one-leg PT--DT correspondence \eqref{eq:PTDT1} on the level of hooks; that is, to decompose every object involved into a hook-length-weighted tableau and place the resulting entries in bijection with one another. To that end, we define a map associating the $n + k$ distinct $n$-hooks of the asymptotic Young diagram $\mathbb{N}^2 \setminus \lambda$ to the $n$ distinct $n$-hooks of the asymptotic Young diagram $\mathbb{N}^2$ and the $k$ distinct $n$-hooks of the Young diagram $\lambda$.

\begin{Def}
  For a box $b \in \mathbb{N}^2 \setminus \lambda$ with hook length $h(b) = n$, let $c$ be the unique corner in a quotient $\mathbb{N}^2 \setminus \lambda_{n, i}$ that corresponds to $b$, and define $\varphi_\lambda(b)$ in the following manner. If $c$ is the upper-right-most corner in $\mathbb{N}^2 \setminus \lambda_{n, i}$, then $\varphi_\lambda(b)$ is the unique box in $\mathbb{N}^2$ with hook length $n$ that corresponds to the corner in the $i$th $n$-quotient of $\mathbb{N}^2$ (i.e.\ the $i$th box in the $n$th off-diagonal of $\mathbb{N}^2$ when reading from bottom-left to top-right). Otherwise, $\varphi_\lambda(b)$ is the unique box in $\lambda$ corresponding to the corner in $\lambda_{n, i}$ immediately following $c$ in the alternating sequence of corners inside and outside of $\lambda_{n, i}$ when reading from bottom-left to top-right.
  \label{def:quotientPartition}
\end{Def}

This bijection is, perhaps unsurprisingly, best explained by example.

\begin{figure}
  \centering
  \begin{tikzpicture}
    
      \begin{scope}[shift={(0,0)}]
        \tableaulabel{0}{0}{}
\tableaulabel{0}{1}{}
\tableaulabel{0}{2}{}
\tableaubox{0}{3}{}
\tableaubox{0}{4}{}
\tableaubox[fill=violet!60]{0}{5}{}
\tableaulabel{1}{0}{}
\tableaulabel{1}{1}{}
\tableaulabel{1}{2}{}
\tableaubox{1}{3}{}
\tableaubox[fill=blue!60]{1}{4}{}
\tableaubox{1}{5}{}
\tableaubox{2}{0}{}
\tableaubox{2}{1}{}
\tableaubox[fill=green!60]{2}{2}{}
\tableaubox{2}{3}{}
\tableaubox{2}{4}{}
\tableaubox{2}{5}{}
\tableaubox{3}{0}{}
\tableaubox[fill=orange!60]{3}{1}{}
\tableaubox{3}{2}{}
\tableaubox{3}{3}{}
\tableaubox{3}{4}{}
\tableaubox{3}{5}{}
\tableaubox[fill=red!60]{4}{0}{}
\tableaubox{4}{1}{}
\tableaubox{4}{2}{}
\tableaubox{4}{3}{}
\tableaubox{4}{4}{}
\tableaubox{4}{5}{}
\tableaubox{5}{0}{}
\tableaubox{5}{1}{}
\tableaubox{5}{2}{}
\tableaubox{5}{3}{}
\tableaubox{5}{4}{}
\tableaubox{5}{5}{}
\tableaulabel{2.5}{6}{$\cdots$}
\tableaulabel{6}{2.5}{$\vdots$}

      \end{scope}
    
      \begin{scope}[shift={(4.5,-0.75)}]
        \tableaubox[fill=green!60]{0}{0}{}
\tableaubox{0}{1}{}
\tableaubox{0}{2}{}
\tableaubox{1}{0}{}
\tableaubox{1}{1}{}
\tableaubox{1}{2}{}
\tableaubox{2}{0}{}
\tableaubox{2}{1}{}
\tableaubox{2}{2}{}
\tableaulabel{1}{3}{$\cdots$}
\tableaulabel{3}{1}{$\vdots$}

      \end{scope}
    
      \begin{scope}[shift={(7.5,-0.75)}]
        \tableaubox[red!30,fill=red!30]{0}{0}{}
\tableaubox[fill=blue!60]{0}{1}{}
\tableaubox{0}{2}{}
\tableaubox[fill=red!60]{1}{0}{}
\tableaubox{1}{1}{}
\tableaubox{1}{2}{}
\tableaubox{2}{0}{}
\tableaubox{2}{1}{}
\tableaubox{2}{2}{}
\tableaulabel{1}{3}{$\cdots$}
\tableaulabel{3}{1}{$\vdots$}

      \end{scope}
    
      \begin{scope}[shift={(10.5,-0.75)}]
        \tableaubox[orange!30,fill=orange!30]{0}{0}{}
\tableaubox[fill=violet!60]{0}{1}{}
\tableaubox{0}{2}{}
\tableaubox[fill=orange!60]{1}{0}{}
\tableaubox{1}{1}{}
\tableaubox{1}{2}{}
\tableaubox{2}{0}{}
\tableaubox{2}{1}{}
\tableaubox{2}{2}{}
\tableaulabel{1}{3}{$\cdots$}
\tableaulabel{3}{1}{$\vdots$}

      \end{scope}
    
  \end{tikzpicture}
  \begin{tikzpicture}
    
      \begin{scope}[shift={(0,-0.5)}]
        \tableaubox{0}{0}{}
\tableaubox[fill=orange!30]{0}{1}{}
\tableaubox{0}{2}{}
\tableaubox[fill=red!30]{1}{0}{}
\tableaubox{1}{1}{}
\tableaubox{1}{2}{}

      \end{scope}
    
      \begin{scope}[shift={(2.5,0)}]
        \tableaubox{0}{0}{}
\tableaubox{0}{1}{}
\tableaubox[fill=violet!60]{0}{2}{}
\tableaubox{0}{3}{}
\tableaubox{1}{0}{}
\tableaubox[fill=blue!60]{1}{1}{}
\tableaubox{1}{2}{}
\tableaubox{1}{3}{}
\tableaubox[fill=green!60]{2}{0}{}
\tableaubox{2}{1}{}
\tableaubox{2}{2}{}
\tableaubox{2}{3}{}
\tableaubox{3}{0}{}
\tableaubox{3}{1}{}
\tableaubox{3}{2}{}
\tableaubox{3}{3}{}
\tableaulabel{1.5}{4}{$\cdots$}
\tableaulabel{4}{1.5}{$\vdots$}

      \end{scope}
    
  \end{tikzpicture} \caption{Pivots of corresponding $3$-hooks in the Young diagram asymptotic to $\lambda = (3, 3)$ and its quotients. In reading order: $\mathbb{N}^2 \setminus \lambda$, $\mathbb{N}^2 \setminus \lambda_{3, 0}$, $\mathbb{N}^2 \setminus \lambda_{3, 1}$, $\mathbb{N}^2 \setminus \lambda_{3, 2}$, $\lambda$, $\mathbb{N}^2$.}
  \label{fig:nQuotientPivots}
\end{figure}
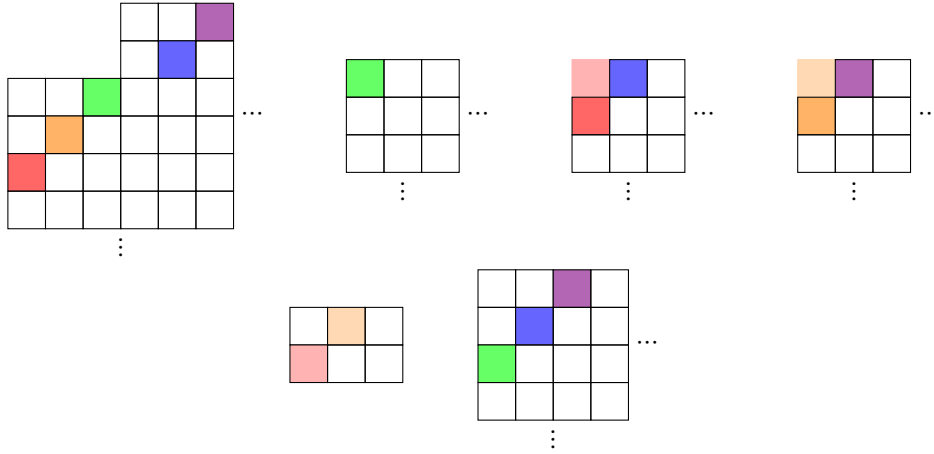

\begin{Ex}
  Let $\lambda = (3, 3)$. The Young diagram asymptotic to $\lambda$ has $5$ distinct $3$-hooks; their pivots are shaded in \figref{fig:nQuotientPivots}.

  The three $3$-quotients of $\lambda$ are $\lambda_{3, 0} = \emptyset$, $\lambda_{3, 1} = (1)$, and $\lambda_{3, 2} = (1)$. Each colored box in $\mathbb{N}^2 \setminus \lambda$ corresponds to a unique corner in $\mathbb{N}^2 \setminus \lambda_{3, i}$ for some $i$; for example, the box $(4, 2)$ is colored orange, and the two edges where its hook intersects the boundary of the diagram are present and adjacent in $\lambda_{3, 2}$ (and therefore no other $3$-quotient). It therefore corresponds to the orange box $(2, 1)$ in $\mathbb{N}^2 \setminus \lambda_{3, 2}$. Similarly, the box $(1, 6) \in \mathbb{N}^2 \setminus \lambda$ is colored purple and corresponds to the purple box $(1, 2) \in \mathbb{N}^2 \setminus \lambda_{3, 2}$.

  These two corners $(2, 1)$ and $(1, 2)$ in $\mathbb{N}^2 \setminus \lambda_{3, 2}$ necessarily have exactly one corner of $\lambda_{3, 2}$ (i.e.\ $(1, 1)$) occurring between them when reading the edge sequence from bottom-left to top right. We therefore associate $(2, 1)$ with $(1, 1)$ and leave $(1, 2)$ temporarily unassociated.

  Repeating this calculation for the remaining quotients $\lambda_{3, 0}$ and $\lambda_{3, 1}$ associates the red box $(2, 1) \in \mathbb{N}^2 \setminus \lambda_{3, 1}$ with $(1, 1) \in \lambda_{3, 1}$, and leaves the green, blue, and purple boxes in the quotients unassociated. The purple box corresponds to a corner in the quotient $\lambda_{3, 2}$; in the asymptotic Young diagram $\mathbb{N}^2 = \mathbb{N}^2 \setminus \emptyset$, the unique box with hook length $3$ that corresponds to a corner in the quotient $\mathbb{N}^2 \setminus \emptyset_{3, 2}$ is $(1, 3)$, so we define $\varphi_\lambda(1, 6) = (1, 3)$, and similarly for the green and blue boxes. By the same logic, the red and orange corners in $\lambda_{1, 3}$ and $\lambda_{2, 3}$ correspond to the boxes $(2, 1)$ and $(1, 2)$ in $\lambda$. In total, the map $\varphi_\lambda$ has the following effect on hook-length-$3$ boxes in $\mathbb{N}^2 \setminus \lambda$:
  $$
    \varphi_\lambda(5, 1) = (2, 1) \in \lambda \hspace{.5in} \varphi_\lambda(4, 2) = (1, 2) \in \lambda
  $$%
  $$
    \varphi_\lambda(3, 3) = (3, 1) \in \mathbb{N}^2 \hspace{.5in} \varphi_\lambda(2, 5) = (2, 2) \in \mathbb{N}^2 \hspace{.5in} \varphi_\lambda(1, 6) = (1, 3) \in \mathbb{N}^2.
  $$
\end{Ex}

At last, we are prepared to bijectivize the one-leg PT--DT correspondence \eqref{eq:PTDT1}. The proof amounts to a generalization of \thmref{thm:macMachonsFunctionByToggling} to SPPs that is equivalent to the Pak--Sulzgruber algorithm \cite{pak, sulzgruber}, combined with an application of the $\varphi_\lambda$ map.

\begin{Thm}
  \label{thm:oneLegDecomposition}
  Let $\lambda \subset \mathbb{N}^2$ be a Young diagram. Then there is a weight-preserving bijection between SPPs $\sigma$ of shape $(\emptyset, \emptyset, \lambda)$ and pairs $(\rho, \pi)$ of RPPs $\rho$ of shape $(\emptyset, \emptyset, \lambda)$ and plane partitions $\pi$, where $|\sigma| = |\rho| + |\pi|$.
\end{Thm}

\begin{Proof}
  Recall that $V_{(\emptyset, \emptyset, \lambda)}(q)$ can be expressed in terms of vertex operators as
  $$
    V_{(\emptyset, \emptyset, \lambda)}(q) = \leftbra{\emptyset}\,\prod_{n \in \mathbb{Z}} \Gamma_{e(n)}\left(q^{p(n)}\right)\,\rightket{\emptyset},
  $$
  where $e(n) = e_\lambda(n)$ and $p(n) = p_\lambda(n)$ are the edge sign and edge power sequences, respectively. Analogous to the proof of \thmref{thm:macMachonsFunctionByToggling}, we follow the algebraic proof of commuting every $\Gamma_-$ to the right and every $\Gamma_+$ to the left, while leaving the order of same-sign operators unchanged (i.e.\ only commuting operators of opposite sign). By \lemref{lem:bijectivizingToggles}, we may associate a bijection to each commutation we perform, where each one toggles a diagonal in the diagram corresponding to a corner. Moreover, each produces a factor of $(1 - q^{h(\square)})^{-1}$ by \thmref{thm:hookLengthsInAsymptoticDiagrams} and the same logic as in the proof of \thmref{thm:macMachonsFunctionByToggling}. Since any given SPP contains only finitely many nonzero entries, the composition of these bijections is well-defined. Algebraically, we have expressed
  $$
    V_{(\emptyset, \emptyset, \lambda)}(q) = \prod_{\square \in \mathbb{N}^2 \setminus \lambda} \frac{1}{1 - q^{h(\square)}},
  $$
  and combinatorially, there is a bijection between SPPs and hook-length-weighted tableaux of the same shape, given by toggling diagonals until the SPP is empty. We denote it $\tau$, since it specializes to the bijection $\tau$ given in \thmref{thm:macMachonsFunctionByToggling} when $\lambda = \emptyset$.

  We may now complete the proof of the theorem. Given an SPP $\sigma$ of shape $(\emptyset, \emptyset, \lambda)$, we apply the bijection $\tau$ to create a hook-length-weighted tableau $\tau(\sigma)$. We then apply the map $\varphi_\lambda$ from \defref{def:quotientPartition} that associates the cells of $\tau(\sigma)$ with the cells of two tableaux $R$ and $P$ of shapes $\lambda$ and $\mathbb{N}^2$, respectively, and preserves hook length.

  What remains is to convert these two hook-length-weighted tableaux back into an RPP of shape $(\emptyset, \emptyset, \lambda)$ and a plane partition. For the latter, we apply $\tau^{-1}$ to produce a plane partition $\pi = \tau^{-1}(P)$. For the former, slightly more care must be taken, since the hooks in $R$ extend down and right rather than up and left. However, we may rotate $R$ by $180^\circ$ to express it as an SPP while preserving cells' hook lengths; it is exactly this latter presentation of the bijection that is the Pak--Sulzgruber algorithm \cite{pak, sulzgruber}. We rotate $R$, apply $\tau^{-1}$ to this new SPP, and rotate back, producing an RPP $\rho$ of shape $(\emptyset, \emptyset, \lambda)$ that satisfies $|\sigma| = |\rho| + |\pi|$.

  As in the proof of \thmref{thm:hookLengthsInAsymptoticDiagrams}, one useful property of this bijectivization is that it cleanly demonstrates the independence of the algorithm from the order in which we toggle diagonals --- each $\Gamma$ commutation produces a commutator that factors directly out of the vertex operator product and does not affect any other operators. This independence is also proven directly in \cite{pak}.
\end{Proof}

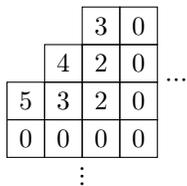
\begin{figure}
  \centering
  \begin{tikzpicture}
    
      \begin{scope}[shift={(0,0)}]
        \tableaulabel{0}{0}{}
\tableaulabel{0}{1}{}
\tableaubox{0}{2}{3}
\tableaubox{0}{3}{0}
\tableaulabel{1}{0}{}
\tableaubox{1}{1}{4}
\tableaubox{1}{2}{2}
\tableaubox{1}{3}{0}
\tableaubox{2}{0}{5}
\tableaubox{2}{1}{3}
\tableaubox{2}{2}{2}
\tableaubox{2}{3}{0}
\tableaubox{3}{0}{0}
\tableaubox{3}{1}{0}
\tableaubox{3}{2}{0}
\tableaubox{3}{3}{0}
\tableaulabel{1.5}{4}{$\cdots$}
\tableaulabel{4}{1.5}{$\vdots$}

      \end{scope}
    
  \end{tikzpicture} \caption{A skew plane partition $\sigma$ of shape $(\emptyset, \emptyset, (2, 1))$.}
  \label{fig:asymptoticPlanePartitionExample}
\end{figure}

\begin{figure}
  \centering
  \begin{tikzpicture}
    
      \begin{scope}[shift={(0,0)}]
        \tableaulabel{0}{0}{}
\tableaulabel{0}{1}{}
\tableaubox{0}{2}{3}
\tableaubox{0}{3}{0}
\tableaulabel{1}{0}{}
\tableaubox{1}{1}{4}
\tableaubox{1}{2}{2}
\tableaubox{1}{3}{0}
\tableaubox{2}{0}{5}
\tableaubox{2}{1}{3}
\tableaubox{2}{2}{2}
\tableaubox{2}{3}{0}
\tableaubox{3}{0}{0}
\tableaubox{3}{1}{0}
\tableaubox{3}{2}{0}
\tableaubox{3}{3}{0}

      \tableaulabel{1.5}{4}{$\mapsto$}
    \tableaudraw[line width=1mm]{(-0.5, 3.5) -- (-0.5, 1.5) -- (0.5, 1.5) -- (0.5, 0.5) -- (1.5, 0.5) -- (1.5, -0.5) -- (2.5, -0.5) -- (2.5, -0.5) -- (3.5, -0.5)}
      \end{scope}
    
      \begin{scope}[shift={(2.5,0)}]
        \tableaulabel{0}{0}{}
\tableaulabel{0}{1}{}
\tableaubox{0}{2}{1}
\tableaubox{0}{3}{0}
\tableaulabel{1}{0}{}
\tableaubox{1}{1}{4}
\tableaubox{1}{2}{2}
\tableaubox{1}{3}{0}
\tableaubox{2}{0}{5}
\tableaubox{2}{1}{3}
\tableaubox{2}{2}{2}
\tableaubox{2}{3}{0}
\tableaubox{3}{0}{0}
\tableaubox{3}{1}{0}
\tableaubox{3}{2}{0}
\tableaubox{3}{3}{0}

      \tableaulabel{1.5}{4}{$\mapsto$}
      \tableaulabel{1.5}{5}{$\cdots$}
      \tableaulabel{1.5}{6}{$\mapsto$}
    \tableaudraw[line width=1mm]{(-0.5, 3.5) -- (-0.5, 2.5) -- (0.5, 2.5) -- (0.5, 0.5) -- (1.5, 0.5) -- (1.5, -0.5) -- (2.5, -0.5) -- (2.5, -0.5) -- (3.5, -0.5)}
      \end{scope}
    
      \begin{scope}[shift={(6,0)}]
        \tableaulabel{0}{0}{}
\tableaulabel{0}{1}{}
\tableaubox{0}{2}{1}
\tableaubox{0}{3}{0}
\tableaulabel{1}{0}{}
\tableaubox{1}{1}{4}
\tableaubox{1}{2}{2}
\tableaubox{1}{3}{0}
\tableaubox{2}{0}{5}
\tableaubox{2}{1}{3}
\tableaubox{2}{2}{2}
\tableaubox{2}{3}{0}
\tableaubox{3}{0}{0}
\tableaubox{3}{1}{0}
\tableaubox{3}{2}{0}
\tableaubox{3}{3}{0}

      \tableaulabel{1.5}{4}{$\mapsto$}
    \tableaudraw[line width=1mm]{(0.5, 3.5) -- (0.5, 0.5) -- (1.5, 0.5) -- (1.5, -0.5) -- (2.5, -0.5) -- (2.5, -0.5) -- (3.5, -0.5)}
      \end{scope}
    
      \begin{scope}[shift={(8.5,0)}]
        \tableaulabel{0}{0}{}
\tableaulabel{0}{1}{}
\tableaubox{0}{2}{1}
\tableaubox{0}{3}{0}
\tableaulabel{1}{0}{}
\tableaubox{1}{1}{1}
\tableaubox{1}{2}{2}
\tableaubox{1}{3}{0}
\tableaubox{2}{0}{5}
\tableaubox{2}{1}{3}
\tableaubox{2}{2}{0}
\tableaubox{2}{3}{0}
\tableaubox{3}{0}{0}
\tableaubox{3}{1}{0}
\tableaubox{3}{2}{0}
\tableaubox{3}{3}{0}

      \tableaulabel{1.5}{4}{$\mapsto$}
    \tableaudraw[line width=1mm]{(0.5, 3.5) -- (0.5, 1.5) -- (1.5, 1.5) -- (1.5, -0.5) -- (2.5, -0.5) -- (2.5, -0.5) -- (3.5, -0.5)}
      \end{scope}
    
      \begin{scope}[shift={(11,0)}]
        \tableaulabel{0}{0}{}
\tableaulabel{0}{1}{}
\tableaubox{0}{2}{1}
\tableaubox{0}{3}{0}
\tableaulabel{1}{0}{}
\tableaubox{1}{1}{1}
\tableaubox{1}{2}{2}
\tableaubox{1}{3}{0}
\tableaubox{2}{0}{5}
\tableaubox{2}{1}{3}
\tableaubox{2}{2}{0}
\tableaubox{2}{3}{0}
\tableaubox{3}{0}{0}
\tableaubox{3}{1}{0}
\tableaubox{3}{2}{0}
\tableaubox{3}{3}{0}

      \tableaulabel{1.5}{4}{$\mapsto$}
      \tableaulabel{1.5}{5}{$\cdots$}
      \tableaulabel{1.5}{6}{$\mapsto$}
    \tableaudraw[line width=1mm]{(0.5, 3.5) -- (0.5, 2.5) -- (1.5, 2.5) -- (1.5, -0.5) -- (2.5, -0.5) -- (2.5, -0.5) -- (3.5, -0.5)}
      \end{scope}
    
      \begin{scope}[shift={(14.5,0)}]
        \tableaulabel{0}{0}{}
\tableaulabel{0}{1}{}
\tableaubox{0}{2}{1}
\tableaubox{0}{3}{0}
\tableaulabel{1}{0}{}
\tableaubox{1}{1}{1}
\tableaubox{1}{2}{2}
\tableaubox{1}{3}{0}
\tableaubox{2}{0}{5}
\tableaubox{2}{1}{3}
\tableaubox{2}{2}{0}
\tableaubox{2}{3}{0}
\tableaubox{3}{0}{0}
\tableaubox{3}{1}{0}
\tableaubox{3}{2}{0}
\tableaubox{3}{3}{0}

      \tableaulabel{1.5}{4}{$\mapsto$}
    \tableaudraw[line width=1mm]{(1.5, 3.5) -- (1.5, -0.5) -- (2.5, -0.5) -- (2.5, -0.5) -- (3.5, -0.5)}
      \end{scope}
    
  \end{tikzpicture}

  \vspace{18pt}

  \begin{tikzpicture}
    
      \begin{scope}[shift={(0,0)}]
        \tableaulabel{0}{0}{}
\tableaulabel{0}{1}{}
\tableaubox{0}{2}{1}
\tableaubox{0}{3}{0}
\tableaulabel{1}{0}{}
\tableaubox{1}{1}{1}
\tableaubox{1}{2}{2}
\tableaubox{1}{3}{0}
\tableaubox{2}{0}{2}
\tableaubox{2}{1}{3}
\tableaubox{2}{2}{0}
\tableaubox{2}{3}{0}
\tableaubox{3}{0}{0}
\tableaubox{3}{1}{0}
\tableaubox{3}{2}{0}
\tableaubox{3}{3}{0}

      \tableaulabel{1.5}{4}{$\mapsto$}
    \tableaudraw[line width=1mm]{(1.5, 3.5) -- (1.5, 0.5) -- (2.5, 0.5) -- (2.5, -0.5) -- (3.5, -0.5)}
      \end{scope}
    
      \begin{scope}[shift={(2.5,0)}]
        \tableaulabel{0}{0}{}
\tableaulabel{0}{1}{}
\tableaubox{0}{2}{1}
\tableaubox{0}{3}{0}
\tableaulabel{1}{0}{}
\tableaubox{1}{1}{1}
\tableaubox{1}{2}{2}
\tableaubox{1}{3}{0}
\tableaubox{2}{0}{2}
\tableaubox{2}{1}{3}
\tableaubox{2}{2}{0}
\tableaubox{2}{3}{0}
\tableaubox{3}{0}{0}
\tableaubox{3}{1}{0}
\tableaubox{3}{2}{0}
\tableaubox{3}{3}{0}

      \tableaulabel{1.5}{4}{$\mapsto$}
      \tableaulabel{1.5}{5}{$\cdots$}
      \tableaulabel{1.5}{6}{$\mapsto$}
    \tableaudraw[line width=1mm]{(1.5, 3.5) -- (1.5, 1.5) -- (2.5, 1.5) -- (2.5, -0.5) -- (3.5, -0.5)}
      \end{scope}
    
      \begin{scope}[shift={(6,0)}]
        \tableaulabel{0}{0}{}
\tableaulabel{0}{1}{}
\tableaubox{0}{2}{1}
\tableaubox{0}{3}{0}
\tableaulabel{1}{0}{}
\tableaubox{1}{1}{1}
\tableaubox{1}{2}{2}
\tableaubox{1}{3}{0}
\tableaubox{2}{0}{2}
\tableaubox{2}{1}{3}
\tableaubox{2}{2}{0}
\tableaubox{2}{3}{0}
\tableaubox{3}{0}{0}
\tableaubox{3}{1}{0}
\tableaubox{3}{2}{0}
\tableaubox{3}{3}{0}

      \end{scope}
    
  \end{tikzpicture} \caption{Decomposing a one-leg SPP into a hook-length-weighted tableau. At each step, we choose a corner of the asymptotic Young diagram (here, we choose them in lexicographic order) and toggle its diagonal.}
  \label{fig:oneLegSPPDecomposition}
\end{figure}
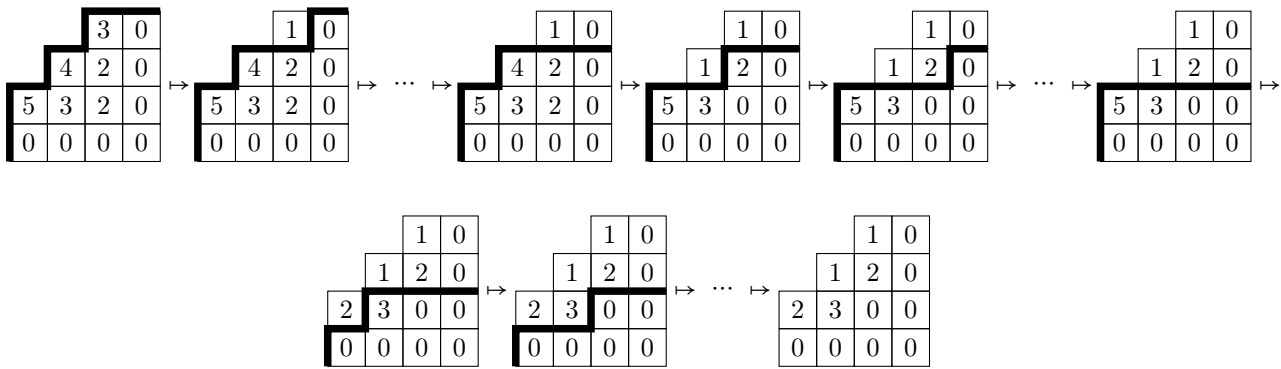

\begin{figure}
  \centering
  \begin{tikzpicture}
    
      \begin{scope}[shift={(0,0)}]
        \tableaulabel{0}{0}{}
\tableaulabel{0}{1}{}
\tableaubox{0}{2}{1}
\tableaubox{0}{3}{0}
\tableaulabel{1}{0}{}
\tableaubox{1}{1}{1}
\tableaubox{1}{2}{2}
\tableaubox{1}{3}{0}
\tableaubox{2}{0}{2}
\tableaubox{2}{1}{3}
\tableaubox{2}{2}{0}
\tableaubox{2}{3}{0}
\tableaubox{3}{0}{0}
\tableaubox{3}{1}{0}
\tableaubox{3}{2}{0}
\tableaubox{3}{3}{0}
\tableaulabel{1.5}{4}{$\cdots$}
\tableaulabel{4}{1.5}{$\vdots$}

      \end{scope}
    
      \begin{scope}[shift={(3.5,-0.5)}]
        \tableaubox{0}{0}{0}
\tableaubox{0}{1}{1}
\tableaubox{1}{0}{2}

      \end{scope}
    
      \begin{scope}[shift={(5.5,0)}]
        \tableaubox{0}{0}{1}
\tableaubox{0}{1}{0}
\tableaubox{0}{2}{0}
\tableaubox{0}{3}{0}
\tableaubox{1}{0}{0}
\tableaubox{1}{1}{2}
\tableaubox{1}{2}{0}
\tableaubox{1}{3}{0}
\tableaubox{2}{0}{3}
\tableaubox{2}{1}{0}
\tableaubox{2}{2}{0}
\tableaubox{2}{3}{0}
\tableaubox{3}{0}{0}
\tableaubox{3}{1}{0}
\tableaubox{3}{2}{0}
\tableaubox{3}{3}{0}
\tableaulabel{1.5}{4}{$\cdots$}
\tableaulabel{4}{1.5}{$\vdots$}

      \end{scope}
    
  \end{tikzpicture} \caption{Rearranging the entries of the tableau $T(\sigma)$ (left) into a tableau $T(\rho)$ (center) of shape $\lambda$ and a tableau $T(\pi)$ (right) of shape $\mathbb{N}^2$.}
  \label{fig:tableauDecomposition}
\end{figure}
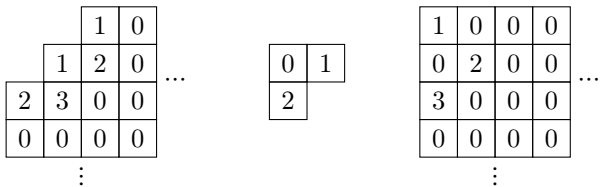

\begin{figure}
  \centering
  \begin{tikzpicture}
    
      \begin{scope}[shift={(0,-0.5)}]
        \tableaubox{0}{0}{0}
\tableaubox{0}{1}{1}
\tableaubox{1}{0}{2}

      \end{scope}
    
      \begin{scope}[shift={(2,0)}]
        \tableaubox{0}{0}{4}
\tableaubox{0}{1}{2}
\tableaubox{0}{2}{0}
\tableaubox{0}{3}{0}
\tableaubox{1}{0}{3}
\tableaubox{1}{1}{2}
\tableaubox{1}{2}{0}
\tableaubox{1}{3}{0}
\tableaubox{2}{0}{3}
\tableaubox{2}{1}{2}
\tableaubox{2}{2}{0}
\tableaubox{2}{3}{0}
\tableaubox{3}{0}{0}
\tableaubox{3}{1}{0}
\tableaubox{3}{2}{0}
\tableaubox{3}{3}{0}
\tableaulabel{1.5}{4}{$\cdots$}
\tableaulabel{4}{1.5}{$\vdots$}

      \end{scope}
    
  \end{tikzpicture} \caption{The untoggled RPP $\rho = T^{-1}(T(\rho))$ (left) and plane partition $\pi = T^{-1}(T(\pi))$.}
  \label{fig:untoggling}
\end{figure}
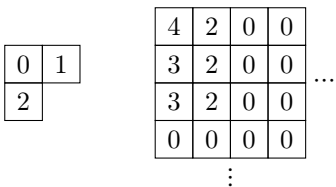

\begin{Ex}
  Let $\lambda = (2, 1)$ and let $\sigma$ be the SPP of shape to $(\emptyset, \emptyset, \lambda)$ in \figref{fig:asymptoticPlanePartitionExample}. To convert $\sigma$ into its corresponding tableau $\tau(\sigma)$, we iteratively toggle diagonals beginning with a corner; one consistent way to accomplish this is to choose those corners in lexicographic order. With this choice of order, we produce the sequence of diagrams in \figref{fig:oneLegSPPDecomposition}. For compactness, we omit the ellipses on the right and bottom of the diagrams, and since we record each entry in the tableau in exactly the location we remove it from the plane partition, we superimpose the two objects, following the convention of \cite{pak} by separating them with a thicker border.

  We now apply the map $\varphi_\lambda$ to $\tau(\sigma)$, producing the tableaux $R$ and $P$ in \figref{fig:tableauDecomposition}. Finally, we apply $\tau^{-1}$ to each of these two tableaux (taking care to rotate $R$ by $180^\circ$ before untoggling so that the usual plane partition inequalities hold), producing a final RPP $\rho$ and plane partition $\pi$ in \figref{fig:untoggling}. Checking the weights, we have $|\rho| + |\pi| = 3 + 16 = 19 = |\sigma|$, as expected.
\end{Ex}

What is new in this section is neither the bijection $\tau$ nor the corresponding algebraic proof using vertex operators, but the connection between the two of them, along with the use of $n$-quotients to bijectivize the one-leg PT--DT correspondence. While the details are more difficult in future sections, the fundamental approach remains similar.

\section{Two-Leg Objects} \label{sec:twoLegObjects}

Both the generating function identity and the bijection of \thmref{thm:oneLegDecomposition} hold in a different and in fact more general setting. We begin by generalizing the definitions of skew and reverse plane partition before stating and proving the new result.

\begin{Def}
  Let $\lambda$ and $\mu$ be partitions. A \textbf{skew plane partition} (SPP) with \textbf{shape $(\lambda, \mu, \emptyset)$} is a plane partition $\sigma$ for which $\sigma(i, j) \geq \max\{\lambda_j, \mu_i\}$ for all $i, j \in \mathbb{N}$, and only finitely many $\sigma(i, j)$ satisfy $\sigma(i, j) > \max\{\lambda_j, \mu_i\}$. The generating function $V_{(\lambda, \mu, \emptyset)}(q)$ is given by
  \begin{align}
    \label{eq:twoLegSppVertexOperators}
    V_{(\lambda, \mu, \emptyset)}(q) = \leftbra{\lambda}\,\prod_{n = -\infty}^{0} \Gamma_-\left(q^{(-n + 1) / 2}\right)\prod_{n = 0}^\infty \Gamma_+\left(q^{(n + 1) / 2}\right)\,\rightket{\mu},
  \end{align}
  where $q$ marks the weight \cite[Equation 3.19]{vertexOperators}. If $V_0(\lambda, \mu, \emptyset)$ denotes the minimal exponent of $q$ in $V_{(\lambda, \mu, \emptyset)}(q)$, i.e.\ the weight of the minimal configuration, then the weight $|\sigma|$ of a general SPP $\sigma$ with shape $(\lambda, \mu, \emptyset)$ is given by
  $$
    |\sigma| = V_0(\lambda, \mu, \emptyset) + \sum_{i, j \in \mathbb{N}} \left( \sigma(i, j) - \max\{\lambda_j, \mu_i\} \right).
  $$
  We call these objects \textbf{two-leg} SPPs in the terminology of \cite{PT}.
\end{Def}

\begin{figure}
  \begin{tabu}{X[m,c]X[m,c]}
    \begin{tikzpicture}
      
      \begin{scope}[shift={(0,0)}]
        \tableaubox{0}{0}{5}
\tableaubox{0}{1}{4}
\tableaubox{0}{2}{3}
\tableaubox{0}{3}{3}
\tableaubox{0}{4}{3}
\tableaubox{1}{0}{5}
\tableaubox{1}{1}{3}
\tableaubox{1}{2}{2}
\tableaubox{1}{3}{1}
\tableaubox{1}{4}{1}
\tableaubox{2}{0}{3}
\tableaubox{2}{1}{2}
\tableaubox{2}{2}{1}
\tableaubox{2}{3}{0}
\tableaubox{2}{4}{0}
\tableaubox{3}{0}{2}
\tableaubox{3}{1}{2}
\tableaubox{3}{2}{0}
\tableaubox{3}{3}{0}
\tableaubox{3}{4}{0}
\tableaubox{4}{0}{2}
\tableaubox{4}{1}{2}
\tableaubox{4}{2}{0}
\tableaubox{4}{3}{0}
\tableaubox{4}{4}{0}
\tableaulabel{2}{5}{$\cdots$}
\tableaulabel{5}{2}{$\vdots$}

      \end{scope}
    
    \end{tikzpicture} & \includegraphics[height=2.5in]{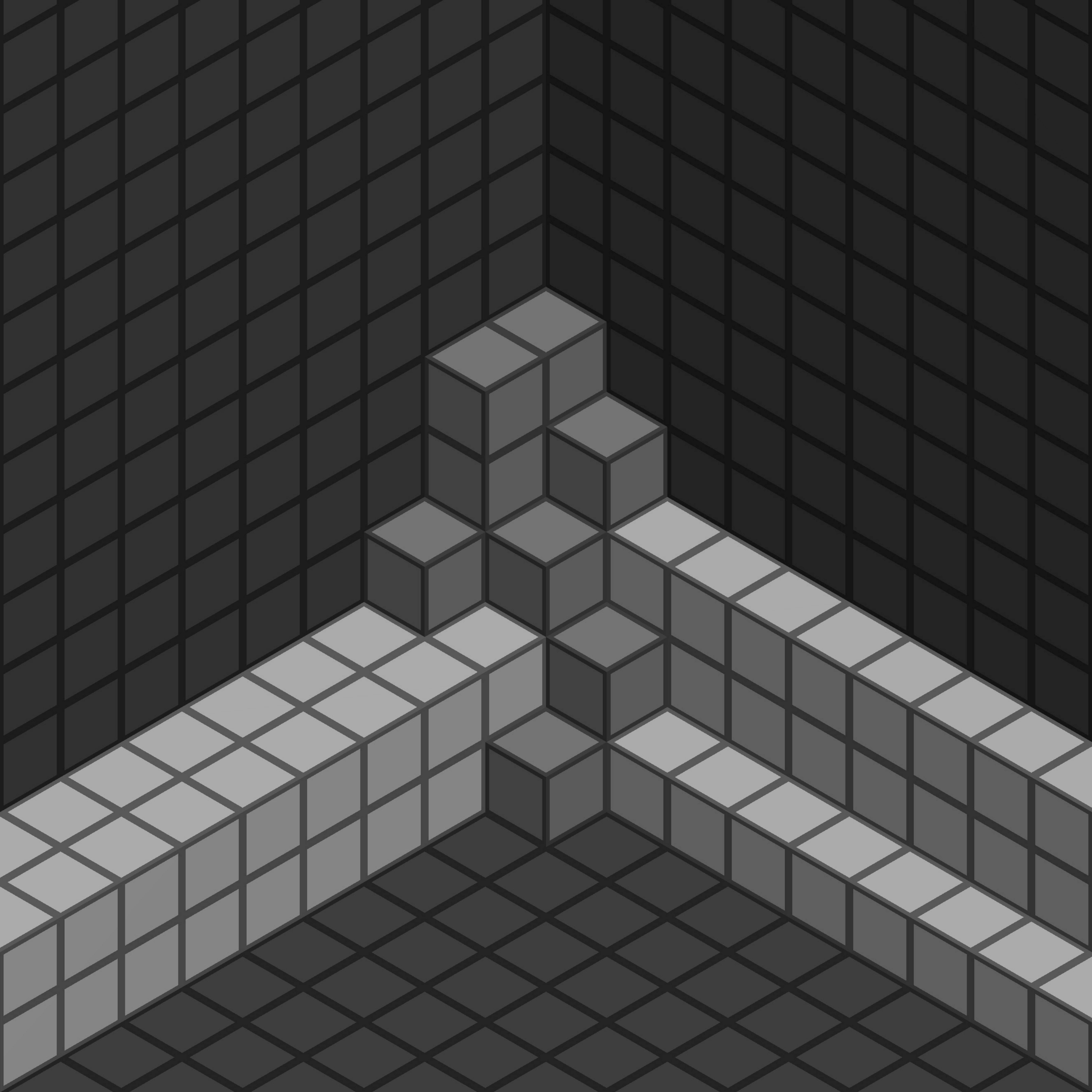}
  \end{tabu} \caption{A two-leg skew plane partition of shape $((2, 2), (3, 1), \emptyset)$, visualized both as a grid of numbers and a stack of 10 weight-contributing blocks (dark gray) on top of a non-removable ``tray'' of blocks that contributes weight 1, as measured by the vertex operator expression \eqref{eq:twoLegSppVertexOperators}.}
  \label{fig:twoLegSPP}
\end{figure}

An example of an SPP of weight $11$ and shape $(\lambda, \mu, \emptyset)$ for $\lambda = (2, 2)$ and $\mu = (3, 1)$ is given in \figref{fig:twoLegSPP}. We remark that our convention differs slightly from \cite{vertexOperators}: their definition uses the conjugate partition $\lambda'$ instead of $\lambda$ itself, in order to match their convention for the fully general three-leg case. In \cite[Equation 3.20]{vertexOperators}, the authors also give an expression for three-leg (and thereby two-leg) objects in terms of Schur functions.

\begin{figure}
  \begin{tabu}{X[m,c]X[m,c]}
    \begin{tikzpicture}
        
      \begin{scope}[shift={(0,0)}]
        \tableaubox{0}{0}{0}
\tableaubox{0}{1}{2}
\tableaubox{0}{2}{2}
\tableaubox{0}{3}{3}
\tableaubox{1}{0}{2}
\tableaubox{1}{1}{3}
\tableaubox{1}{2}{4}
\tableaubox{2}{0}{3}

      \end{scope}
    
    \end{tikzpicture} & \includegraphics[height=2.5in]{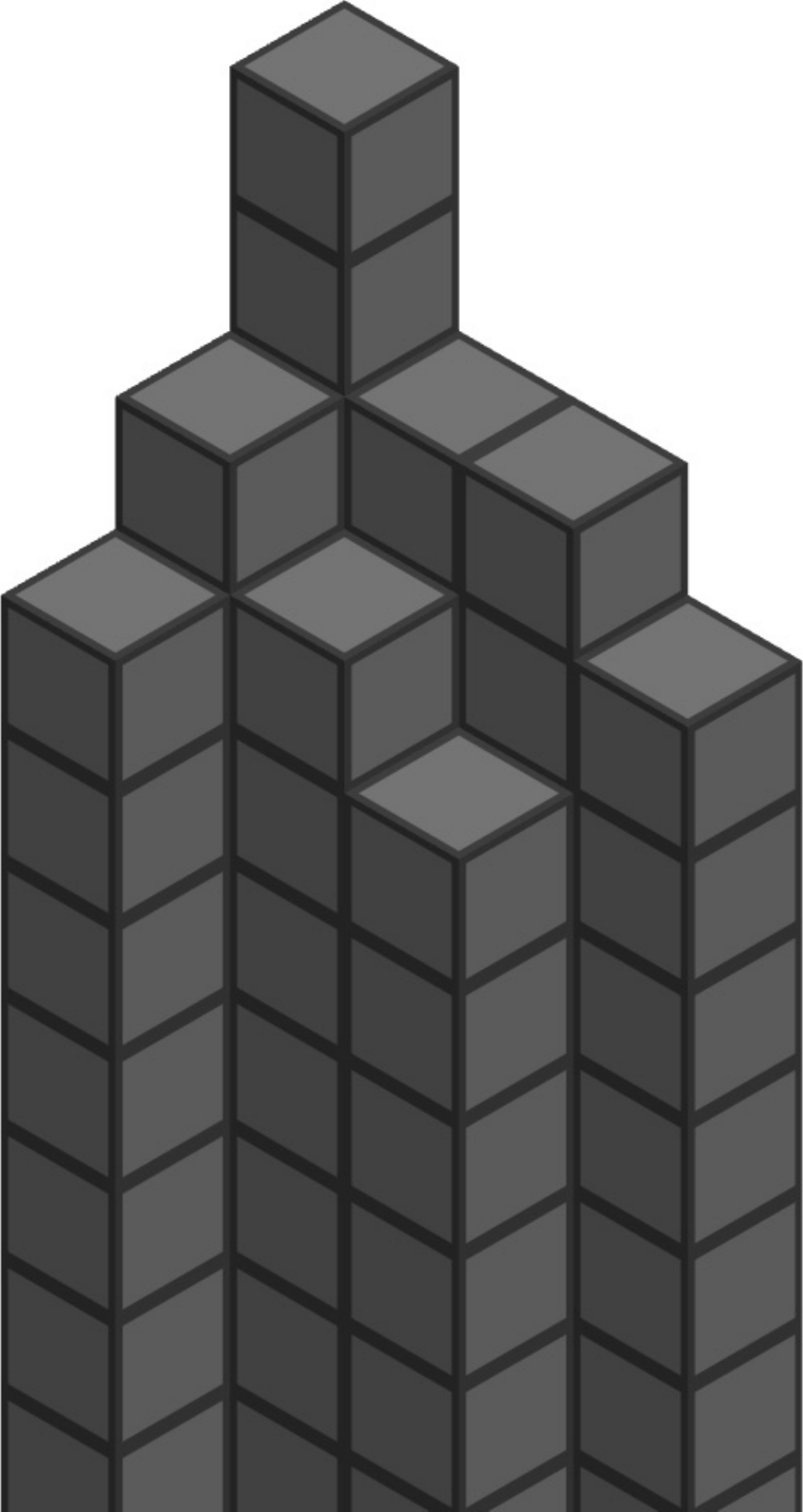}
  \end{tabu} \caption{A one-leg RPP of shape $(\emptyset, \emptyset, (4, 3, 1))$ and weight 19, visualized as $19$ boxes removed from an infinite vertical tower.}
  \label{fig:oneLegDownwardRPP}
\end{figure}

The notion of a two-leg RPP is somewhat more complicated. Observe that the definition of a one-leg RPP (\defref{def:rppAndSpp}) can be directly translated into the language of an SPP without rotating $180^\circ$ by treating the entries as counting the number of blocks \textit{removed} from the top of a tower that descends downward infinitely (rather than added to an empty floor), as in \figref{fig:oneLegDownwardRPP}. This interpretation generalizes more easily --- just as the one leg extends infinitely downward, a two-leg RPP will have legs extending horizontally backward from which blocks are removed, effectively placing a cap on the maximum value of each entry.

\begin{Def}
  Let $\lambda$ and $\mu$ be partitions. An \textbf{reverse plane partition with shape $(\lambda, \mu, \emptyset)$} is a function $\rho : \left( \mathbb{Z} \times \mathbb{N} \cup \mathbb{N} \times \mathbb{Z} \right) \to \mathbb{N}_{\geq 0}$ such that for all $i$ and $j$ for which $\rho(i, j)$ is defined,
  \begin{enumerate}[label=\arabic*.]
    \item $\rho(i,j) \geq \max\{\rho(i + 1, j), \rho(i, j + 1)\}$ (the usual plane partition inequalities).
    \item $\rho(i,j) \leq \min\{\lambda_j, \mu_i\}$, and the inequality is strict for only finitely many pairs $(i, j)$.
  \end{enumerate}
  Similarly to two-leg SPPs, the generating function $W_{(\lambda, \mu, \emptyset)}(q)$ for RPPs with shape $(\lambda, \mu, \emptyset)$ is given by
  \begin{align}
    \label{eq:twoLegRppVertexOperators}
    W_{(\lambda, \mu, \emptyset)}(q) = \leftbra{\mu}\,\prod_{n = -\infty}^{0} \Gamma_+\left(q^{(-n + 1) / 2}\right)\prod_{n = 0}^\infty \Gamma_-\left(q^{(n + 1) / 2}\right)\,\rightket{\lambda},
  \end{align}
  and if $W_0(\lambda, \mu, \emptyset)$ is the minimal exponent of $q$ in this generating function, then the weight of a general RPP $\rho$ with shape $(\lambda, \mu, \emptyset)$ is given by
  $$
    |\rho| = W_0(\lambda, \mu, \emptyset) + \sum \left( \min\{\lambda_j, \mu_i\} - \rho(i, j) \right),
  $$
  where the sum ranges over the entire diagram.
\end{Def}

\begin{figure}
  \begin{tabu}{X[m,c]X[m,c]}
    \begin{tikzpicture}
          
      \begin{scope}[shift={(0,0)}]
        \tableaulabel{0}{0}{}
\tableaulabel{0}{1}{}
\tableaulabel{0}{2}{}
\tableaubox{0}{3}{3}
\tableaubox{0}{4}{1}
\tableaubox{0}{5}{0}
\tableaulabel{1}{0}{}
\tableaulabel{1}{1}{}
\tableaulabel{1}{2}{}
\tableaubox{1}{3}{3}
\tableaubox{1}{4}{1}
\tableaubox{1}{5}{0}
\tableaulabel{2}{0}{}
\tableaulabel{2}{1}{}
\tableaulabel{2}{2}{}
\tableaubox{2}{3}{1}
\tableaubox{2}{4}{1}
\tableaubox{2}{5}{0}
\tableaubox{3}{0}{2}
\tableaubox{3}{1}{2}
\tableaubox{3}{2}{2}
\tableaubox{3}{3}{1}
\tableaubox{3}{4}{0}
\tableaubox{3}{5}{0}
\tableaubox{4}{0}{2}
\tableaubox{4}{1}{2}
\tableaubox{4}{2}{1}
\tableaubox{4}{3}{1}
\tableaubox{4}{4}{0}
\tableaubox{4}{5}{0}
\tableaubox{5}{0}{0}
\tableaubox{5}{1}{0}
\tableaubox{5}{2}{0}
\tableaubox{5}{3}{0}
\tableaubox{5}{4}{0}
\tableaubox{5}{5}{0}
 
\tableaulabel{-1}{4}{$\vdots$}
\tableaulabel{4}{-1}{$\cdots$}
\tableaulabel{4}{6}{$\cdots$}
\tableaulabel{6}{4}{$\vdots$}
      \end{scope}
    
    \end{tikzpicture} & \includegraphics[height=2.5in]{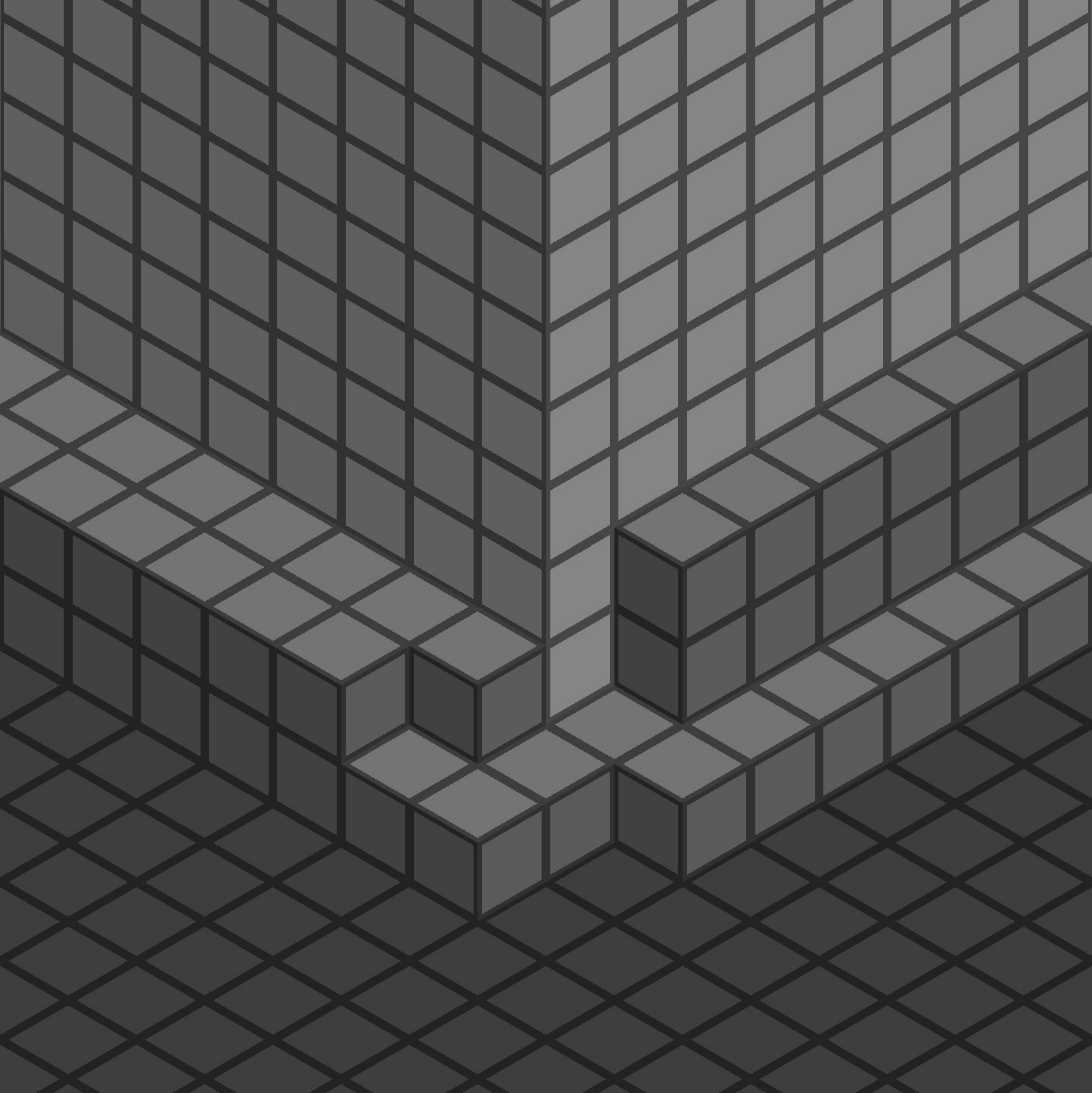}
  \end{tabu} \caption{A two-leg RPP of shape $((3, 1), (2, 2), \emptyset)$ and weight $7$, visualized as a tray with $6$ boxes removed (the minimal configuration has weight $1$ as measured by the vertex operators).}
  \label{fig:twoLegRPP}
\end{figure}

A two-leg RPP has the appearance of a tray pressed up against an infinitely tall building, from which we remove blocks near the corner --- \figref{fig:twoLegRPP} gives an example.

Both two-leg objects are in fact generalizations of their one-leg counterparts: for example, $V_{(\lambda, \emptyset, \emptyset)}(q) = q^k V_{(\emptyset, \emptyset, \lambda)}(q)$ for some $k$ by a $120^\circ$ rotation of the 3-dimensional box diagrams, and similarly for $W_{(\lambda, \emptyset, \emptyset)}(q)$. The factor of $q^k$ arises since the weights of the minimal configurations may be different. Specifically, $V_0(\lambda, \mu, \emptyset)$ is not zero in general; as an example, $V_0((1), \emptyset, \emptyset) = -\frac{1}{2}$.

The vertex operator expressions for two-leg SPPs and RPPs hint at a relationship between the two: $W_{(\lambda, \mu, \emptyset)}(q)$ is obtained from $V_{(\lambda, \mu, \emptyset)}(q)$ by commuting every $\Gamma_-$ past every $\Gamma_+$, then commuting the operators of the same sign so that their order is reversed, and finally swapping the positions of $\lambda$ and $\mu$. In the proofs of \thmref{thm:macMachonsFunctionByToggling} and \thmref{thm:oneLegDecomposition}, the first step alone of these three was sufficient, since regardless of the order of the operators relative to others of the same sign, the $\bra{\emptyset}$ and $\ket{\emptyset}$ bounding the operators ensured that the entire factor counted only the empty partition. As always, commuting a $\Gamma_+$ with a $\Gamma_-$ has the effect of toggling the diagonal containing the relevant edges and therefore changing the shape of the objects counted by the generating function. By \propref{prop:sameSignGammaCommutation} and \lemref{lem:bijectivizingToggles}, commuting two $\Gamma$ operators of the same sign is still a toggling operation, but no longer one that changes the shape. Moreover, without the vertical leg, we also no longer need to handle the bookkeeping of $n$-quotients and various hook lengths.

Before bijectivizing the two-leg version of the PT--DT correspondence, we require an additional technical result. If we limit the size of two-leg SPPs or RPPs to some finite box, effectively limiting $V_{(\lambda, \mu, \emptyset)}(q)$ and $W_{(\lambda, \mu, \emptyset)}(q)$ to finitely many $\Gamma$ operators, then constructing a bijection between the two collections of objects is possible directly from the previous paragraph with little additional work. However, we encounter an issue attempting to define such a bijection for all two-leg SPPs and RPPs: while performing an infinite number of commutations of $\Gamma$ operators does not necessarily present a problem, performing an infinite number of toggles on a plane-partition-like object \textit{does}. To ensure our map is well-defined, we must be able to guarantee that a limiting behavior exists and can be determined with a finite amount of computation; this is the content of the following proposition.

\begin{Prop}
  Let $\lambda, \mu \subset \mathbb{N}^2$ be two Young diagrams and let $\sigma$ be an SPP of shape $(\lambda, \mu, \emptyset)$. Given $n \in \mathbb{N}$, let $\sigma_n$ be the result of toggling diagonals of $\sigma$ until exactly the cells in $[1, n]^2$ have been popped off. Then:

  \begin{enumerate}[label=\arabic*]
    \item There is an $N \in \mathbb{N}$ such that all future toggles of $\sigma_N$ pop off zeros.

    \item Fix $n \geq N$ and let $\alpha(n, i)$ be the partition given by the diagonal of $\sigma_n$ beginning at $(i, n + 1)$. Then for all $i \in \{1, 2, ..., n + 1\}$, $\alpha(n + 1, i) = \alpha(n, i)$. Similarly, if $\beta(n, i)$ is the partition given by the diagonal of $\sigma_n$ beginning at $(n + 1, i)$, then for all $i \in \{1, 2, ..., n + 1\}$, $\beta(n + 1, i) = \beta(n, i)$.

    \item Let $\gamma = \alpha(n, n + 1) = \beta(n, n + 1)$ be the partition on the main diagonal of $\sigma_n$; then $\gamma = \alpha(n, n) = \beta(n, n)$.
  \end{enumerate}
  \label{prop:twoLegLimitBoundary}
\end{Prop}

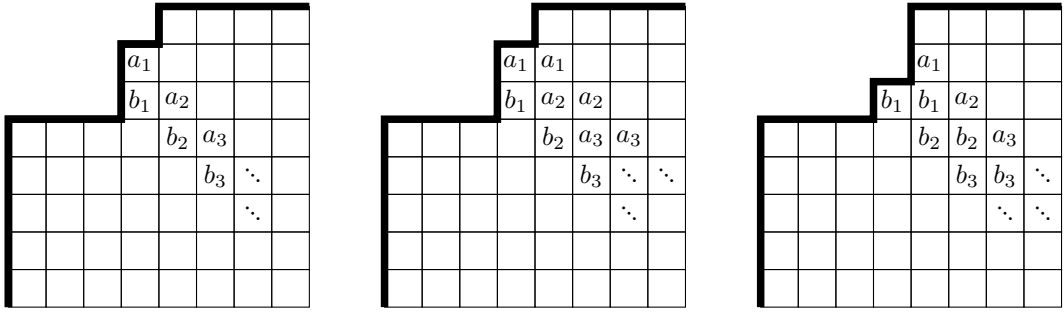
\begin{figure}
  \centering
  \begin{tikzpicture}
    
      \begin{scope}[shift={(0,0)}]
        \tableaulabel{0}{0}{}
\tableaulabel{0}{1}{}
\tableaulabel{0}{2}{}
\tableaulabel{0}{3}{}
\tableaubox{0}{4}{}
\tableaubox{0}{5}{}
\tableaubox{0}{6}{}
\tableaubox{0}{7}{}
\tableaulabel{1}{0}{}
\tableaulabel{1}{1}{}
\tableaulabel{1}{2}{}
\tableaubox{1}{3}{$a_1$}
\tableaubox{1}{4}{}
\tableaubox{1}{5}{}
\tableaubox{1}{6}{}
\tableaubox{1}{7}{}
\tableaulabel{2}{0}{}
\tableaulabel{2}{1}{}
\tableaulabel{2}{2}{}
\tableaubox{2}{3}{$b_1$}
\tableaubox{2}{4}{$a_2$}
\tableaubox{2}{5}{}
\tableaubox{2}{6}{}
\tableaubox{2}{7}{}
\tableaubox{3}{0}{}
\tableaubox{3}{1}{}
\tableaubox{3}{2}{}
\tableaubox{3}{3}{}
\tableaubox{3}{4}{$b_2$}
\tableaubox{3}{5}{$a_3$}
\tableaubox{3}{6}{}
\tableaubox{3}{7}{}
\tableaubox{4}{0}{}
\tableaubox{4}{1}{}
\tableaubox{4}{2}{}
\tableaubox{4}{3}{}
\tableaubox{4}{4}{}
\tableaubox{4}{5}{$b_3$}
\tableaubox{4}{6}{$\ddots$}
\tableaubox{4}{7}{}
\tableaubox{5}{0}{}
\tableaubox{5}{1}{}
\tableaubox{5}{2}{}
\tableaubox{5}{3}{}
\tableaubox{5}{4}{}
\tableaubox{5}{5}{}
\tableaubox{5}{6}{$\ddots$}
\tableaubox{5}{7}{}
\tableaubox{6}{0}{}
\tableaubox{6}{1}{}
\tableaubox{6}{2}{}
\tableaubox{6}{3}{}
\tableaubox{6}{4}{}
\tableaubox{6}{5}{}
\tableaubox{6}{6}{}
\tableaubox{6}{7}{}
\tableaubox{7}{0}{}
\tableaubox{7}{1}{}
\tableaubox{7}{2}{}
\tableaubox{7}{3}{}
\tableaubox{7}{4}{}
\tableaubox{7}{5}{}
\tableaubox{7}{6}{}
\tableaubox{7}{7}{}
\tableaudraw[line width=1mm]{(-0.5, 7.5) -- (-0.5, 3.5) -- (0.5, 3.5) -- (0.5, 2.5) -- (1.5, 2.5) -- (1.5, 2.5) -- (2.5, 2.5) -- (2.5, -0.5) -- (3.5, -0.5) -- (3.5, -0.5) -- (4.5, -0.5) -- (4.5, -0.5) -- (5.5, -0.5) -- (5.5, -0.5) -- (6.5, -0.5) -- (6.5, -0.5) -- (7.5, -0.5)}
\tableaudraw[line width=1mm]{(7.5, -0.5) -- (2.5, -0.5) -- (2.5, 2.5) -- (0.5, 2.5) -- (0.5, 3.5) -- (-0.5, 3.5)}
      \end{scope}
    
      \begin{scope}[shift={(5,0)}]
        \tableaulabel{0}{0}{}
\tableaulabel{0}{1}{}
\tableaulabel{0}{2}{}
\tableaulabel{0}{3}{}
\tableaubox{0}{4}{}
\tableaubox{0}{5}{}
\tableaubox{0}{6}{}
\tableaubox{0}{7}{}
\tableaulabel{1}{0}{}
\tableaulabel{1}{1}{}
\tableaulabel{1}{2}{}
\tableaubox{1}{3}{$a_1$}
\tableaubox{1}{4}{$a_1$}
\tableaubox{1}{5}{}
\tableaubox{1}{6}{}
\tableaubox{1}{7}{}
\tableaulabel{2}{0}{}
\tableaulabel{2}{1}{}
\tableaulabel{2}{2}{}
\tableaubox{2}{3}{$b_1$}
\tableaubox{2}{4}{$a_2$}
\tableaubox{2}{5}{$a_2$}
\tableaubox{2}{6}{}
\tableaubox{2}{7}{}
\tableaubox{3}{0}{}
\tableaubox{3}{1}{}
\tableaubox{3}{2}{}
\tableaubox{3}{3}{}
\tableaubox{3}{4}{$b_2$}
\tableaubox{3}{5}{$a_3$}
\tableaubox{3}{6}{$a_3$}
\tableaubox{3}{7}{}
\tableaubox{4}{0}{}
\tableaubox{4}{1}{}
\tableaubox{4}{2}{}
\tableaubox{4}{3}{}
\tableaubox{4}{4}{}
\tableaubox{4}{5}{$b_3$}
\tableaubox{4}{6}{$\ddots$}
\tableaubox{4}{7}{$\ddots$}
\tableaubox{5}{0}{}
\tableaubox{5}{1}{}
\tableaubox{5}{2}{}
\tableaubox{5}{3}{}
\tableaubox{5}{4}{}
\tableaubox{5}{5}{}
\tableaubox{5}{6}{$\ddots$}
\tableaubox{5}{7}{}
\tableaubox{6}{0}{}
\tableaubox{6}{1}{}
\tableaubox{6}{2}{}
\tableaubox{6}{3}{}
\tableaubox{6}{4}{}
\tableaubox{6}{5}{}
\tableaubox{6}{6}{}
\tableaubox{6}{7}{}
\tableaubox{7}{0}{}
\tableaubox{7}{1}{}
\tableaubox{7}{2}{}
\tableaubox{7}{3}{}
\tableaubox{7}{4}{}
\tableaubox{7}{5}{}
\tableaubox{7}{6}{}
\tableaubox{7}{7}{}
\tableaudraw[line width=1mm]{(-0.5, 7.5) -- (-0.5, 3.5) -- (0.5, 3.5) -- (0.5, 2.5) -- (1.5, 2.5) -- (1.5, 2.5) -- (2.5, 2.5) -- (2.5, -0.5) -- (3.5, -0.5) -- (3.5, -0.5) -- (4.5, -0.5) -- (4.5, -0.5) -- (5.5, -0.5) -- (5.5, -0.5) -- (6.5, -0.5) -- (6.5, -0.5) -- (7.5, -0.5)}
      \end{scope}
    
      \begin{scope}[shift={(10,0)}]
        \tableaulabel{0}{0}{}
\tableaulabel{0}{1}{}
\tableaulabel{0}{2}{}
\tableaulabel{0}{3}{}
\tableaubox{0}{4}{}
\tableaubox{0}{5}{}
\tableaubox{0}{6}{}
\tableaubox{0}{7}{}
\tableaulabel{1}{0}{}
\tableaulabel{1}{1}{}
\tableaulabel{1}{2}{}
\tableaulabel{1}{3}{}
\tableaubox{1}{4}{$a_1$}
\tableaubox{1}{5}{}
\tableaubox{1}{6}{}
\tableaubox{1}{7}{}
\tableaulabel{2}{0}{}
\tableaulabel{2}{1}{}
\tableaulabel{2}{2}{}
\tableaubox{2}{3}{$b_1$}
\tableaubox{2}{4}{$b_1$}
\tableaubox{2}{5}{$a_2$}
\tableaubox{2}{6}{}
\tableaubox{2}{7}{}
\tableaubox{3}{0}{}
\tableaubox{3}{1}{}
\tableaubox{3}{2}{}
\tableaubox{3}{3}{}
\tableaubox{3}{4}{$b_2$}
\tableaubox{3}{5}{$b_2$}
\tableaubox{3}{6}{$a_3$}
\tableaubox{3}{7}{}
\tableaubox{4}{0}{}
\tableaubox{4}{1}{}
\tableaubox{4}{2}{}
\tableaubox{4}{3}{}
\tableaubox{4}{4}{}
\tableaubox{4}{5}{$b_3$}
\tableaubox{4}{6}{$b_3$}
\tableaubox{4}{7}{$\ddots$}
\tableaubox{5}{0}{}
\tableaubox{5}{1}{}
\tableaubox{5}{2}{}
\tableaubox{5}{3}{}
\tableaubox{5}{4}{}
\tableaubox{5}{5}{}
\tableaubox{5}{6}{$\ddots$}
\tableaubox{5}{7}{$\ddots$}
\tableaubox{6}{0}{}
\tableaubox{6}{1}{}
\tableaubox{6}{2}{}
\tableaubox{6}{3}{}
\tableaubox{6}{4}{}
\tableaubox{6}{5}{}
\tableaubox{6}{6}{}
\tableaubox{6}{7}{}
\tableaubox{7}{0}{}
\tableaubox{7}{1}{}
\tableaubox{7}{2}{}
\tableaubox{7}{3}{}
\tableaubox{7}{4}{}
\tableaubox{7}{5}{}
\tableaubox{7}{6}{}
\tableaubox{7}{7}{}
\tableaudraw[line width=1mm]{(-0.5, 7.5) -- (-0.5, 3.5) -- (0.5, 3.5) -- (0.5, 3.5) -- (1.5, 3.5) -- (1.5, 2.5) -- (2.5, 2.5) -- (2.5, -0.5) -- (3.5, -0.5) -- (3.5, -0.5) -- (4.5, -0.5) -- (4.5, -0.5) -- (5.5, -0.5) -- (5.5, -0.5) -- (6.5, -0.5) -- (6.5, -0.5) -- (7.5, -0.5)}
      \end{scope}
    
    \end{tikzpicture} \caption{Three steps of the inductive portion of the proof of \propref{prop:twoLegLimitBoundary}. Left: $a = \alpha(3, 2)$ and $b = \alpha(3, 3)$ in a diagram that is $\sigma_3$ with one toggled corner.}
    \label{fig:toggleInduction}
\end{figure}

\begin{Proof}
  \begin{enumerate}[label=\arabic*]
    \item Toggling every diagonal of $\sigma$ can produce only finitely many nonzero numbers that are popped off; otherwise, $\sigma$ would have infinite weight. Therefore, an $N \in \mathbb{N}$ exists that satisfies the first condition of the proposition. We may also choose $N$ large enough so that $\sigma(i, j) = \max\{\lambda_j, \mu_i\}$ (i.e.\ the minimum possible value) whenever $(i, j) \in \mathbb{N}^2 \setminus [1, N]^2$, and also large enough so that $N \geq \max\{\lambda'_1, \mu'_1\}$. Now if $n \in \mathbb{N}$ and $(i, j) \in \mathbb{N}^2 \setminus [1, n]^2$ satisfies $\sigma_n(i, j) \neq \sigma(i, j)$, then $|j - i|$ must be at least $n$ --- that is, $(i, j)$ must lie outside the diagonals that were toggled to produce $\sigma_n$. Therefore, for all $n \geq N$ and all $(i, j)$ with $|j - i| \geq n$, $\sigma_n(i, j) = \max\{\lambda_j, \mu_i\}$. In particular, $\sigma_n(i, j) = \mu_i$ whenever $n \geq N$ and $j - i \geq N$.

    \item For fixed $n \geq N$, we prove the second claim of the proposition by induction on $i$. The base case is shown by the previous paragraph: $\alpha(n, 1) = \alpha(n + 1, 1) = \mu$. For the induction step, suppose $i \geq 2$ and $\alpha(n, i - 1) = \alpha(n + 1, i - 1)$. For ease of notation, set $a = \alpha(n, i - 1)$ and $b = \alpha(n, i)$; if we begin with $\sigma_n$ and toggle the diagonals beginning with cells $(k, n + 1)$ for all $k \in \{1, ..., i - 2\}$, the result is then the leftmost diagram of \figref{fig:toggleInduction}. Since the diagonal immediately above $a$ will not change from any of the remaining toggles that produce $\sigma_{n + 1}$, it is equal to $\alpha(n + 1, i - 1)$; therefore, the induction hypothesis guarantees that it is in fact equal to $a$, as in the middle diagram of \figref{fig:toggleInduction}. When we toggle the next diagonal down, the plane partition inequalities guarantee that $a_k \geq b_k$ for all $k \in \mathbb{N}$. Therefore, each $a_k$ toggles to
    $$
      \min\{a_{k - 1}, b_{k - 1}\} + \max\{a_k, b_k\} - a_k = b_{k - 1} + a_k - a_k = b_{k - 1},
    $$
    resulting in the rightmost diagram of \figref{fig:toggleInduction} and proving the claim. An exactly symmetric argument shows the result for the $\beta$ partitions.

    \item The weight of $\sigma_n$ is given by the generating function
    $$
      \leftbra{\lambda}\,\prod_{k = -\infty}^{-n} \Gamma_-\left(q^{(-k + 1) / 2}\right) \prod_{k = 0}^{n - 1} \Gamma_+\left(q^{(k + 1) / 2}\right) \prod_{k = -n + 1}^0 \Gamma_-\left(q^{(-k + 1) / 2}\right) \prod_{k = n}^\infty \Gamma_+\left(q^{(k + 1) / 2}\right)\,\rightket{\mu},
    $$
    so the contribution to the weight from the main diagonal $\gamma$ is given by the middle two $\Gamma$ operators; that is, 
    $$
      \frac{n}{2} \left( |\alpha(n, n)| - |\gamma| \right) + \frac{n}{2} \left( |\beta(n, n)| - |\gamma| \right).
    $$
    By the previous part of the proposition, $\alpha(n + 1, n + 1) = \alpha(n, n + 1) = \gamma$ and $\beta(n + 1, n + 1) = \alpha(n, n + 1) = \gamma$, so if $\gamma' = \alpha(n + 1, n + 2) = \beta(n + 1, n + 2)$ is the partition on the main diagonal of $\sigma_{n + 1}$, then its contribution to the weight of $\sigma_{n + 1}$ is
    \begin{align*}
      \frac{n + 1}{2} \left( |\alpha(n + 1, n + 1)| - |\gamma'| \right) + \frac{n + 1}{2} \left( |\beta(n + 1, n + 1)| - |\gamma'| \right) &= \frac{n + 1}{2} \left( |\gamma| - |\gamma'| \right) + \frac{n + 1}{2} \left( |\gamma| - |\gamma'| \right)\\
      &= (n + 1)\left( |\gamma| - |\gamma'| \right).
    \end{align*}
    However, the weights of $\sigma_n$ and $\sigma_{n + 1}$ are equal since no nonzero numbers are popped off in the toggles producing $\sigma_{n + 1}$. Comparing the two generating functions, the only difference in weight is exactly the contribution from $\gamma'$, so $|\gamma| = |\gamma'|$. Since $\gamma \succ \gamma'$, the only way for this to occur is if $\gamma = \gamma'$.
  \end{enumerate}
\end{Proof}

\begin{Thm}
  Let $\lambda, \mu \subset \mathbb{N}^2$ be two Young diagrams. Then there is a weight-preserving bijection between SPPs $\sigma$ of shape $(\lambda, \mu, \emptyset)$ and pairs $(\rho, \pi)$ of RPPs $\rho$ of shape $(\lambda, \mu, \emptyset)$ and plane partitions $\pi$, where $|\sigma| = |\rho| + |\pi|$.
  \label{thm:twoLegDecomposition}
\end{Thm}

\begin{Proof}
  As in the proofs of \thmref{thm:macMachonsFunctionByToggling} and \thmref{thm:oneLegDecomposition}, our bijection follows the algebraic proof. Two-leg RPPs and SPPs do not occur frequently in the literature, but the existence of an algebraic proof that $V_{(\lambda, \mu, \emptyset)}(q) = W_{(\lambda, \mu, \emptyset)}(q)M(q)$ was mentioned in \cite{PT}; we provide such a proof in full detail here. We begin by commuting the $\Gamma$ operators of opposite sign in the expression
  $$
    V_{(\lambda, \mu, \emptyset)}(q) = \leftbra{\mu}\,\prod_{n = -\infty}^0 \Gamma_-\left(q^{(-n + 1) / 2}\right)\prod_{n = 0}^\infty \Gamma_+\left(q^{(n + 1) / 2}\right)\,\rightket{\lambda},
  $$
  producing
  \begin{align}
    V_{(\lambda, \mu, \emptyset)}(q) = M(q)\leftbra{\mu}\,\prod_{n = 0}^\infty \Gamma_+\left(q^{(n + 1) / 2}\right)\prod_{n = -\infty}^0 \Gamma_-\left(q^{(-n + 1) / 2}\right)\,\rightket{\lambda}.
    \label{eq:twoLegVertexOperatorsBeforePalindrome}
  \end{align}
  Since the only difference between this vertex operator product and \eqref{eq:macMahonWithArguments} is that it contains $\bra{\mu}$ and $\ket{\lambda}$ instead of $\bra{\emptyset}$ and $\ket{\emptyset}$, these commutations produce a factor of $M(q)$.

  We now ``palindromically'' commute the operators of the same sign (i.e\ reverse their order). This produces no factors --- it only reweights the objects that are counted to be closer to our usual definitions. The result is
  $$
    V_{(\lambda, \mu, \emptyset)}(q) = M(q) W_{(\mu, \lambda, \emptyset)}(q).
  $$
  However, $W_{(\mu, \lambda, \emptyset)}(q) = W_{(\lambda, \mu, \emptyset)}(q)$; the transpose is a weight-preserving bijection between the sets counted by these two generating functions.

  Exactly as in the proof of \thmref{thm:macMachonsFunctionByToggling}, we define a map $\tau$ on SPPs of shape $(\lambda, \mu, \emptyset)$ by iteratively toggling diagonals. In that proof, we could justify that the map was well-defined since every plane partition contains only finitely many nonzero entries, but here we must be more careful. Let $\sigma$ be an SPP of shape $(\lambda, \mu, \emptyset)$; with the $N \in \mathbb{N}$ guaranteed by \propref{prop:twoLegLimitBoundary}, we may perform toggles until only every cell in $[1, N]^2$ is popped off, then use just those values to create an $N \times N$ tableau that we can untoggle into a plane partition $\pi$ using \thmref{thm:macMachonsFunctionByToggling}.
  
  The remaining toggled object ($\sigma_N$ in the parlance of \propref{prop:twoLegLimitBoundary}) is constant on diagonals sufficiently close the main diagonal, and by that proposition, it continues to have that property as we toggle more and more diagonals. Therefore, we do not need to perform any further toggles to determine the limiting object, an RPP-like remnant infinitely far from the origin. We then perform the toggles corresponding to the same-sign commutations. We first commute $\Gamma_-\left( q^{1/2} \right)$ to the left past every other $\Gamma_-$ operator, then $\Gamma_-\left( q^{3/2} \right)$ to the left past every other $\Gamma_-$ except for $\Gamma_-\left( q^{1/2} \right)$, and so on. Each of these commutations nominally involves an infinite number of toggles, but since the RPP-like object is eventually constant on diagonals, we may determine the output after only a finite number of toggles. Similarly, we need only commute a finite number of the $\Gamma_-\left( q^{k/2} \right)$ to the left in total, since the final RPP must also eventually be constant on diagonals far enough away from the origin in order to have finite weight. The result is a two-leg RPP of shape $(\mu, \lambda, \emptyset)$; exactly as in the algebraic proof, we then transpose the diagram to produce an RPP $\rho$ of shape $(\lambda, \mu, \emptyset)$, as required.
\end{Proof}

\begin{figure}
  \centering
  \begin{tikzpicture}
    
      \begin{scope}[shift={(0,0)}]
        \tableaubox{0}{0}{6}
\tableaubox{0}{1}{5}
\tableaubox{0}{2}{3}
\tableaubox{0}{3}{3}
\tableaubox{0}{4}{3}
\tableaubox{0}{5}{3}
\tableaubox{1}{0}{5}
\tableaubox{1}{1}{3}
\tableaubox{1}{2}{3}
\tableaubox{1}{3}{1}
\tableaubox{1}{4}{1}
\tableaubox{1}{5}{1}
\tableaubox{2}{0}{3}
\tableaubox{2}{1}{3}
\tableaubox{2}{2}{2}
\tableaubox{2}{3}{0}
\tableaubox{2}{4}{0}
\tableaubox{2}{5}{0}
\tableaubox{3}{0}{2}
\tableaubox{3}{1}{2}
\tableaubox{3}{2}{0}
\tableaubox{3}{3}{0}
\tableaubox{3}{4}{0}
\tableaubox{3}{5}{0}
\tableaubox{4}{0}{2}
\tableaubox{4}{1}{2}
\tableaubox{4}{2}{0}
\tableaubox{4}{3}{0}
\tableaubox{4}{4}{0}
\tableaubox{4}{5}{0}
\tableaubox{5}{0}{2}
\tableaubox{5}{1}{2}
\tableaubox{5}{2}{0}
\tableaubox{5}{3}{0}
\tableaubox{5}{4}{0}
\tableaubox{5}{5}{0}
\tableaulabel{2.5}{6}{$\cdots$}
\tableaulabel{6}{2.5}{$\vdots$}

      \end{scope}
    
    \end{tikzpicture} \caption{An SPP of shape $((2, 2), (3, 1), \emptyset)$ and weight $16$.}
    \label{fig:twoLegSPPExample}
\end{figure}

\begin{figure}
  \centering
  \begin{tikzpicture}
    
      \begin{scope}[shift={(0,0)}]
        \tableaubox{0}{0}{6}
\tableaubox{0}{1}{5}
\tableaubox{0}{2}{3}
\tableaubox{0}{3}{3}
\tableaubox{0}{4}{3}
\tableaubox{0}{5}{3}
\tableaubox{1}{0}{5}
\tableaubox{1}{1}{3}
\tableaubox{1}{2}{3}
\tableaubox{1}{3}{1}
\tableaubox{1}{4}{1}
\tableaubox{1}{5}{1}
\tableaubox{2}{0}{3}
\tableaubox{2}{1}{3}
\tableaubox{2}{2}{2}
\tableaubox{2}{3}{0}
\tableaubox{2}{4}{0}
\tableaubox{2}{5}{0}
\tableaubox{3}{0}{2}
\tableaubox{3}{1}{2}
\tableaubox{3}{2}{0}
\tableaubox{3}{3}{0}
\tableaubox{3}{4}{0}
\tableaubox{3}{5}{0}
\tableaubox{4}{0}{2}
\tableaubox{4}{1}{2}
\tableaubox{4}{2}{0}
\tableaubox{4}{3}{0}
\tableaubox{4}{4}{0}
\tableaubox{4}{5}{0}
\tableaubox{5}{0}{2}
\tableaubox{5}{1}{2}
\tableaubox{5}{2}{0}
\tableaubox{5}{3}{0}
\tableaubox{5}{4}{0}
\tableaubox{5}{5}{0}

      \tableaulabel{2.5}{6}{$\mapsto$}
    \tableaudraw[line width=1mm]{(-0.5, 5.5) -- (-0.5, -0.5) -- (0.5, -0.5) -- (0.5, -0.5) -- (1.5, -0.5) -- (1.5, -0.5) -- (2.5, -0.5) -- (2.5, -0.5) -- (3.5, -0.5) -- (3.5, -0.5) -- (4.5, -0.5) -- (4.5, -0.5) -- (5.5, -0.5)}
      \end{scope}
    
      \begin{scope}[shift={(3.5,0)}]
        \tableaubox{0}{0}{1}
\tableaubox{0}{1}{5}
\tableaubox{0}{2}{3}
\tableaubox{0}{3}{3}
\tableaubox{0}{4}{3}
\tableaubox{0}{5}{3}
\tableaubox{1}{0}{5}
\tableaubox{1}{1}{5}
\tableaubox{1}{2}{3}
\tableaubox{1}{3}{1}
\tableaubox{1}{4}{1}
\tableaubox{1}{5}{1}
\tableaubox{2}{0}{3}
\tableaubox{2}{1}{3}
\tableaubox{2}{2}{1}
\tableaubox{2}{3}{0}
\tableaubox{2}{4}{0}
\tableaubox{2}{5}{0}
\tableaubox{3}{0}{2}
\tableaubox{3}{1}{2}
\tableaubox{3}{2}{0}
\tableaubox{3}{3}{0}
\tableaubox{3}{4}{0}
\tableaubox{3}{5}{0}
\tableaubox{4}{0}{2}
\tableaubox{4}{1}{2}
\tableaubox{4}{2}{0}
\tableaubox{4}{3}{0}
\tableaubox{4}{4}{0}
\tableaubox{4}{5}{0}
\tableaubox{5}{0}{2}
\tableaubox{5}{1}{2}
\tableaubox{5}{2}{0}
\tableaubox{5}{3}{0}
\tableaubox{5}{4}{0}
\tableaubox{5}{5}{0}

      \tableaulabel{2.5}{6}{$\mapsto$}
    \tableaudraw[line width=1mm]{(-0.5, 5.5) -- (-0.5, 0.5) -- (0.5, 0.5) -- (0.5, -0.5) -- (1.5, -0.5) -- (1.5, -0.5) -- (2.5, -0.5) -- (2.5, -0.5) -- (3.5, -0.5) -- (3.5, -0.5) -- (4.5, -0.5) -- (4.5, -0.5) -- (5.5, -0.5)}
      \end{scope}
    
      \begin{scope}[shift={(7,0)}]
        \tableaubox{0}{0}{1}
\tableaubox{0}{1}{0}
\tableaubox{0}{2}{3}
\tableaubox{0}{3}{3}
\tableaubox{0}{4}{3}
\tableaubox{0}{5}{3}
\tableaubox{1}{0}{0}
\tableaubox{1}{1}{5}
\tableaubox{1}{2}{1}
\tableaubox{1}{3}{1}
\tableaubox{1}{4}{1}
\tableaubox{1}{5}{1}
\tableaubox{2}{0}{3}
\tableaubox{2}{1}{2}
\tableaubox{2}{2}{1}
\tableaubox{2}{3}{1}
\tableaubox{2}{4}{0}
\tableaubox{2}{5}{0}
\tableaubox{3}{0}{2}
\tableaubox{3}{1}{2}
\tableaubox{3}{2}{1}
\tableaubox{3}{3}{0}
\tableaubox{3}{4}{0}
\tableaubox{3}{5}{0}
\tableaubox{4}{0}{2}
\tableaubox{4}{1}{2}
\tableaubox{4}{2}{0}
\tableaubox{4}{3}{0}
\tableaubox{4}{4}{0}
\tableaubox{4}{5}{0}
\tableaubox{5}{0}{2}
\tableaubox{5}{1}{2}
\tableaubox{5}{2}{0}
\tableaubox{5}{3}{0}
\tableaubox{5}{4}{0}
\tableaubox{5}{5}{0}

      \tableaulabel{2.5}{6}{$\mapsto$}
    \tableaudraw[line width=1mm]{(-0.5, 5.5) -- (-0.5, 1.5) -- (0.5, 1.5) -- (0.5, 0.5) -- (1.5, 0.5) -- (1.5, -0.5) -- (2.5, -0.5) -- (2.5, -0.5) -- (3.5, -0.5) -- (3.5, -0.5) -- (4.5, -0.5) -- (4.5, -0.5) -- (5.5, -0.5)}
      \end{scope}
    
      \begin{scope}[shift={(10.5,0)}]
        \tableaubox{0}{0}{1}
\tableaubox{0}{1}{0}
\tableaubox{0}{2}{3}
\tableaubox{0}{3}{3}
\tableaubox{0}{4}{3}
\tableaubox{0}{5}{3}
\tableaubox{1}{0}{0}
\tableaubox{1}{1}{3}
\tableaubox{1}{2}{1}
\tableaubox{1}{3}{1}
\tableaubox{1}{4}{1}
\tableaubox{1}{5}{1}
\tableaubox{2}{0}{3}
\tableaubox{2}{1}{2}
\tableaubox{2}{2}{1}
\tableaubox{2}{3}{1}
\tableaubox{2}{4}{0}
\tableaubox{2}{5}{0}
\tableaubox{3}{0}{2}
\tableaubox{3}{1}{2}
\tableaubox{3}{2}{1}
\tableaubox{3}{3}{1}
\tableaubox{3}{4}{0}
\tableaubox{3}{5}{0}
\tableaubox{4}{0}{2}
\tableaubox{4}{1}{2}
\tableaubox{4}{2}{0}
\tableaubox{4}{3}{0}
\tableaubox{4}{4}{0}
\tableaubox{4}{5}{0}
\tableaubox{5}{0}{2}
\tableaubox{5}{1}{2}
\tableaubox{5}{2}{0}
\tableaubox{5}{3}{0}
\tableaubox{5}{4}{0}
\tableaubox{5}{5}{0}

      \tableaulabel{2.5}{6}{$\mapsto$}
    \tableaudraw[line width=1mm]{(-0.5, 5.5) -- (-0.5, 1.5) -- (0.5, 1.5) -- (0.5, 1.5) -- (1.5, 1.5) -- (1.5, -0.5) -- (2.5, -0.5) -- (2.5, -0.5) -- (3.5, -0.5) -- (3.5, -0.5) -- (4.5, -0.5) -- (4.5, -0.5) -- (5.5, -0.5)}
      \end{scope}
    
  \end{tikzpicture}

  \vspace{18pt}

  \begin{tikzpicture}
    
      \begin{scope}[shift={(0,0)}]
        \tableaubox{0}{0}{1}
\tableaubox{0}{1}{0}
\tableaubox{0}{2}{0}
\tableaubox{0}{3}{3}
\tableaubox{0}{4}{3}
\tableaubox{0}{5}{3}
\tableaubox{1}{0}{0}
\tableaubox{1}{1}{3}
\tableaubox{1}{2}{1}
\tableaubox{1}{3}{1}
\tableaubox{1}{4}{1}
\tableaubox{1}{5}{1}
\tableaubox{2}{0}{1}
\tableaubox{2}{1}{2}
\tableaubox{2}{2}{1}
\tableaubox{2}{3}{1}
\tableaubox{2}{4}{1}
\tableaubox{2}{5}{0}
\tableaubox{3}{0}{2}
\tableaubox{3}{1}{2}
\tableaubox{3}{2}{1}
\tableaubox{3}{3}{1}
\tableaubox{3}{4}{0}
\tableaubox{3}{5}{0}
\tableaubox{4}{0}{2}
\tableaubox{4}{1}{2}
\tableaubox{4}{2}{1}
\tableaubox{4}{3}{0}
\tableaubox{4}{4}{0}
\tableaubox{4}{5}{0}
\tableaubox{5}{0}{2}
\tableaubox{5}{1}{2}
\tableaubox{5}{2}{0}
\tableaubox{5}{3}{0}
\tableaubox{5}{4}{0}
\tableaubox{5}{5}{0}

      \tableaulabel{2.5}{6}{$\mapsto$}
    \tableaudraw[line width=1mm]{(-0.5, 5.5) -- (-0.5, 2.5) -- (0.5, 2.5) -- (0.5, 1.5) -- (1.5, 1.5) -- (1.5, 0.5) -- (2.5, 0.5) -- (2.5, -0.5) -- (3.5, -0.5) -- (3.5, -0.5) -- (4.5, -0.5) -- (4.5, -0.5) -- (5.5, -0.5)}
      \end{scope}
    
      \begin{scope}[shift={(3.5,0)}]
        \tableaubox{0}{0}{1}
\tableaubox{0}{1}{0}
\tableaubox{0}{2}{0}
\tableaubox{0}{3}{3}
\tableaubox{0}{4}{3}
\tableaubox{0}{5}{3}
\tableaubox{1}{0}{0}
\tableaubox{1}{1}{3}
\tableaubox{1}{2}{0}
\tableaubox{1}{3}{1}
\tableaubox{1}{4}{1}
\tableaubox{1}{5}{1}
\tableaubox{2}{0}{1}
\tableaubox{2}{1}{0}
\tableaubox{2}{2}{1}
\tableaubox{2}{3}{1}
\tableaubox{2}{4}{1}
\tableaubox{2}{5}{0}
\tableaubox{3}{0}{2}
\tableaubox{3}{1}{2}
\tableaubox{3}{2}{1}
\tableaubox{3}{3}{1}
\tableaubox{3}{4}{1}
\tableaubox{3}{5}{0}
\tableaubox{4}{0}{2}
\tableaubox{4}{1}{2}
\tableaubox{4}{2}{1}
\tableaubox{4}{3}{1}
\tableaubox{4}{4}{0}
\tableaubox{4}{5}{0}
\tableaubox{5}{0}{2}
\tableaubox{5}{1}{2}
\tableaubox{5}{2}{0}
\tableaubox{5}{3}{0}
\tableaubox{5}{4}{0}
\tableaubox{5}{5}{0}

      \tableaulabel{2.5}{6}{$\mapsto$}
    \tableaudraw[line width=1mm]{(-0.5, 5.5) -- (-0.5, 2.5) -- (0.5, 2.5) -- (0.5, 2.5) -- (1.5, 2.5) -- (1.5, 1.5) -- (2.5, 1.5) -- (2.5, -0.5) -- (3.5, -0.5) -- (3.5, -0.5) -- (4.5, -0.5) -- (4.5, -0.5) -- (5.5, -0.5)}
      \end{scope}
    
      \begin{scope}[shift={(7,0)}]
        \tableaubox{0}{0}{1}
\tableaubox{0}{1}{0}
\tableaubox{0}{2}{0}
\tableaubox{0}{3}{3}
\tableaubox{0}{4}{3}
\tableaubox{0}{5}{3}
\tableaubox{1}{0}{0}
\tableaubox{1}{1}{3}
\tableaubox{1}{2}{0}
\tableaubox{1}{3}{1}
\tableaubox{1}{4}{1}
\tableaubox{1}{5}{1}
\tableaubox{2}{0}{1}
\tableaubox{2}{1}{0}
\tableaubox{2}{2}{0}
\tableaubox{2}{3}{1}
\tableaubox{2}{4}{1}
\tableaubox{2}{5}{0}
\tableaubox{3}{0}{2}
\tableaubox{3}{1}{2}
\tableaubox{3}{2}{1}
\tableaubox{3}{3}{1}
\tableaubox{3}{4}{1}
\tableaubox{3}{5}{0}
\tableaubox{4}{0}{2}
\tableaubox{4}{1}{2}
\tableaubox{4}{2}{1}
\tableaubox{4}{3}{1}
\tableaubox{4}{4}{1}
\tableaubox{4}{5}{0}
\tableaubox{5}{0}{2}
\tableaubox{5}{1}{2}
\tableaubox{5}{2}{0}
\tableaubox{5}{3}{0}
\tableaubox{5}{4}{0}
\tableaubox{5}{5}{0}

      \tableaulabel{2.5}{6}{$\mapsto$}
      \tableaulabel{2.5}{7}{$\cdots$}
      \tableaulabel{2.5}{8}{$\mapsto$}
    \tableaudraw[line width=1mm]{(-0.5, 5.5) -- (-0.5, 2.5) -- (0.5, 2.5) -- (0.5, 2.5) -- (1.5, 2.5) -- (1.5, 2.5) -- (2.5, 2.5) -- (2.5, -0.5) -- (3.5, -0.5) -- (3.5, -0.5) -- (4.5, -0.5) -- (4.5, -0.5) -- (5.5, -0.5)}
      \end{scope}
    
      \begin{scope}[shift={(11.5,0)}]
        \tableaubox{0}{0}{1}
\tableaubox{0}{1}{0}
\tableaubox{0}{2}{0}
\tableaubox{0}{3}{0}
\tableaubox{0}{4}{3}
\tableaubox{0}{5}{3}
\tableaubox{1}{0}{0}
\tableaubox{1}{1}{3}
\tableaubox{1}{2}{0}
\tableaubox{1}{3}{0}
\tableaubox{1}{4}{1}
\tableaubox{1}{5}{1}
\tableaubox{2}{0}{1}
\tableaubox{2}{1}{0}
\tableaubox{2}{2}{0}
\tableaubox{2}{3}{0}
\tableaubox{2}{4}{1}
\tableaubox{2}{5}{1}
\tableaubox{3}{0}{0}
\tableaubox{3}{1}{0}
\tableaubox{3}{2}{0}
\tableaubox{3}{3}{0}
\tableaubox{3}{4}{1}
\tableaubox{3}{5}{1}
\tableaubox{4}{0}{2}
\tableaubox{4}{1}{2}
\tableaubox{4}{2}{1}
\tableaubox{4}{3}{1}
\tableaubox{4}{4}{1}
\tableaubox{4}{5}{1}
\tableaubox{5}{0}{2}
\tableaubox{5}{1}{2}
\tableaubox{5}{2}{1}
\tableaubox{5}{3}{1}
\tableaubox{5}{4}{1}
\tableaubox{5}{5}{1}
\tableaudraw[line width=1mm]{(-0.5, 5.5) -- (-0.5, 3.5) -- (0.5, 3.5) -- (0.5, 3.5) -- (1.5, 3.5) -- (1.5, 3.5) -- (2.5, 3.5) -- (2.5, 3.5) -- (3.5, 3.5) -- (3.5, -0.5) -- (4.5, -0.5) -- (4.5, -0.5) -- (5.5, -0.5)}
      \end{scope}
    
  \end{tikzpicture} \caption{Iteratively toggling the diagonals of a two-leg SPP.}
  \label{fig:twoLegSPPToggle}
\end{figure}
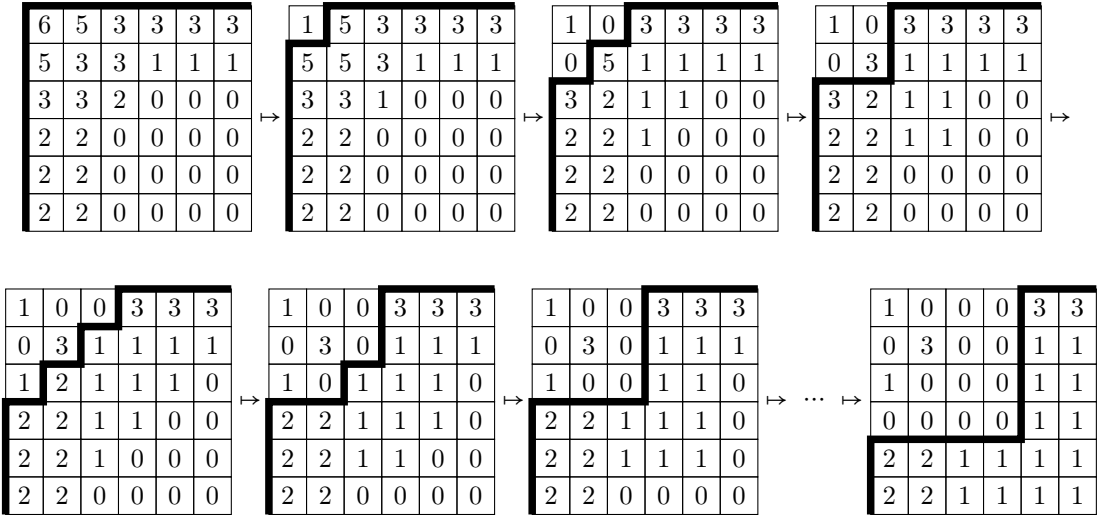

\begin{figure}
  \centering
  \begin{tikzpicture}
    
      \begin{scope}[shift={(0,0)}]
        \tableaubox{0}{0}{1}
\tableaubox{0}{1}{0}
\tableaubox{0}{2}{0}
\tableaubox{0}{3}{0}
\tableaulabel{0}{4}{}
\tableaubox{0}{5}{3}
\tableaulabel{0}{6}{}
\tableaubox{1}{0}{0}
\tableaubox{1}{1}{3}
\tableaubox{1}{2}{0}
\tableaubox{1}{3}{0}
\tableaulabel{1}{4}{}
\tableaubox{1}{5}{1}
\tableaubox{1}{6}{1}
\tableaubox{2}{0}{1}
\tableaubox{2}{1}{0}
\tableaubox{2}{2}{0}
\tableaubox{2}{3}{0}
\tableaulabel{2}{4}{}
\tableaubox{2}{5}{1}
\tableaubox{2}{6}{1}
\tableaubox{3}{0}{0}
\tableaubox{3}{1}{0}
\tableaubox{3}{2}{0}
\tableaubox{3}{3}{0}
\tableaulabel{3}{4}{}
\tableaubox{3}{5}{1}
\tableaubox{3}{6}{1}
\tableaulabel{4}{0}{}
\tableaulabel{4}{1}{}
\tableaulabel{4}{2}{}
\tableaulabel{4}{3}{}
\tableaulabel{4}{4}{}
\tableaulabel{4}{5}{}
\tableaulabel{4}{6}{}
\tableaubox{5}{0}{2}
\tableaubox{5}{1}{2}
\tableaubox{5}{2}{1}
\tableaubox{5}{3}{1}
\tableaulabel{5}{4}{}
\tableaubox{5}{5}{1}
\tableaubox{5}{6}{1}
\tableaulabel{6}{0}{}
\tableaubox{6}{1}{2}
\tableaubox{6}{2}{1}
\tableaubox{6}{3}{1}
\tableaulabel{6}{4}{}
\tableaubox{6}{5}{1}
\tableaubox{6}{6}{1}

\tableaulabel{1.5}{4}{$\cdots$}
\tableaulabel{4}{1.5}{$\vdots$}

\tableaulabel{5.5}{4}{$\cdots$}
\tableaulabel{4}{5.5}{$\vdots$}

\tableaudraw[line width=1mm]{(4.5, -0.5) -- (4.5, 4.5) -- (-0.5, 4.5)}
      \end{scope}
    
  \end{tikzpicture} \caption{The limiting diagram after toggling the SPP in \figref{fig:twoLegSPPToggle}.}
  \label{fig:twoLegSPPToggleResult}
\end{figure}

\begin{figure}
  \centering
  \begin{tikzpicture}
    
      \begin{scope}[shift={(0,0)}]
        \tableaulabel{0}{0}{}
\tableaulabel{0}{1}{}
\tableaulabel{0}{2}{}
\tableaulabel{0}{3}{}
\tableaulabel{0}{4}{}
\tableaubox{0}{5}{3}
\tableaulabel{0}{6}{}
\tableaulabel{1}{0}{}
\tableaulabel{1}{1}{}
\tableaulabel{1}{2}{}
\tableaulabel{1}{3}{}
\tableaulabel{1}{4}{}
\tableaubox{1}{5}{1}
\tableaubox{1}{6}{1}
\tableaulabel{2}{0}{}
\tableaulabel{2}{1}{}
\tableaulabel{2}{2}{}
\tableaulabel{2}{3}{}
\tableaulabel{2}{4}{}
\tableaubox{2}{5}{1}
\tableaubox{2}{6}{1}
\tableaulabel{3}{0}{}
\tableaulabel{3}{1}{}
\tableaulabel{3}{2}{}
\tableaulabel{3}{3}{}
\tableaulabel{3}{4}{}
\tableaulabel{3}{5}{}
\tableaulabel{3}{6}{}
\tableaulabel{4}{0}{}
\tableaulabel{4}{1}{}
\tableaulabel{4}{2}{}
\tableaulabel{4}{3}{}
\tableaulabel{4}{4}{}
\tableaubox{4}{5}{1}
\tableaubox{4}{6}{1}
\tableaubox{5}{0}{2}
\tableaubox{5}{1}{2}
\tableaubox{5}{2}{1}
\tableaulabel{5}{3}{}
\tableaubox{5}{4}{1}
\tableaubox{5}{5}{1}
\tableaubox{5}{6}{1}
\tableaulabel{6}{0}{}
\tableaubox{6}{1}{2}
\tableaubox{6}{2}{1}
\tableaulabel{6}{3}{}
\tableaubox{6}{4}{1}
\tableaubox{6}{5}{1}
\tableaubox{6}{6}{1}

      \tableaulabel{3}{7}{$\mapsto$}
    
\tableaulabel{5.5}{3}{$\cdots$}
\tableaulabel{3}{5.5}{$\vdots$}

\tableaudraw[line width=1mm]{(4.5, -0.5) -- (4.5, 4.5) -- (-0.5, 4.5)}
      \end{scope}
    
      \begin{scope}[shift={(4,0)}]
        \tableaulabel{0}{0}{}
\tableaulabel{0}{1}{}
\tableaulabel{0}{2}{}
\tableaulabel{0}{3}{}
\tableaulabel{0}{4}{}
\tableaubox{0}{5}{3}
\tableaulabel{0}{6}{}
\tableaulabel{1}{0}{}
\tableaulabel{1}{1}{}
\tableaulabel{1}{2}{}
\tableaulabel{1}{3}{}
\tableaulabel{1}{4}{}
\tableaubox{1}{5}{3}
\tableaubox{1}{6}{1}
\tableaulabel{2}{0}{}
\tableaulabel{2}{1}{}
\tableaulabel{2}{2}{}
\tableaulabel{2}{3}{}
\tableaulabel{2}{4}{}
\tableaubox{2}{5}{1}
\tableaubox{2}{6}{1}
\tableaulabel{3}{0}{}
\tableaulabel{3}{1}{}
\tableaulabel{3}{2}{}
\tableaulabel{3}{3}{}
\tableaulabel{3}{4}{}
\tableaulabel{3}{5}{}
\tableaulabel{3}{6}{}
\tableaulabel{4}{0}{}
\tableaulabel{4}{1}{}
\tableaulabel{4}{2}{}
\tableaulabel{4}{3}{}
\tableaulabel{4}{4}{}
\tableaubox{4}{5}{1}
\tableaubox{4}{6}{1}
\tableaubox{5}{0}{2}
\tableaubox{5}{1}{2}
\tableaubox{5}{2}{1}
\tableaulabel{5}{3}{}
\tableaubox{5}{4}{1}
\tableaubox{5}{5}{1}
\tableaubox{5}{6}{1}
\tableaulabel{6}{0}{}
\tableaubox{6}{1}{2}
\tableaubox{6}{2}{1}
\tableaulabel{6}{3}{}
\tableaubox{6}{4}{1}
\tableaubox{6}{5}{1}
\tableaubox{6}{6}{1}

      \tableaulabel{3}{7}{$\mapsto$}
      \tableaulabel{3}{8}{$\cdots$}
      \tableaulabel{3}{9}{$\mapsto$}
    
\tableaulabel{5.5}{3}{$\cdots$}
\tableaulabel{3}{5.5}{$\vdots$}

\tableaudraw[line width=1mm]{(4.5, -0.5) -- (4.5, 4.5) -- (-0.5, 4.5)}
      \end{scope}
    
      \begin{scope}[shift={(9,0)}]
        \tableaulabel{0}{0}{}
\tableaulabel{0}{1}{}
\tableaulabel{0}{2}{}
\tableaulabel{0}{3}{}
\tableaulabel{0}{4}{}
\tableaubox{0}{5}{3}
\tableaulabel{0}{6}{}
\tableaulabel{1}{0}{}
\tableaulabel{1}{1}{}
\tableaulabel{1}{2}{}
\tableaulabel{1}{3}{}
\tableaulabel{1}{4}{}
\tableaubox{1}{5}{3}
\tableaubox{1}{6}{1}
\tableaulabel{2}{0}{}
\tableaulabel{2}{1}{}
\tableaulabel{2}{2}{}
\tableaulabel{2}{3}{}
\tableaulabel{2}{4}{}
\tableaubox{2}{5}{3}
\tableaubox{2}{6}{1}
\tableaulabel{3}{0}{}
\tableaulabel{3}{1}{}
\tableaulabel{3}{2}{}
\tableaulabel{3}{3}{}
\tableaulabel{3}{4}{}
\tableaulabel{3}{5}{}
\tableaulabel{3}{6}{}
\tableaulabel{4}{0}{}
\tableaulabel{4}{1}{}
\tableaulabel{4}{2}{}
\tableaulabel{4}{3}{}
\tableaulabel{4}{4}{}
\tableaubox{4}{5}{3}
\tableaubox{4}{6}{1}
\tableaubox{5}{0}{2}
\tableaubox{5}{1}{2}
\tableaubox{5}{2}{1}
\tableaulabel{5}{3}{}
\tableaubox{5}{4}{1}
\tableaubox{5}{5}{1}
\tableaubox{5}{6}{1}
\tableaulabel{6}{0}{}
\tableaubox{6}{1}{2}
\tableaubox{6}{2}{1}
\tableaulabel{6}{3}{}
\tableaubox{6}{4}{1}
\tableaubox{6}{5}{1}
\tableaubox{6}{6}{1}

\tableaulabel{5.5}{3}{$\cdots$}
\tableaulabel{3}{5.5}{$\vdots$}

\tableaudraw[line width=1mm]{(4.5, -0.5) -- (4.5, 4.5) -- (-0.5, 4.5)}
      \end{scope}
    
  \end{tikzpicture} \caption{Commuting the operator $\Gamma_-\left(q^{1/2}\right)$ past every other $\Gamma_-$ has the effect of toggling every diagonal on the right side of the diagram, with the exception of the topmost and bottommost.}
  \label{fig:palindromicToggle}
\end{figure}
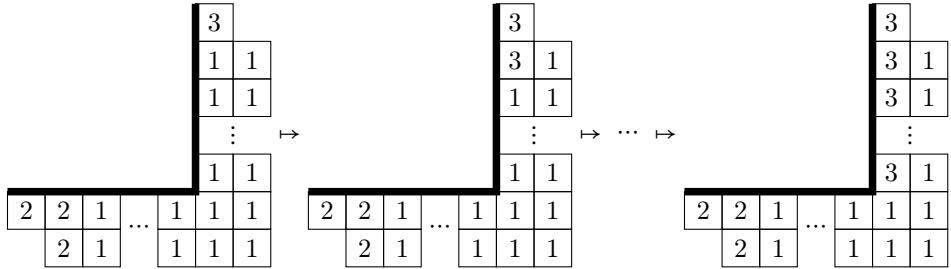

\begin{figure}
  \centering
  \begin{tikzpicture}
    
      \begin{scope}[shift={(0,0)}]
        \tableaulabel{0}{0}{}
\tableaulabel{0}{1}{}
\tableaulabel{0}{2}{}
\tableaulabel{0}{3}{}
\tableaulabel{0}{4}{}
\tableaubox{0}{5}{3}
\tableaulabel{0}{6}{}
\tableaulabel{1}{0}{}
\tableaulabel{1}{1}{}
\tableaulabel{1}{2}{}
\tableaulabel{1}{3}{}
\tableaulabel{1}{4}{}
\tableaubox{1}{5}{3}
\tableaubox{1}{6}{1}
\tableaulabel{2}{0}{}
\tableaulabel{2}{1}{}
\tableaulabel{2}{2}{}
\tableaulabel{2}{3}{}
\tableaulabel{2}{4}{}
\tableaubox{2}{5}{3}
\tableaubox{2}{6}{1}
\tableaulabel{3}{0}{}
\tableaulabel{3}{1}{}
\tableaulabel{3}{2}{}
\tableaulabel{3}{3}{}
\tableaulabel{3}{4}{}
\tableaulabel{3}{5}{}
\tableaulabel{3}{6}{}
\tableaulabel{4}{0}{}
\tableaulabel{4}{1}{}
\tableaulabel{4}{2}{}
\tableaulabel{4}{3}{}
\tableaulabel{4}{4}{}
\tableaubox{4}{5}{3}
\tableaubox{4}{6}{1}
\tableaubox{5}{0}{2}
\tableaubox{5}{1}{2}
\tableaubox{5}{2}{2}
\tableaulabel{5}{3}{}
\tableaubox{5}{4}{2}
\tableaubox{5}{5}{1}
\tableaubox{5}{6}{1}
\tableaulabel{6}{0}{}
\tableaubox{6}{1}{2}
\tableaubox{6}{2}{2}
\tableaulabel{6}{3}{}
\tableaubox{6}{4}{2}
\tableaubox{6}{5}{1}
\tableaubox{6}{6}{1}

\tableaulabel{5.5}{3}{$\cdots$}
\tableaulabel{3}{5.5}{$\vdots$}

\tableaudraw[line width=1mm]{(4.5, -0.5) -- (4.5, 4.5) -- (-0.5, 4.5)}
      \end{scope}
    
  \end{tikzpicture} \caption{The result of palindromically commuting the $\Gamma$ operators of the same sign.}
  \label{fig:palindromicToggleResult}
\end{figure}
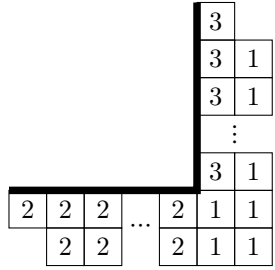

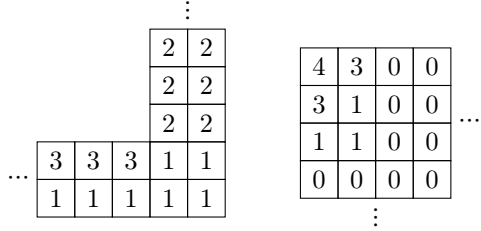
\begin{figure}
  \centering
  \begin{tikzpicture}
    
      \begin{scope}[shift={(0,0)}]
        \tableaulabel{0}{0}{}
\tableaulabel{0}{1}{}
\tableaulabel{0}{2}{}
\tableaubox{0}{3}{2}
\tableaubox{0}{4}{2}
\tableaulabel{1}{0}{}
\tableaulabel{1}{1}{}
\tableaulabel{1}{2}{}
\tableaubox{1}{3}{2}
\tableaubox{1}{4}{2}
\tableaulabel{2}{0}{}
\tableaulabel{2}{1}{}
\tableaulabel{2}{2}{}
\tableaubox{2}{3}{2}
\tableaubox{2}{4}{2}
\tableaubox{3}{0}{3}
\tableaubox{3}{1}{3}
\tableaubox{3}{2}{3}
\tableaubox{3}{3}{1}
\tableaubox{3}{4}{1}
\tableaubox{4}{0}{1}
\tableaubox{4}{1}{1}
\tableaubox{4}{2}{1}
\tableaubox{4}{3}{1}
\tableaubox{4}{4}{1}

\tableaulabel{3.5}{-1}{$\cdots$}
\tableaulabel{-1}{3.5}{$\vdots$}
      \end{scope}
    
      \begin{scope}[shift={(3.5,-0.25)}]
        \tableaubox{0}{0}{4}
\tableaubox{0}{1}{3}
\tableaubox{0}{2}{0}
\tableaubox{0}{3}{0}
\tableaubox{1}{0}{3}
\tableaubox{1}{1}{1}
\tableaubox{1}{2}{0}
\tableaubox{1}{3}{0}
\tableaubox{2}{0}{1}
\tableaubox{2}{1}{1}
\tableaubox{2}{2}{0}
\tableaubox{2}{3}{0}
\tableaubox{3}{0}{0}
\tableaubox{3}{1}{0}
\tableaubox{3}{2}{0}
\tableaubox{3}{3}{0}
\tableaulabel{1.5}{4}{$\cdots$}
\tableaulabel{4}{1.5}{$\vdots$}

      \end{scope}
    
  \end{tikzpicture} \caption{The two-leg RPP $\rho$ (left; after transposing) and the plane partition $\pi$ (right; after untoggling) corresponding to the SPP $\sigma$ from \figref{fig:twoLegSPPExample}.}
  \label{fig:twoLegBijectionOutput}
\end{figure}

\begin{Ex}
  Let $\lambda = (2, 2)$ and $\mu = (3, 1)$, and let $\sigma$ be the SPP of shape $(\lambda, \mu, \emptyset)$ and weight $16$ given in \figref{fig:twoLegSPPExample}. To associate $\sigma$ with a pair $(\rho, \pi)$ of an RPP $\rho$ of shape $(\lambda, \mu, \emptyset)$ and a plane partition $\pi$, we begin as in the one-leg case by iteratively toggling the diagonals of $\pi$ that begin with corners. In \figref{fig:twoLegSPPToggle}, we perform the toggles to produce squares instead of in lexicographic order --- the limiting behavior is easier to observe with this approach.

  In this example, the value of $N$ guaranteed by \propref{prop:twoLegLimitBoundary} is $N = 3$, as demonstrated by the final two objects in the sequence in \figref{fig:twoLegSPPToggle}. All future toggles pop off zeros, and the diagonals are constantly $(1, 1)$ sufficiently close to the main diagonal. The limiting diagram resulting from performing all of the toggles is then given by \figref{fig:twoLegSPPToggleResult}. We draw the bounding partitions $\lambda$ and $\mu$ as diagonals in accordance with the sequence of vertex operators --- the diagram begins at a diagonal of $(2, 2)$ and ends at a diagonal of $(3, 1)$. We can easily check at this halfway point that the weights are correct: the tableau in the top-left is weighted by hook length, so it has weight $13$. On the other hand, the RPP-like object left over from the toggling is counted by the generating function \eqref{eq:twoLegRppVertexOperators},
  meaning its weight is $\frac{1}{2} \cdot 2 + \frac{1}{2} + \frac{3}{2} = 3$, as expected.

  To complete the bijection, we palindromically toggle the remaining diagonals of the RPP-like object, as described in the proof of \thmref{thm:twoLegDecomposition}. Beginning with the right side of the vertex operator expression, we commute the $\Gamma_-\left(q^{1/2}\right)$ to the left past every other $\Gamma_-$, but no $\Gamma_+$. This corresponds to toggling every diagonal on the right of the diagram from top to bottom, except the topmost and bottommost, as shown in \figref{fig:palindromicToggle}. We next commute the $\Gamma_-\left(q^{3/2}\right)$ left past every $\Gamma_-$ except the $\Gamma_-\left(q^{1/2}\right)$, in effect toggling the diagonals from top to bottom but ignoring the second-to-bottom diagonal as well as the bottom one. This has no effect on our example, and in fact the rest of the $\Gamma_-$ commutations proceed without changing the diagram further. Commuting the $\Gamma_+$ is an analogous task, and the final diagram is given in \figref{fig:palindromicToggleResult}.

  Comparing this object to the minimal configuration, exactly two boxes have been removed, and the weight of the base RPP of shape $(\lambda, \mu, \emptyset)$ is 1, so this RPP is still weight 3. All that remains is to transpose the RPP and to convert the tableau produced in the first step back into a plane partition via \thmref{thm:macMachonsFunctionByToggling}, and we have successfully mapped the SPP $\sigma$ to an RPP $\rho$ of the same shape and a plane partition $\pi$ (\figref{fig:twoLegBijectionOutput}).
\end{Ex}

\section{Discussion and Future Directions} \label{sec:futureDirections}

\begin{figure}
  \centering
  \begin{tabu}{X[m,c]X[m,c]}
    \begin{tikzpicture}
      
      \begin{scope}[shift={(0,0)}]
        \tableaulabel{0}{0}{}
\tableaulabel{0}{1}{}
\tableaulabel{0}{2}{}
\tableaulabel{0}{3}{}
\tableaubox{0}{4}{5}
\tableaubox{0}{5}{5}
\tableaubox{0}{6}{3}
\tableaubox{0}{7}{3}
\tableaulabel{1}{0}{}
\tableaulabel{1}{1}{}
\tableaubox{1}{2}{6}
\tableaubox{1}{3}{4}
\tableaubox{1}{4}{3}
\tableaubox{1}{5}{2}
\tableaubox{1}{6}{2}
\tableaubox{1}{7}{2}
\tableaulabel{2}{0}{}
\tableaubox{2}{1}{3}
\tableaubox{2}{2}{2}
\tableaubox{2}{3}{2}
\tableaubox{2}{4}{1}
\tableaubox{3}{0}{4}
\tableaubox{3}{1}{3}
\tableaubox{3}{2}{2}
\tableaubox{3}{3}{1}
\tableaubox{3}{4}{1}
\tableaubox{4}{0}{3}
\tableaubox{4}{1}{1}
\tableaubox{4}{2}{1}
\tableaubox{4}{3}{1}
\tableaubox{5}{0}{2}
\tableaubox{5}{1}{1}
\tableaubox{5}{2}{1}
\tableaubox{6}{0}{2}
\tableaubox{6}{1}{1}
\tableaubox{6}{2}{1}
\tableaubox{7}{0}{2}
\tableaubox{7}{1}{1}
\tableaubox{7}{2}{1}

\tableaulabel{0.5}{8}{$\cdots$}
\tableaulabel{8}{1}{$\vdots$}
      \end{scope}
    
    \end{tikzpicture} & \includegraphics[height=2.5in]{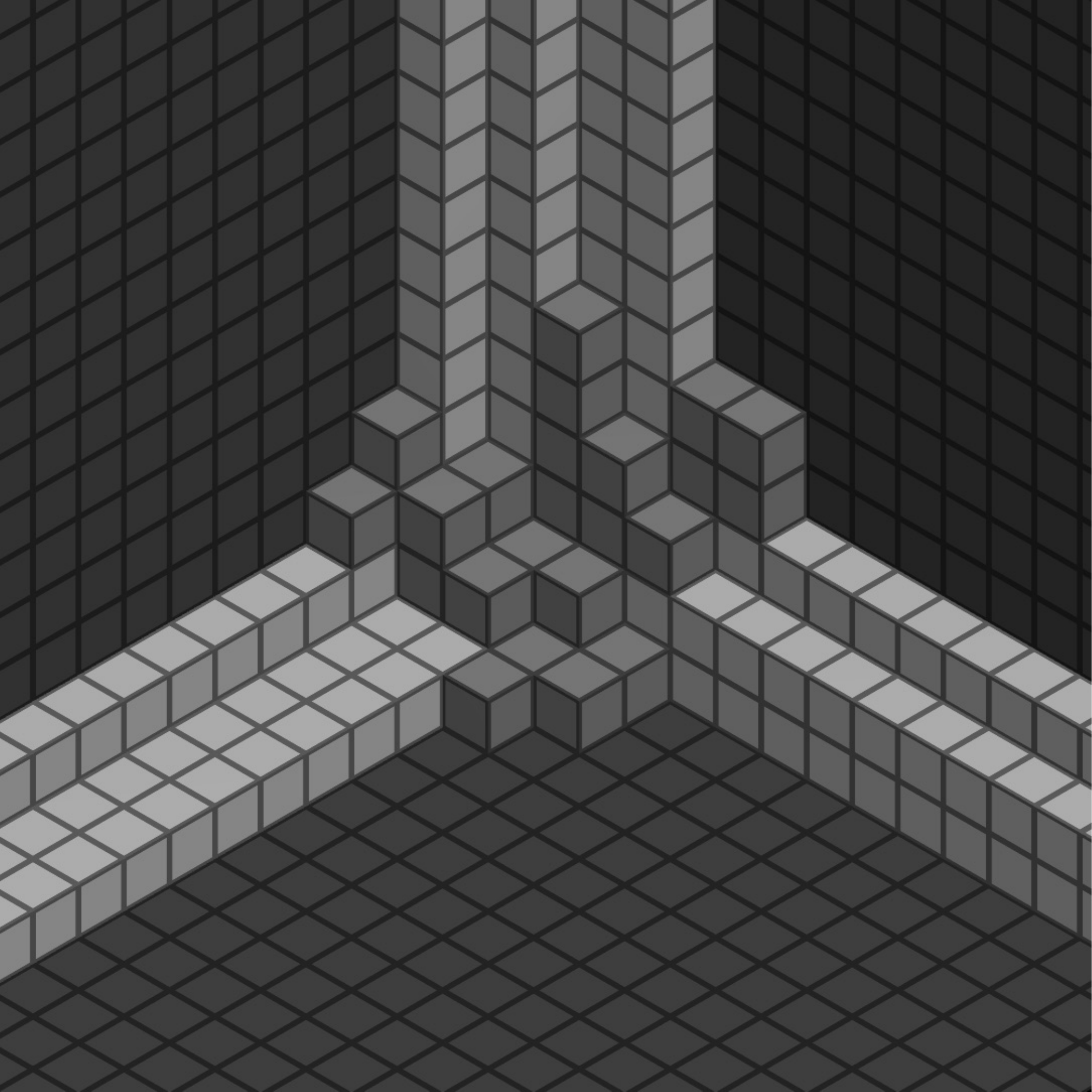}
  \end{tabu}
  \caption{A three-leg SPP of shape $((2, 1, 1), (3, 2), (4, 2, 1))$ and weight $\frac{17}{2}$, visualized as a grid of numbers (left) and a stack of 25 weight-contributing boxes (right; the minimal configuration has weight $-\frac{33}{2}$).}
  \label{fig:threeLegSpp}
\end{figure}

\begin{figure}
  \begin{tabu}{X[m,c]X[m,c]X[m,c]}
    \includegraphics[height=2in]{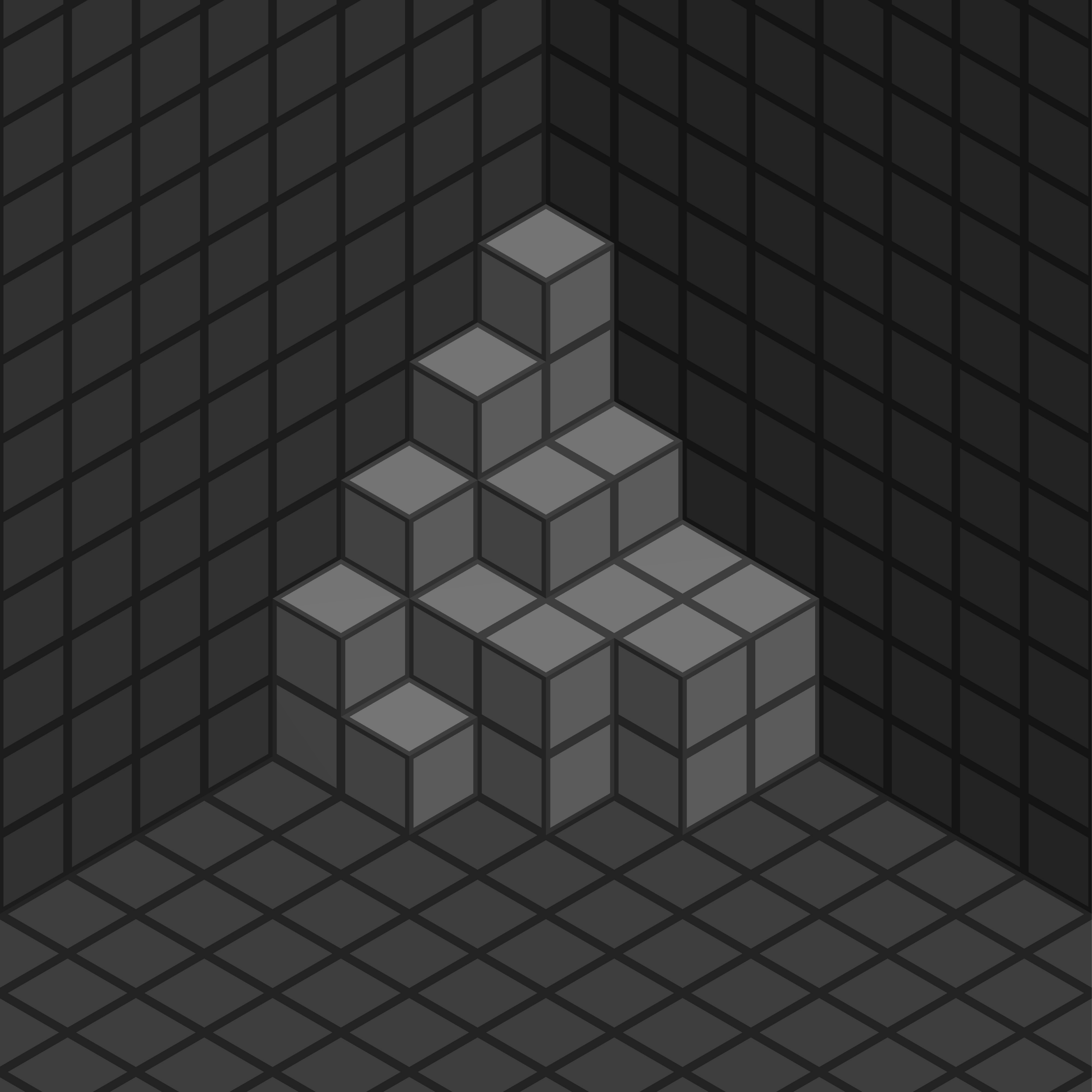} & \includegraphics[height=2in]{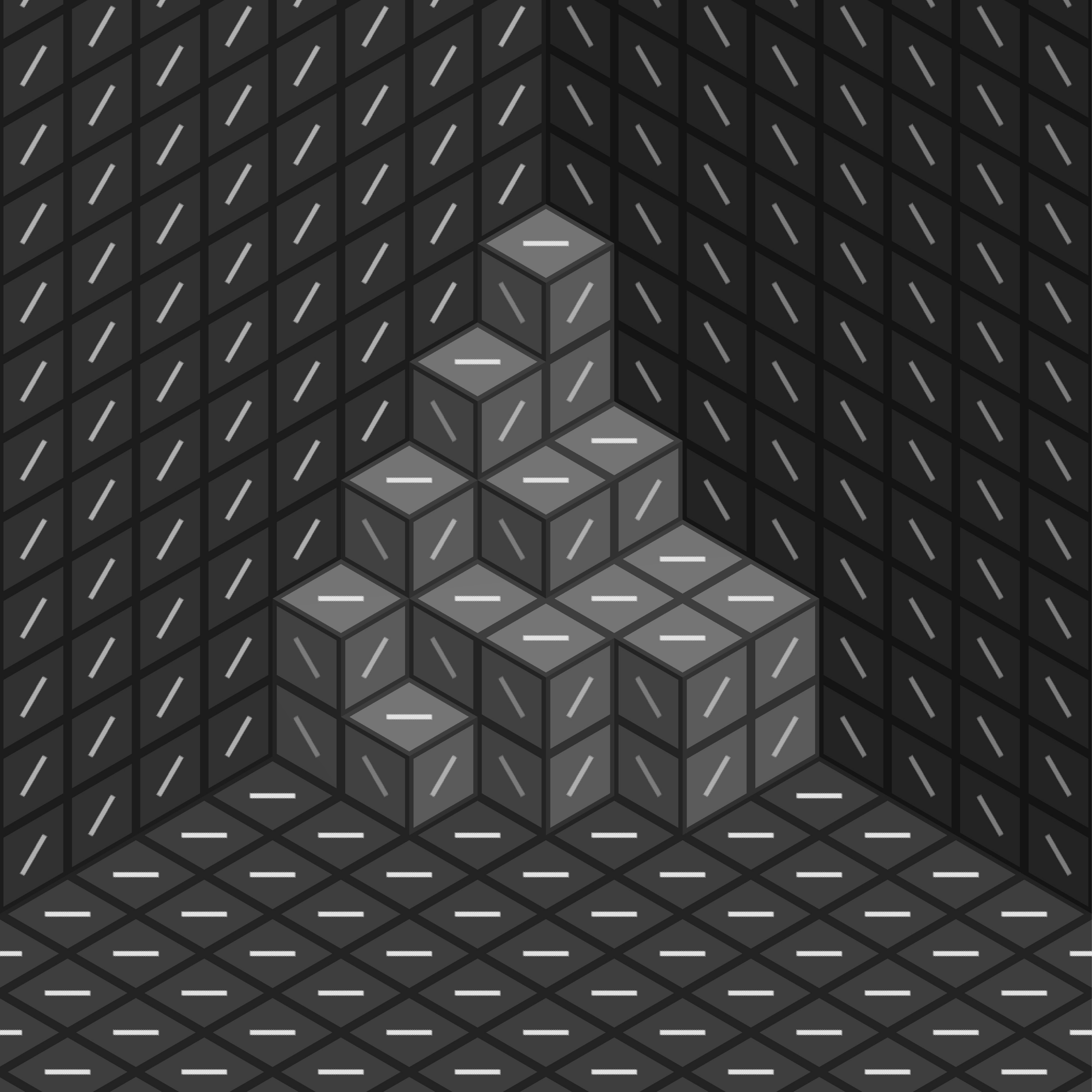} & \includegraphics[height=2in]{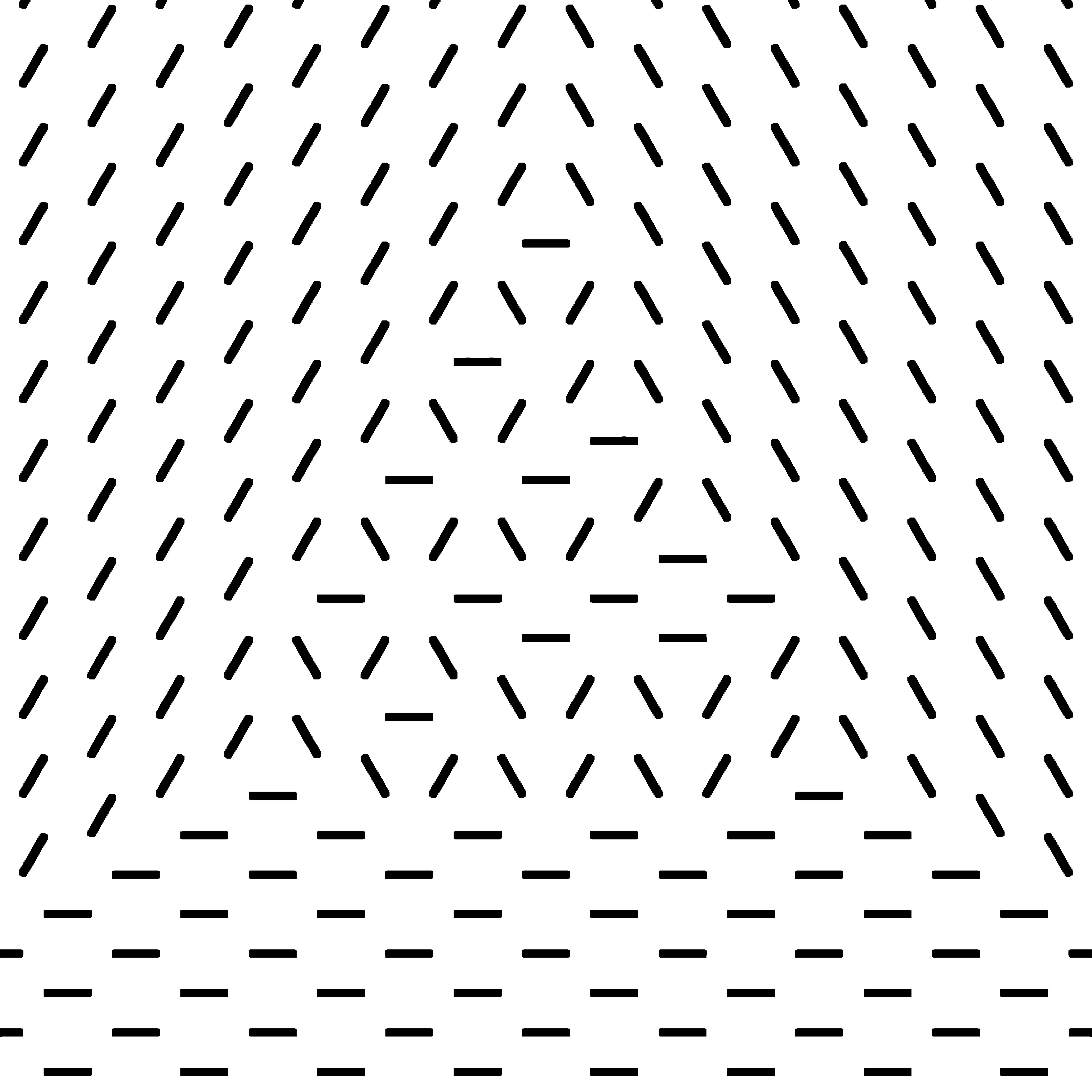}
  \end{tabu} \caption{A plane partition with walls shown (left), each rhombus painted with a dimer (center), and the configuration formed from those dimers (right).}
  \label{fig:folkloreBijection}
\end{figure}

\begin{figure}
  \centering
  \begin{tabu}{X[m,c]X[m,c]X[m,c]}
    \includegraphics[height=2in]{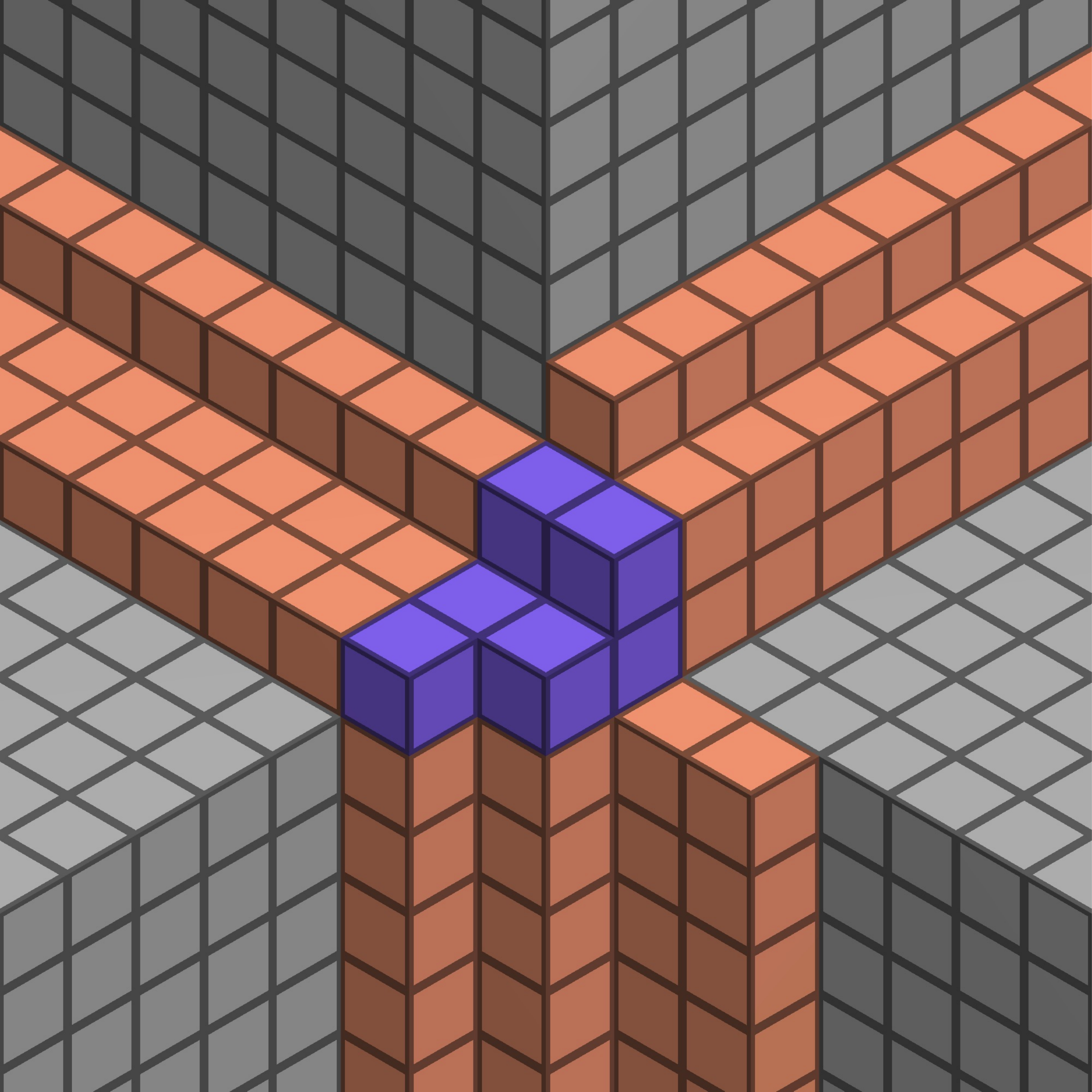} & \includegraphics[height=2in]{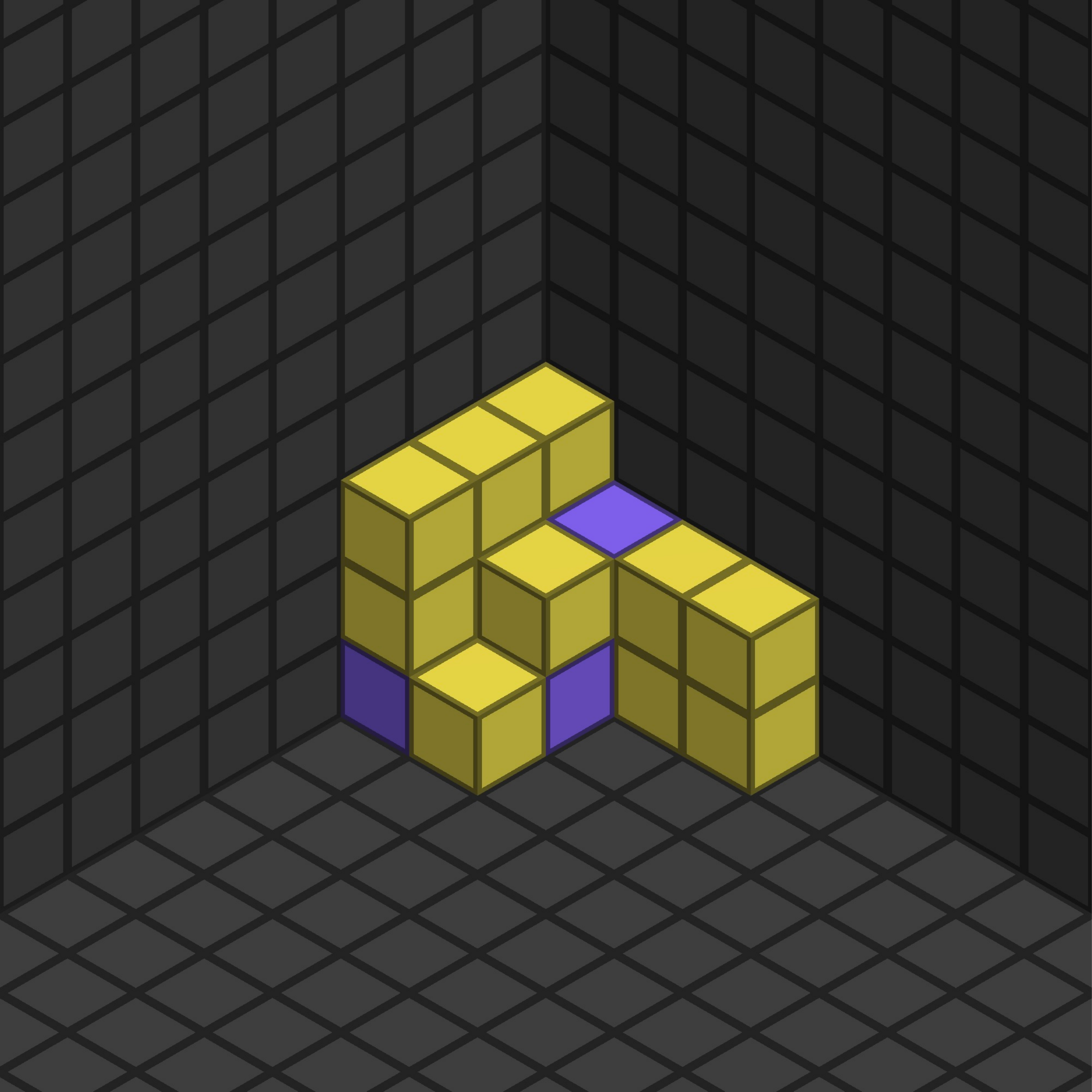} & \includegraphics[height=2in]{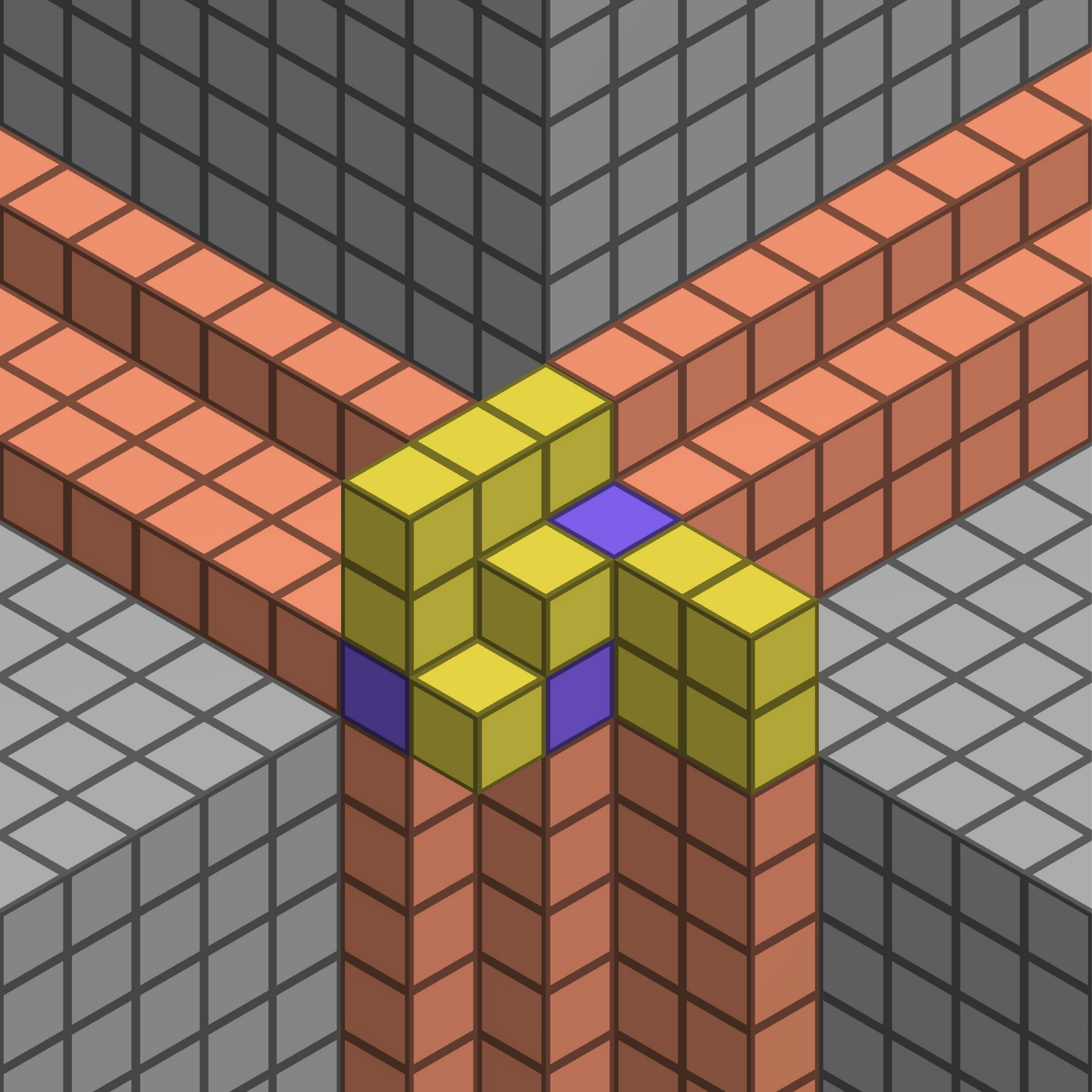}
  \end{tabu} \caption{The minimal configuration for an RPP of shape $((3, 2), (2, 1, 1), (4, 2, 1))$. The left and middle diagrams are superimposed to create the rightmost diagram; the purple boxes overlap and are all present with multiplicity two.}
  \label{fig:regionsIandIIandIII}
\end{figure}

\begin{figure}
  \centering
  \begin{center}
    \includegraphics[height=2in]{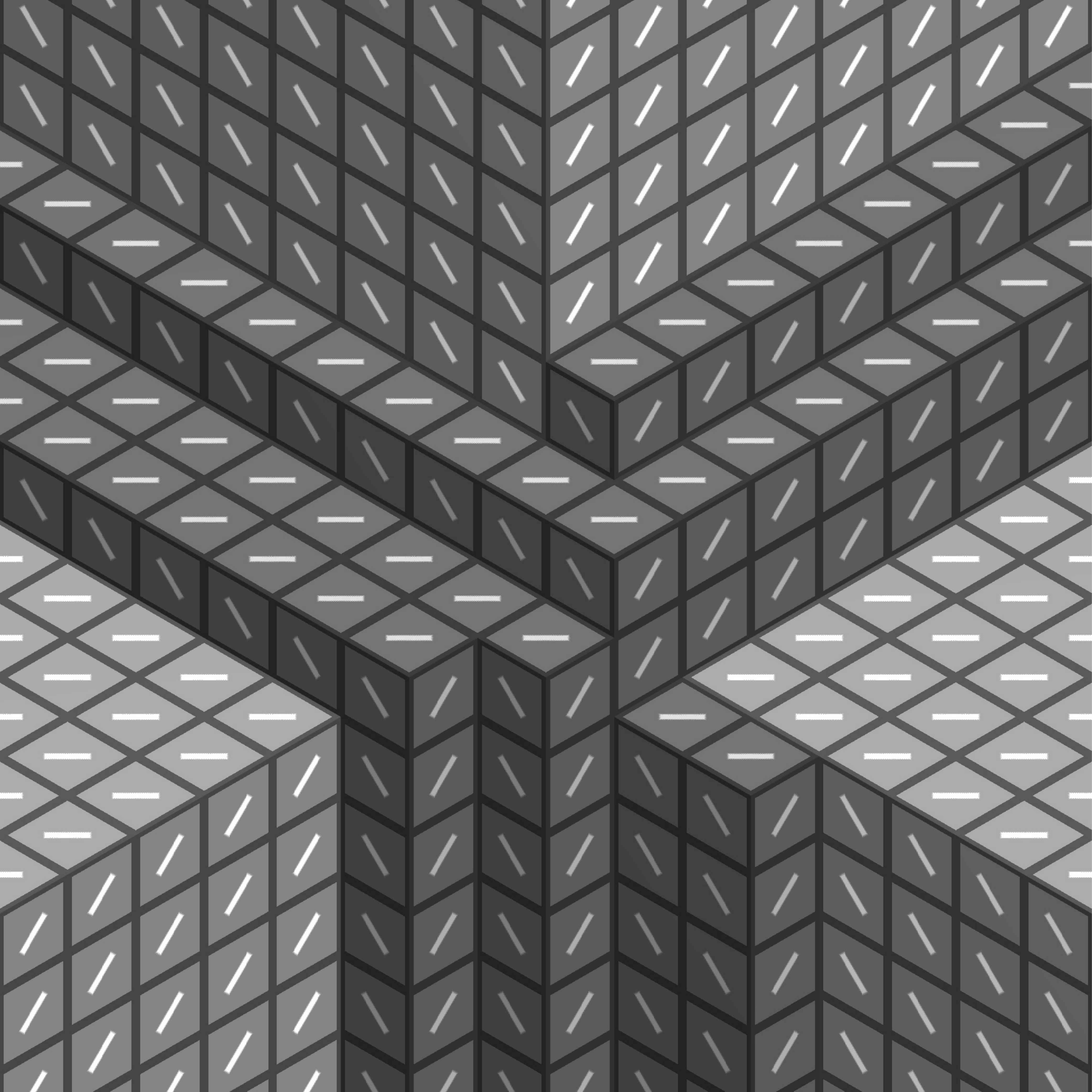} \hspace{4pt} \includegraphics[height=2in]{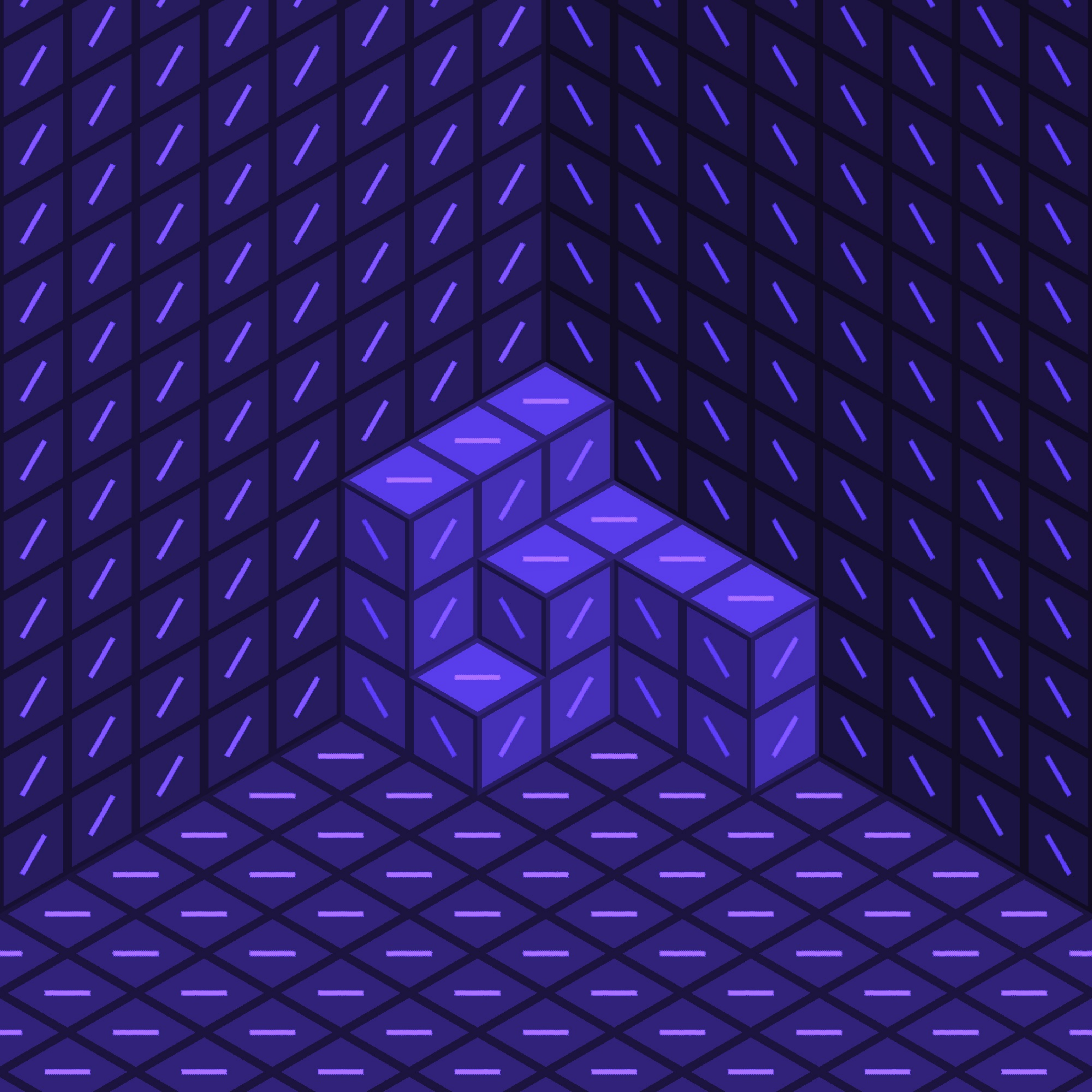}
  \end{center}
  \begin{center}
    \includegraphics[height=2in]{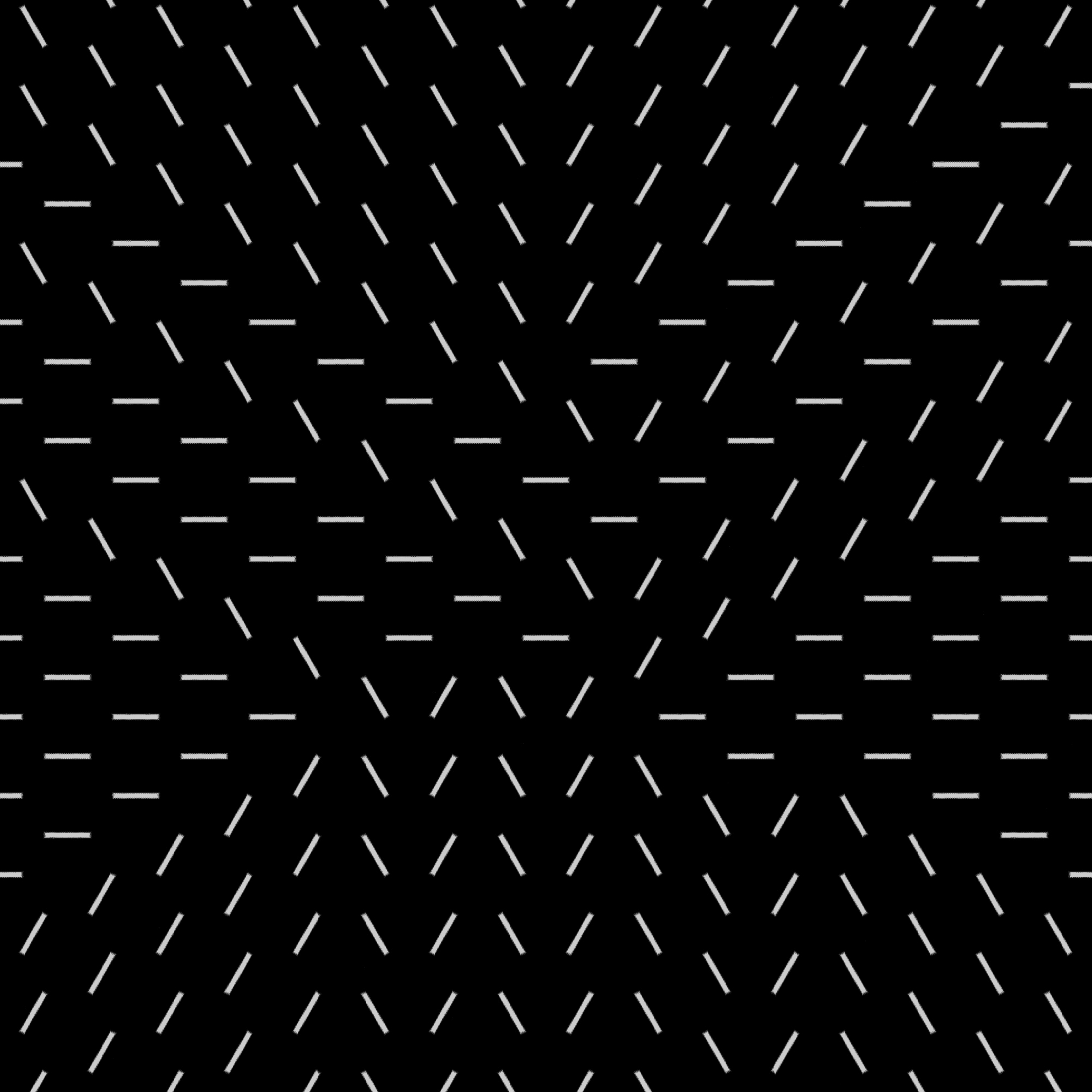} \hspace{4pt} \includegraphics[height=2in]{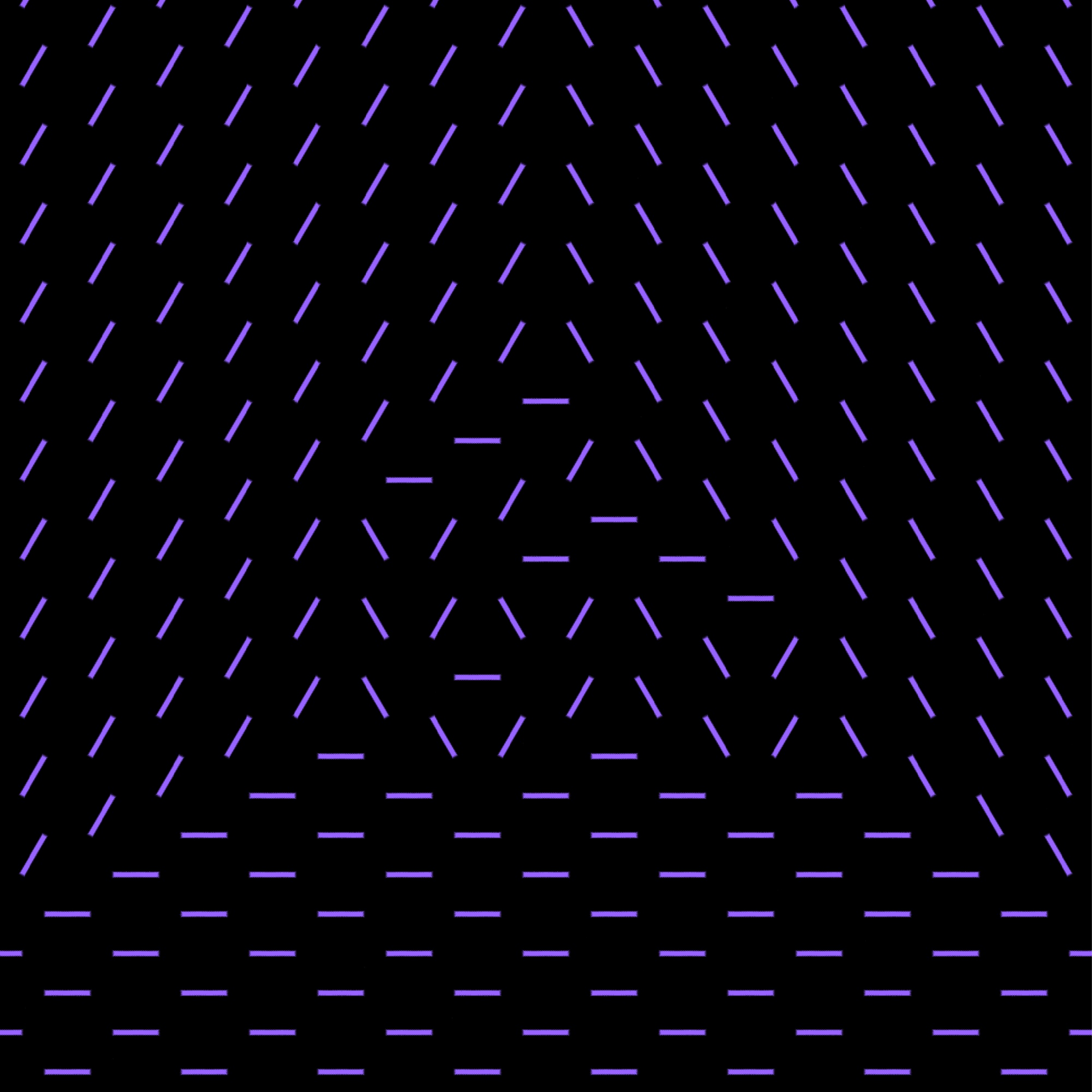}
  \end{center}
  \begin{center}
    \includegraphics[height=3in]{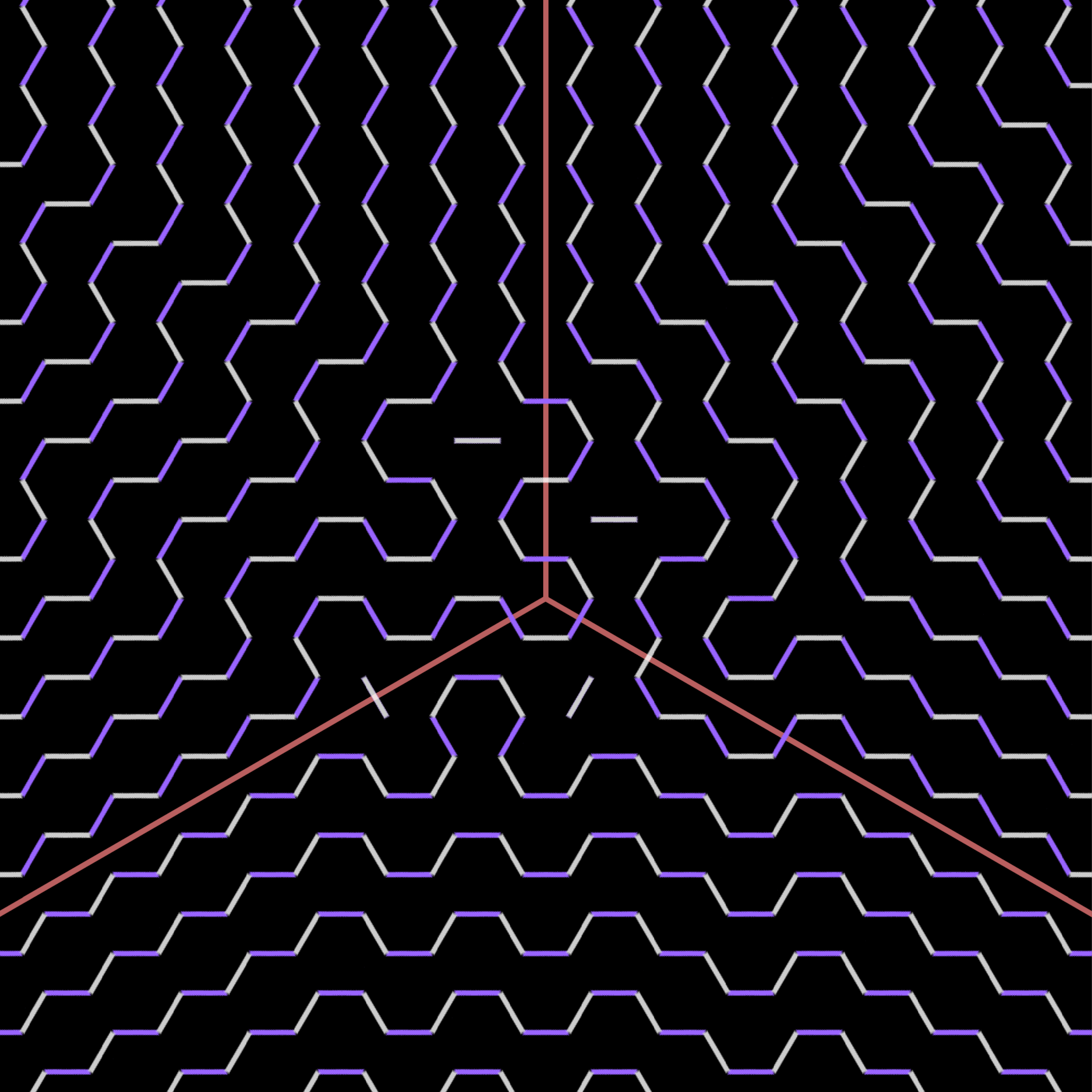}
  \end{center} \caption{By drawing dimers on the faces of the objects in \figref{fig:regionsIandIIandIII}, we produce the corresponding double-dimer configuration.}
  \label{fig:abConfigsAndDoubleDimers}
\end{figure}

The fully general case of the PT--DT correspondence involves objects that are substantively more complicated than the two special cases we have bijectivized. While three-leg SPPs are a direct generalization of one- and two-leg ones, as in \figref{fig:threeLegSpp} \cite{vertexOperators}, three-leg RPPs are a significant departure. Plane partitions are famously in correspondence with perfect matchings (or \textbf{dimer configurations}) on a hexagon lattice in $\mathbb{R}^2$ by drawing a line in the middle of each rhombus face when the 3-dimensional block expression of a plane partition is viewed from an isometric perspective, as in \figref{fig:folkloreBijection}. We refer the interested reader to \cite{dimerModelIntro} for further details on dimer configurations. Similar bijections hold for all SPPs and one- and two-leg RPPs. Since each vertex is matched with exactly one other, we call these \textbf{single-dimer} objects. In contrast, three-leg RPPs are \textbf{double-dimer} objects \cite{kenyonWilson, jenneWebbYoung, jenne}, meaning every vertex is matched to exactly two others. While a definition in terms of double-dimer objects is certainly more natural, a description in terms of boxes is much more conducive to our existing results. Such a description was introduced in \cite{jenneWebbYoung}, but it is nuanced and technical, involving boxes that can be present with multiplicity two and a labeling condition that permits only certain adjacency relationships between boxes in various regions. We sketch the minimal configuration for a three-leg RPP in \figref{fig:regionsIandIIandIII}, coloring the various regions, and depict the corresponding double-dimer object in \figref{fig:abConfigsAndDoubleDimers}.

Our methods of proving \thmref{thm:oneLegDecomposition} and \thmref{thm:twoLegDecomposition} only partially extend to three-leg objects. The generating function for three-leg SPPs \cite{vertexOperators} is given by
\begin{align}
  \label{eq:threeLegAppVertexOperators}
  V_{(\lambda, \mu, \nu)}(q) = \leftbra{\lambda}\,\prod_{n \in \mathbb{Z}} \Gamma_{e(n)}\left(q^{p(n)}\right)\,\rightket{\mu}
\end{align}
for the edge sign and power sequences $e(n) = e_\nu(n)$ and $p(n) = p_\nu(n)$; i.e.\ a combination of \eqref{eq:oneLegSppVertexOperators} and \eqref{eq:twoLegSppVertexOperators}. Iteratively toggling the diagonals of an SPP of shape $(\lambda, \mu, \nu)$ then produces two objects: a hook-length-weighted tableau that we may decompose and untoggle into a plane partition and a one-leg RPP of shape $(\emptyset, \emptyset, \nu)$, and an object that resembles a two-leg RPP of shape $(\lambda, \mu, \emptyset)$, but whose weight differs due to the edge sign sequence $e_\nu(n)$.

What remains is to define a bijection from three-leg RPPs to pairs of one-leg RPPs and these two-leg-RPP-esque objects. However, we are unaware of a vertex-operator description of three-leg RPPs --- the global labeling conditions make such operators difficult to define. In a paper in preparation \cite{godarYoung2}, we prove results that characterize the poset of valid entries for certain individual cells in a three-leg RPP when all others remain constant, and we conjecture a minor generalization that would enable us to define such operators. However, our results suggest that it may be difficult to define a general involution on the poset of valid entries for a three-leg RPP, precluding a naive generalization of the toggle map. We plan to define these vertex operators, find and bijectivize their commutation relations, and define a more general toggle using them, with the goal of defining the bijective decomposition of three-leg RPPs that we seek.

\printbibliography

@article{pak,
      author = {Igor Pak},
      title = {Hook Length Formula and Geometric Combinatorics},
      journal = {Séminaire Lotharingien de Combinatoire},
      volume = {46},
      year = {2001},
      URL = {https://www.mat.univie.ac.at/~slc/wpapers/s46pak.html},
      eprint = {https://www.mat.univie.ac.at/~slc/wpapers/s46pak.html}
}

@misc{hopkinsNotes,
      author = {Samuel Hopkins},
      title = {RSK via local transformations},
      howpublished = {\url{https://www.samuelfhopkins.com/docs/rsk.pdf}},
      note = {Accessed 2024-06-22},
      year = {2014}
}

@misc{jenneWebbYoung,
      title={Double-dimer condensation and the PT-DT correspondence}, 
      author={Helen Jenne and Gautam Webb and Benjamin Young},
      year={2022},
      eprint={2109.11773},
      archivePrefix={arXiv},
      primaryClass={math.CO},
      url={https://arxiv.org/abs/2109.11773}, 
}

@article{jenne,
  title={Combinatorics of the double-dimer model},
  author={Helen Jenne},
  journal={Advances in Mathematics},
  year={2019},
  url={https://api.semanticscholar.org/CorpusID:207853447}
}

@article{kenyonWilson,
  title={Boundary Partitions in Trees and Dimers},
  author={Richard W. Kenyon and David Bruce Wilson},
  journal={Transactions of the American Mathematical Society},
  year={2006},
  volume={363},
  pages={1325-1364},
  url={https://api.semanticscholar.org/CorpusID:18950209}
}

@article{dimerModelIntro,
  title={An introduction to the dimer model},
  author={Richard W. Kenyon},
  journal={arXiv: Combinatorics},
  year={2003},
  url={https://api.semanticscholar.org/CorpusID:3083216}
}

@article{PT,
      author = {Rahul Pandharipande and Richard P Thomas},
      title = {{The $3$–fold vertex via stable pairs}},
      volume = {13},
      journal = {Geometry \& Topology},
      number = {4},
      publisher = {MSP},
      pages = {1835 -- 1876},
      year = {2009},
      doi = {10.2140/gt.2009.13.1835},
      URL = {https://doi.org/10.2140/gt.2009.13.1835}
}

@inproceedings{formalPowerSeries,
  title={Enumerative Combinatorics: Volume 1},
  author={Richard P. Stanley},
  year={2011},
  url={https://api.semanticscholar.org/CorpusID:198489053}
}

@book{stanley,
      title     = "Enumerative Combinatorics",
      volume    = {2},
      author    = "Stanley, Richard",
      year      = {2023},
      publisher = "Cambridge University Press",
}

@article{okounkovReshetikhin,
      author = {Andrei Okounkov and Nicolai Reshetikhin},
      title = {{Random Skew Plane Partitions and the Pearcey Process}},
      journal = {Communications in Mathematical Physics},
      number = {269},
      pages = {571 -- 609},
      year = {2007},
      doi = {10.1007/s00220-006-0128-8},
      URL = {https://doi.org/10.1007/s00220-006-0128-8}
}

@Inbook{vertexOperators,
      author="Okounkov, Andrei
      and Reshetikhin, Nikolai
      and Vafa, Cumrun",
      editor="Etingof, Pavel
      and Retakh, Vladimir
      and Singer, I. M.",
      title="Quantum Calabi-Yau and Classical Crystals",
      bookTitle="The Unity of Mathematics: In Honor of the Ninetieth Birthday of I.M. Gelfand",
      year="2006",
      publisher="Birkh{\"a}user Boston",
      address="Boston, MA",
      pages="597--618",
      isbn="978-0-8176-4467-3",
      doi="10.1007/0-8176-4467-9_16",
      url="https://doi.org/10.1007/0-8176-4467-9_16"
}

@article{stanleyThesis,
  title={Theory and Application of Plane Partitions. Part 2},
  author={Richard P. Stanley},
  journal={Studies in Applied Mathematics},
  year={1971},
  volume={50},
  pages={259-279},
  url={https://api.semanticscholar.org/CorpusID:126649639}
}

@article{hillmanGrassl,
  title={Reverse Plane Partitions and Tableau Hook Numbers},
  author={A. P. Hillman and R. M. Grassl},
  journal={J. Comb. Theory A},
  year={1976},
  volume={21},
  pages={216-221},
  url={https://api.semanticscholar.org/CorpusID:30944428}
}

@article{sulzgruber,
  title={Inserting rim-hooks into reverse plane partitions},
  author={Robin Sulzgruber},
  journal={Journal of Combinatorics},
  year={2017},
  url={https://api.semanticscholar.org/CorpusID:119621454}
}

@book{kac, place={Cambridge}, edition={3}, title={Infinite-Dimensional Lie Algebras}, publisher={Cambridge University Press}, author={Kac, Victor G.}, year={1990}}

@article{steepTilings,
  title={From Aztec diamonds to pyramids: steep tilings},
  author={J{\'e}r{\'e}mie Bouttier and Guillaume Chapuy and Sylvie Corteel},
  journal={arXiv: Combinatorics},
  year={2014},
  url={https://api.semanticscholar.org/CorpusID:119164929}
}

@inproceedings{bressoud,
  title={Proofs and Confirmations: The Story of the Alternating-Sign Matrix Conjecture},
  author={David M. Bressoud},
  year={1999},
  url={https://api.semanticscholar.org/CorpusID:51454856}
}

@book{andrews,
      place={Cambridge},
      series={Encyclopedia of Mathematics and its Applications},
      title={The Theory of Partitions},
      publisher={Cambridge University Press},
      author={Andrews, George E.},
      year={1984},
      collection={Encyclopedia of Mathematics and its Applications}
}

@article{okounkovReshetikhin2,
  title={Correlation function of Schur process with application to local geometry of a random 3-dimensional Young diagram},
  author={Andrei Okounkov and Nikolai Reshetikhin},
  journal={Journal of the American Mathematical Society},
  year={2001},
  volume={16},
  pages={581-603},
  url={https://api.semanticscholar.org/CorpusID:11575202}
}

@article{randomPartitions1,
  title={Symmetric Functions and Random Partitions},
  author={Andrei Okounkov},
  journal={arXiv: Combinatorics},
  year={2003},
  pages={223-252},
  url={https://api.semanticscholar.org/CorpusID:15558365}
}

@misc{randomPartitions2,
      title={Infinite wedge and random partitions}, 
      author={Andrei Okounkov},
      year={2000},
      eprint={math/9907127},
      archivePrefix={arXiv},
      primaryClass={math.RT},
      url={https://arxiv.org/abs/math/9907127}, 
}

@inproceedings{okounkovReshetikhinVafa,
  title={Quantum Calabi-Yau and Classical Crystals},
  author={Andrei Okounkov and Nikolai Reshetikhin and Cumrun Vafa},
  year={2003},
  url={https://api.semanticscholar.org/CorpusID:16711409}
}

@article{railYardGraphs,
  title={Dimers on Rail Yard Graphs},
  author={C{\'e}dric Boutillier and J{\'e}r{\'e}mie Bouttier and Guillaume Chapuy and Sylvie Corteel and Sanjay Ramassamy},
  journal={arXiv: Mathematical Physics},
  year={2015},
  url={https://api.semanticscholar.org/CorpusID:54520004}
}

@article{orbifoldTopologicalVertex,
  title={The Orbifold Topological Vertex},
  author={J. M. Bryan and Charles Cadman and Benjamin Young},
  journal={arXiv: Algebraic Geometry},
  year={2010},
  url={https://api.semanticscholar.org/CorpusID:50069951}
}

@article{quotientsWriteUp,
  title={Tilings of benzels via the abacus bijection},
  author={Colin Defant and Rupert Li and James Gary Propp and Benjamin Young},
  journal={Combinatorial Theory},
  year={2022},
  url={https://api.semanticscholar.org/CorpusID:252531638}
}

@article{MNOP,
  title={Gromov–Witten theory and Donaldson–Thomas theory, I},
  author={Davesh Maulik and Nikita Nekrasov and Andrei Okounkov and Rahul Pandharipande},
  journal={Compositio Mathematica},
  year={2003},
  volume={142},
  pages={1263 - 1285},
  url={https://api.semanticscholar.org/CorpusID:5760317}
}

@misc{godarYoung2,
  title={Entry Posets and Vertex Operators for labeled $AB$ configurations},
  author={Cruz Godar and Benjamin Young},
  year={In preparation},
}

\end{document}